\documentclass[11pt,oneside]{amsart}
\usepackage{typesnew-arxiv}      

\def\IImIVExlabel{\redtype{IV\e{3-3}(3)II}}

\def\IImIVExeqn{y^3\!+\!y=\pi ^4x^5\!+\!\pi x^3}

\def\IImIVExfibre{\begin{tikzpicture}[xscale=1,yscale=0.9]
\draw[thick,shorten >=2] (-1/3,0)--(31/6,0) node[inner sep=2pt,blue,above left,scale=0.8] {3} node[inner sep=2pt,below left,scale=0.8] {$\Gamma_1$};
\draw[shorten <=-3pt] (0,0)--node[inner sep=2pt,scale=0.8,blue,right] {1} (0,3/5);
\draw[shorten <=-3pt] (4/5,0)--node[inner sep=2pt,scale=0.8,blue,right] {1} (4/5,3/5);
\draw[shorten <=-3pt] (8/5,0)--node[inner sep=2pt,scale=0.8,blue,right] {1} (8/5,3/5);
\draw[shorten <=-3pt, shorten >=-3pt] (10/3,0)--node[inner sep=2pt,scale=0.8,blue,right] {3} (10/3,3/5);
\draw[shorten <=-3pt, shorten >=-3pt] (10/3,3/5)--node[inner sep=2pt,scale=0.8,blue,above] {3} (41/15,3/5);
\draw[shorten <=-3pt, shorten >=-3pt] (41/15,3/5)--node[inner sep=2pt,scale=0.8,blue,left] {3} (41/15,6/5);
\draw[thick,shorten >=2] (12/5,6/5)--(179/30,6/5) node[inner sep=2pt,blue,above left,scale=0.8] {6} node[inner sep=2pt,below left,scale=0.8] {$\Gamma_2$};
\draw[shorten <=-3pt] (10/3,6/5)--node[inner sep=2pt,scale=0.8,blue,right] {1} (10/3,9/5);
\draw[shorten <=-3pt] (62/15,6/5)--node[inner sep=2pt,scale=0.8,blue,right] {2} (62/15,9/5);
\end{tikzpicture} }

\hypersetup{
    colorlinks=false,
    pdfborder={0 0 0},
    pdfborderstyle={/S/U/W 0},
}

\begin{document}

\title{A classification of reduction types of curves}

\author{Tim Dokchitser}
\address{Department of Mathematics, University Walk, Bristol BS8 1TW, UK}
\email{tim.dokchitser@bristol.ac.uk}

\begin{abstract}
The aim of this paper is to classify reduction types of algebraic curves. Reduction types capture the discrete
invariants of fibres in one-dimensional families of curves, and they have been described in genus 1, 2 and 3.
For fixed genus $g>1$, they form finitely many families, and we explain how to construct them,
and introduce a naming convention.
\end{abstract}

\maketitle
\tableofcontents
\newpage

\section{Introduction}
\label{sintro}

The aim of this paper is to classify reduction types for curves of arbitrary genus.
Suppose

\begin{itemize}
\item
$O$ is an excellent\footnote{For example, if $O$ is complete or if $\vchar K=0$, then $O$ is automatically excellent.}
discrete valuation ring with field of fractions $K$ and algebraically
closed residue field~$k$,
\item
$C/K$ is a complete, smooth, connected curve of genus $g(C)>0$,
\item
$\cC/O$ is a \textbf{model} of $C$, by which we mean a proper, flat scheme
with generic fibre $\cC_K\iso C$ and connected 1-dimensional special fibre $\cC_k$, not necessarily reduced
or irreducible.
\end{itemize}

\noindent
Among models of $C$ with $\cC$ regular and $\cC_k\subset\cC$
normal crossings divisor (\textbf{rnc} models), there exists a unique minimal one
(\textbf{mrnc} model), see Lipman \cite{Lip} or Saito \cite[3.2.2]{Saito}.

\smallskip

\begin{center}
\tikzset{ecpic/.pic={                 
  \path (3,0)  coordinate (p1);
  \path (1,3)  coordinate (p2);
  \path (0,4) coordinate (p3);
  \path (1,5) coordinate (p4);
  \path (3,8) coordinate (p5);
  \draw[thick]
    (p1) edge[out=100,in=40] (p2)
    (p2) edge[out=220,in=270] (p3)
    (p3) edge[out=90,in=140] (p4)
    (p4) edge[out=320,in=260] (p5);
}}
\tikzset{IIpic/.pic={
  \path (0,0) coordinate (p1);
  \path (-0.5,6) coordinate (p2);
  \draw[thick,-] (p2)--(p1) node[right,mult] {6};
  \draw[lin] ($(p1)!1/5!(p2)$) -- ++(-2,0) node[below,mult] {1};
  \draw[lin] ($(p1)!1/2!(p2)$) -- ++(-2,0) node[below,mult] {2};
  \draw[lin] ($(p1)!4/5!(p2)$) -- ++(-2,0) node[below,mult] {3};
}}
\begin{tikzpicture}[
    scale=0.3,
    mult/.style={blue,pos=0.94,inner sep=2pt,scale=0.8},
    lin/.style={draw,-,thick,shorten <=-5pt,shorten >=-5pt}]
  \path[rounded corners=4pt, top color=gray!5, bottom color=gray!20, draw=black!40]  
    (-0.75,1.3) to[out=5,in=180] (10,1.3)     
                to[out=85,in=-90] (9.5,9.7)   
                to[out=175,in=0] (-1.75,9.7)  
                to[out=-95,in=90] (-1.25,1.3) 
    -- cycle;
  \pic[rotate=3,scale=0.27] at (0,2) {ecpic};                  
  \pic[xscale=0.3,yscale=0.35] at (7.5,2.1) {IIpic};           
  \draw (-0.5,0)--                                             
    node [circle,scale=0.8,ball color=black!50,pos=0.2] {}        
    node [below=5pt,pos=0.2,scale=0.8] {$\Spec K$}
    node [circle,scale=0.4,ball color=black!50,pos=0.75] {}       
    node [below=5pt,pos=0.75,scale=0.8] {$\Spec k$} (10,0);
  \node[align=center,scale=0.8] (genlabel) at (-10.5,6.3) {Generic fibre\\ $\cC_K\iso_K C$};
  \draw[->,shorten <=1pt] (genlabel) to[out=-10,in=180] (-1,5.5);
  \node[align=center,scale=0.8,anchor=west] (splabel) at (17,5.5)
     {Special fibre $\cC_k/k$, a divisor\\
       with normal crossings: here, four\\genus 0 curves with multiplicities\\
      1,2,3 and 6 (Kodaira type II)};
  \draw[->,shorten <=1pt] (splabel) to[out=175,in=0] (8,4.5);
  \node (cC) at (11,9) {$\cC$};                                
  \node (cO) at (13,0) {$\Spec O$};                            
\end{tikzpicture}
\end{center}

The \textbf{reduction type} of $C/K$ encodes the discrete invariants of the special fibre of the mrnc model:
multiplicities of irreducible components (marked {\color{blue}1}, {\color{blue}2}, ...),
their geometric genera (marked {\color{blue}g1}, {\color{blue}g2}, ..., omitted if zero),
and the intersection number for every pair of components (see~\ref{defredtype}).

\subsection{Historical background}
\label{ssHistorical}

In 1960 Kodaira \cite{Ko} studied complex elliptic fibrations of relatively minimal surfaces over curves.
He classified possible configurations of singular fibres, now called Kodaira types, and
N\'eron~\cite{Ne2} extended this to elliptic curves over DVRs.
Stated in terms of mrnc models, the possible reduction types of elliptic curves are (cf.~\ref{exellmrnc}):

\label{kodairatypes}
\begin{center}
\begin{tabular}{c@{\qquad}c@{\qquad}c@{\qquad}c@{\qquad}c}
\begin{tikzpicture}[xscale=0.5,yscale=0.7] \draw[thick,shorten >=2] (-1/3,0)--(11/6,0) node[inner sep=2pt,blue,above left,scale=0.8] {1 \smash{g}1} node[inner sep=2pt,below left,scale=0.8] {\phantom{$\int^X$}\redtype{I_{g1}} (=\IZ)}; \end{tikzpicture}  & \begin{tikzpicture}[xscale=0.5,yscale=0.7] \draw[thick,shorten >=2] (-0.09,0)--(2.67,0) node[inner sep=2pt,blue,above left,scale=0.8] {1} node[inner sep=2pt,below left,scale=0.8] {\phantom{$\int^X$}\In{n}}; \draw[scale=0.8] (0,0) ++(0.05,0.3) node[anchor=east,blue,scale=0.7] {1} ++(0.3,-0.45)--++(115:0.75)++(295:0.15)++(245:0.15)--++(65:0.75)++(245:0.15)++(180:0.15)--++(0:0.22) ++(0:0.15) node{$\cdot$} ++(0:0.15) node{$\cdot$} node[below=0.1,anchor=north,scale=0.7]{$n$}++(0:0.15) node{$\cdot$} ++(0:0.15) --++(0:0.22)++(180:0.15)++(115:0.15)--++(295:0.75)++(115:0.15)++(65:0.15)--++(245:0.75)++(0.3,0.45) node[anchor=west,blue,scale=0.7] {1};\end{tikzpicture}  & \begin{tikzpicture}[xscale=0.5,yscale=0.7] \draw[thick,shorten >=2] (-1/3,0)--(103/30,0) node[inner sep=2pt,blue,above left,scale=0.8] {3} node[inner sep=2pt,below left,scale=0.8] {\phantom{$\int^X$}\redtype{IV}}; \draw[shorten <=-3pt] (0,0)--node[inner sep=2pt,scale=0.8,blue,right] {1} (0,3/5); \draw[shorten <=-3pt] (4/5,0)--node[inner sep=2pt,scale=0.8,blue,right] {1} (4/5,3/5); \draw[shorten <=-3pt] (8/5,0)--node[inner sep=2pt,scale=0.8,blue,right] {1} (8/5,3/5); \end{tikzpicture}  & \begin{tikzpicture}[xscale=0.5,yscale=0.7] \draw[thick,shorten >=2] (-1/3,0)--(103/30,0) node[inner sep=2pt,blue,above left,scale=0.8] {4} node[inner sep=2pt,below left,scale=0.8] {\phantom{$\int^X$}\redtype{III}}; \draw[shorten <=-3pt] (0,0)--node[inner sep=2pt,scale=0.8,blue,right] {1} (0,3/5); \draw[shorten <=-3pt] (4/5,0)--node[inner sep=2pt,scale=0.8,blue,right] {1} (4/5,3/5); \draw[shorten <=-3pt] (8/5,0)--node[inner sep=2pt,scale=0.8,blue,right] {2} (8/5,3/5); \end{tikzpicture}  & \begin{tikzpicture}[xscale=0.5,yscale=0.7] \draw[thick,shorten >=2] (-1/3,0)--(103/30,0) node[inner sep=2pt,blue,above left,scale=0.8] {6} node[inner sep=2pt,below left,scale=0.8] {\phantom{$\int^X$}\redtype{II}}; \draw[shorten <=-3pt] (0,0)--node[inner sep=2pt,scale=0.8,blue,right] {1} (0,3/5); \draw[shorten <=-3pt] (4/5,0)--node[inner sep=2pt,scale=0.8,blue,right] {2} (4/5,3/5); \draw[shorten <=-3pt] (8/5,0)--node[inner sep=2pt,scale=0.8,blue,right] {3} (8/5,3/5); \end{tikzpicture} \cr \begin{tikzpicture}[xscale=0.5,yscale=0.7] \draw[thick,shorten >=2] (-1/3,0)--(127/30,0) node[inner sep=2pt,blue,above left,scale=0.8] {2} node[inner sep=2pt,below left,scale=0.8] {\phantom{$\int^X$}\redtype{I^*_0}}; \draw[shorten <=-3pt] (0,0)--node[inner sep=2pt,scale=0.8,blue,right] {1} (0,3/5); \draw[shorten <=-3pt] (4/5,0)--node[inner sep=2pt,scale=0.8,blue,right] {1} (4/5,3/5); \draw[shorten <=-3pt] (8/5,0)--node[inner sep=2pt,scale=0.8,blue,right] {1} (8/5,3/5); \draw[shorten <=-3pt] (12/5,0)--node[inner sep=2pt,scale=0.8,blue,right] {1} (12/5,3/5); \end{tikzpicture}  & \begin{tikzpicture}[xscale=0.5,yscale=0.7] \draw[thick,shorten >=2] (-1/3,0)--(53/15,0) node[inner sep=2pt,blue,above left,scale=0.8] {2} node[inner sep=2pt,below left,scale=0.8] {\phantom{$\int^X$}\InS{n}}; \draw[shorten <=-3pt] (0,0)--node[inner sep=2pt,scale=0.8,blue,right] {1} (0,3/5); \draw[shorten <=-3pt] (4/5,0)--node[inner sep=2pt,scale=0.8,blue,right] {1} (4/5,3/5); \draw (17/10,0) ++(0.07,-0.14)--++(-0.26,0.53)++(0.12,-0.2) ++(0,0.02) node[blue,scale=0.7,inner sep=1,anchor=west]{2} ++(0,-0.02) ++(-0.16,0.02)--++(0.43,0.34)++(-0.13,0.07) node{$\cdot$} ++(-0.06,0.13) node{$\cdot$} ++(-0.06,0.13) node{$\cdot$}++(0.17,-0.1) node[auto,inner sep=5,scale=0.7,anchor=west]{$n\!-\!1$} ++(-0.25,0.19)--++(0.43,0.34)++(-0.18,0.04) ++(0,-0.04) node[blue,scale=0.7,inner sep=1,anchor=east]{2} ++(0,0.04) ++(0.18,-0.04)++(-0.03,-0.16)--++(-0.26,0.53); \draw[thick,shorten >=2] (41/30,8/5)--(71/15,8/5) node[inner sep=2pt,blue,above left,scale=0.8] {2}; \draw[shorten <=-3pt] (23/10,8/5)--node[inner sep=2pt,scale=0.8,blue,right] {1} (23/10,11/5); \draw[shorten <=-3pt] (31/10,8/5)--node[inner sep=2pt,scale=0.8,blue,right] {1} (31/10,11/5); \end{tikzpicture}  & \begin{tikzpicture}[xscale=0.5,yscale=0.7] \draw[thick,shorten >=2] (4/15,0)--(157/30,0) node[inner sep=2pt,blue,above left,scale=0.8] {3} node[inner sep=2pt,below left,scale=0.8] {\phantom{$\int^X$}\redtype{IV^*}}; \draw[shorten <=-3pt, shorten >=-3pt] (3/5,0)--node[inner sep=2pt,scale=0.8,blue,right] {2} (3/5,3/5); \draw[shorten <=-3pt] (3/5,3/5)--node[inner sep=2pt,scale=0.8,blue,above] {1} (0,3/5); \draw[shorten <=-3pt, shorten >=-3pt] (2,0)--node[inner sep=2pt,scale=0.8,blue,right] {2} (2,3/5); \draw[shorten <=-3pt] (2,3/5)--node[inner sep=2pt,scale=0.8,blue,above] {1} (7/5,3/5); \draw[shorten <=-3pt, shorten >=-3pt] (17/5,0)--node[inner sep=2pt,scale=0.8,blue,right] {2} (17/5,3/5); \draw[shorten <=-3pt] (17/5,3/5)--node[inner sep=2pt,scale=0.8,blue,above] {1} (14/5,3/5); \end{tikzpicture}  & \begin{tikzpicture}[xscale=0.5,yscale=0.7] \draw[thick,shorten >=2] (4/15,0)--(139/30,0) node[inner sep=2pt,blue,above left,scale=0.8] {4} node[inner sep=2pt,below left,scale=0.8] {\phantom{$\int^X$}\redtype{III^*}}; \draw[shorten <=-3pt, shorten >=-3pt] (3/5,0)--node[inner sep=2pt,scale=0.8,blue,right] {3} (3/5,3/5); \draw[shorten <=-3pt, shorten >=-3pt] (3/5,3/5)--node[inner sep=2pt,scale=0.8,blue,above] {2} (0,3/5); \draw[shorten <=-3pt] (0,3/5)--node[inner sep=2pt,scale=0.8,blue,left] {1} (0,6/5); \draw[shorten <=-3pt, shorten >=-3pt] (2,0)--node[inner sep=2pt,scale=0.8,blue,right] {3} (2,3/5); \draw[shorten <=-3pt, shorten >=-3pt] (2,3/5)--node[inner sep=2pt,scale=0.8,blue,above] {2} (7/5,3/5); \draw[shorten <=-3pt] (7/5,3/5)--node[inner sep=2pt,scale=0.8,blue,left] {1} (7/5,6/5); \draw[shorten <=-3pt] (14/5,0)--node[inner sep=2pt,scale=0.8,blue,right] {2} (14/5,3/5); \end{tikzpicture}  & \begin{tikzpicture}[xscale=0.5,yscale=0.7] \draw[thick,shorten >=2] (4/15,0)--(139/30,0) node[inner sep=2pt,blue,above left,scale=0.8] {6} node[inner sep=2pt,below left,scale=0.8] {\phantom{$\int^X$}\redtype{II^*}}; \draw[shorten <=-3pt, shorten >=-3pt] (3/5,0)--node[inner sep=2pt,scale=0.8,blue,right] {5} (3/5,3/5); \draw[shorten <=-3pt, shorten >=-3pt] (3/5,3/5)--node[inner sep=2pt,scale=0.8,blue,above] {4} (0,3/5); \draw[shorten <=-3pt, shorten >=-3pt] (0,3/5)--node[inner sep=2pt,scale=0.8,blue,left] {3} (0,6/5); \draw[shorten <=-3pt, shorten >=-3pt] (0,6/5)--node[inner sep=2pt,scale=0.8,blue,above] {2} (3/5,6/5); \draw[shorten <=-3pt] (3/5,6/5)--node[inner sep=2pt,scale=0.8,blue,right] {1} (3/5,9/5); \draw[shorten <=-3pt, shorten >=-3pt] (2,0)--node[inner sep=2pt,scale=0.8,blue,right] {4} (2,3/5); \draw[shorten <=-3pt] (2,3/5)--node[inner sep=2pt,scale=0.8,blue,above] {2} (7/5,3/5); \draw[shorten <=-3pt] (14/5,0)--node[inner sep=2pt,scale=0.8,blue,right] {3} (14/5,3/5); \end{tikzpicture}
\end{tabular}
\figuretitle Reduction types of elliptic curves
\end{center}

\noindent
Two of these, $\In{n}$ and $\InS{n}$, have a chain of $\P^1$s of any `depth' $n$,
so they are actually infinite families.
We reserve the term `reduction type' for individual types such as \In{3} or \InS{2},
and say that there are 10 \textbf{reduction families} of elliptic curves (of which 8 have just one type).
Artin and Winters \cite{AW} showed that when
$g\ge 2$, reduction types of curves with $g(C)=g$ come in finitely many\footnote
{This is not quite true in genus 1. Reduction types of (non-elliptic) curves of genus 1 are arbitrary multiples of Kodaira types,
so they technically fall into infinitely many families \cite{Saito}.}
such families
(see \S\ref{sssChains} below), and Winters \cite{Wi} gave simple combinatorial restrictions
that guarantee that a reduction type can be realised (\S\ref{ssgeotypes}).

In genus 2, Ogg \cite{Og} gave a `numerical' classification of possible fibres.
It is slightly coarser, identifying irreducible curves of the same arithmetic genus;
for example, an elliptic curve and a genus 0 curve with a node or a cusp all have type 0 in \cite{Og}.
Namikawa-Ueno \cite{NU} gave full classification over~$\C$,
and Uematsu \cite{Ue} and Ashikaga-Ichizaka \cite{AI} extended it to genus~3.

Kodaira-N\'eron classification and the finiteness of the number of types play a role
in many contexts, such as understanding the moduli space of curves \cite{DeRa, AW},
Birch--Swinnerton-Dyer conjecture \cite{Ta}, the theory of elliptic surfaces \cite{ShSurf},
and constructing curves of large rank \cite{Ne1, ShRank, El, Ul}.
It is natural to try to extend it to any genus,
especially in light of recent progress on the algorithmic side:
see, for example, \cite{newton, FN, Mu} and Donnelly's algorithm in Magma \cite{Magma},
in addition to the classical algorithms by Tate \cite{Ta} in genus 1 and Liu \cite{Liu} in genus 2.

\subsection{Outline}
\label{ssOutline}

The number of reduction families grows very fast with the genus (\ref{thm10g},~\ref{thmgg},~\ref{numfamilies}):

\begin{theoremABC}[Growth]
Write $N(g)$ for the number of reduction families in genus $g$. Then
\begin{enumerate}
\item
$N(g)\ge 10^g$ for every $g$.
\item
$N(g)\ge g^{g-o(g)}$ as $g\to \infty$.
\item
$N(2)\!=\!\gtwonumtypes$, $N(3)\!=\!\gthreenumtypes$, $N(4)\!=\!\gfournumtypes$, $N(5)\!=\!\gfivenumtypes$ and $N(6)\!=\!\gsixnumtypes$.
\end{enumerate}
\end{theoremABC}

\noindent
We will decompose reduction types into pieces called \textbf{principal types}, and explain how to classify them
and how they link.
Their number grows much slower, subexponentially in $g$.

In a break with tradition, we only work with mrnc models rather than minimal regular ones,
as they turn out to be more suitable for the purpose.
We review them in \S\ref{smodels}, and define reduction types and recall why their classification
is a combinatorial problem in~\S\ref{sredtypes}.

\subsubsection{Principal components and chains {\rm (\S\ref{sprinchains})}}
\label{sssPrinChains}

A component of the special fibre $\cC_k$ is \textbf{principal} if it has positive genus
or intersects the rest
of $\cC_k$ in $\ge 3$ points \cite{Xi}; these will usually be denoted by $\G$, $\G'$ or $\G_i$.
All other components form \textbf{chains of $\P^1$s}
either \textbf{inner} (meeting two principal components) or \textbf{outer} (meeting one),
usually unnamed:
\begin{center}
\begin{tikzpicture}[xscale=1.1,yscale=0.8]
\draw[thick,shorten >=2] (-1/3,0)--(107/30,0) node[inner sep=2pt,blue,above left,scale=0.8] {2 \smash{g}1} node[inner sep=2pt,below left,scale=0.8] {$\Gamma_1$}; \draw[shorten <=-3pt] (0,0)--node[inner sep=2pt,scale=0.8,blue,right] {1} (0,3/5); \draw[shorten <=-3pt, shorten >=-3pt] (26/15,0)--node[inner sep=2pt,scale=0.8,blue,right] {1} (26/15,3/5); \draw[shorten <=-3pt, shorten >=-3pt] (26/15,3/5)--node[inner sep=2pt,scale=0.8,blue,above] {1} (17/15,3/5); \draw[shorten <=-3pt, shorten >=-3pt] (17/15,3/5)--node[inner sep=2pt,scale=0.8,blue,left] {2} (17/15,6/5); \draw[thick,shorten >=2] (4/5,6/5)--(149/30,6/5) node[inner sep=2pt,blue,above left,scale=0.8] {3} node[inner sep=2pt,below left,scale=0.8] {$\Gamma_2$}; \draw[shorten <=-3pt, shorten >=-3pt] (26/15,6/5)--node[inner sep=2pt,scale=0.8,blue,right] {2} (26/15,9/5); \draw[shorten <=-3pt] (26/15,9/5)--node[inner sep=2pt,scale=0.8,blue,above] {1} (17/15,9/5); \draw[shorten <=-3pt, shorten >=-3pt] (47/15,6/5)--node[inner sep=2pt,scale=0.8,blue,right] {2} (47/15,9/5); \draw[shorten <=-3pt] (47/15,9/5)--node[inner sep=2pt,scale=0.8,blue,above] {1} (38/15,9/5);
\node[inner sep=2pt,align=center] at (-0.6,1) (I) {inner chain};
\draw[->,black!50] (I) edge[out=-5, in=180] (0.8,0.8);
\node[inner sep=4pt,align=center] at (-3,0.8) (O) {outer chains};
\draw[->,black!50] (O) edge[out=-20, in=180] (-0.3,0.3) (O) edge[out=20, in=180] (0.9,1.7);
\node[inner sep=4pt,align=center] at (7.5,0.6) (P) {two principal components};
\draw[->,black!50] (P) edge[out=160, in=0] (5.2,1.2) (P) edge[out=200, in=0] (3.9,0);
\end{tikzpicture}\\
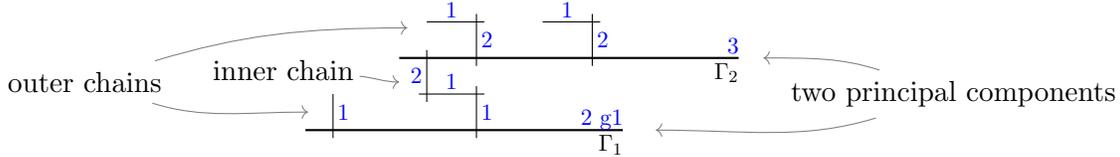
\figuretitle Principal components and chains of $\P^1$s in a special fibre.
\end{center}

\subsubsection{Structure of chains of {\rm $\P^1$s (\S\ref{schains})}}
\label{sssChains}
Classical Hirzebruch--Jung's work on rational singularities of surfaces \cite{Hi}
shows that outer chains are determined by the multiplicities of their outermost components,
and our first step is to extend this to inner chains.\\[-18pt]

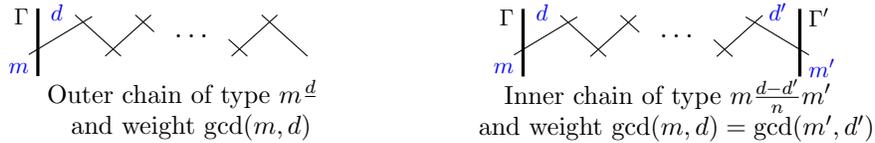
\begin{center}
\def\skewshort{edge ++(0.1,-0.2) ++(0.1,-0.2)}
\def\skewlong#1#2#3{edge node[#2,inner sep=0pt,scale=0.7]{} ++(#3) ++(#3)}
\hbox to 2.5em{\hfill}
\pbox[c]{10cm}{\begin{tikzpicture}[TPMain,scale=1.0,yscale=0.6,clabel/.style={align=center,scale=0.9}]
  \draw[thinedge] (0,0)
  edge[very thick,shorten <= 0,shorten >= 0] 
    node [very near end,left,scale=0.8]{$\G$}
    node [very near start,left,scale=0.8,blue]{$m$}
  node[above right=5pt and 1.5pt,scale=0.8,blue] {$d$}
  ++(0,5)
  ++(0,2)   
  \skewlong{\G_1}{pos=0.4,below right=1pt}{1.5,2}    
  \skewlong{\G_2}{above right}{1,-2}
  \skewlong{\G_3}{below right}{1,2}
  \skewshort ++(1.5,-0.8) node (CX) {$\ldots$} ++(1.5,-0.8) \skewshort    
  \skewlong{\G_{r\!-\!1}}{above left}{1,2}                  
  \skewlong{\G_r}{pos=0.6,below left}{1,-2};
  \coordinate (CY) at (0,-2.5);
  \node[clabel] at (CX|-CY) {Outer chain of type $m\frac{d}{}$\ \ \ \\and weight $\gcd(m,d)$};
\end{tikzpicture}}
\qquad\qquad
\pbox[c]{10cm}{\begin{tikzpicture}[TPMain,scale=1.0,yscale=0.6,clabel/.style={align=center,scale=0.9}]
  \draw[thinedge] (0,0)
  edge[very thick,shorten <= 0,shorten >= 0] 
    node [very near end,left,scale=0.8]{$\G$}
    node [very near start,left,scale=0.8,blue]{$m$}
  node[above right=5pt and 1.5pt,scale=0.8,blue] {$d$}
  ++(0,5) ++(0,2)   
  \skewlong{\G_1}{pos=0.4,below right=1pt}{1.5,2}    
  \skewlong{\G_2}{above right}{1,-2}
  \skewlong{\G_3}{below right}{1,2}
  \skewshort ++(1.5,-0.8) node (CX) {$\ldots$} ++(1.5,-0.8) \skewshort    
  \skewlong{\G_{r\!-\!1}}{above left}{1,2}                  
  \skewlong{\G_r}{pos=0.6,below left}{1.5,-2} ++(0,-2)
  edge[very thick,shorten <= 0,shorten >= 0] 
    node [very near end,right,scale=0.8]{$\G'$}
    node [very near start,right,scale=0.8,blue]{$m'$}
    node[above left=5pt and 1.5pt,scale=0.8,blue] {$d'$}
    ++(0,5);   
  \coordinate (CY) at (0,-2.5);
  \node[clabel] at (CX|-CY) {Inner chain of type $m\edge{d}{d'}{n}m'$\ \ \ \\and weight $\gcd(m,d)=\gcd(m',d')$};
\end{tikzpicture}}\end{center}

The \textbf{weight} of a chain is the gcd of multiplicities of all of its components, including the
principal component(s) to which the chain is attached \cite{Lor}.
We call chains as in the picture above an \textbf{outer chain of type $m\frac{d}{}$} and an
\textbf{inner chain of type $m\edge{d}{d'}{n}m'$}, where $d,d'$ are viewed as elements of $\Z/m\Z$ and $\Z/m'\Z$,
and the \textbf{depth} $n$ is defined as
\begin{center}
$n$ = $-1\>+$ number of times the weight occurs in the sequence $m,d,...,d',m'$.
\end{center}
Two reduction types are in the same \textbf{family}
if they coincide except possibly for their inner chain depths.
For example, Kodaira type \InS{n} has one inner chain of weight 2 and depth $n$ between
two principal components, and four outer chains of weight 1.

\begin{center}
\begin{tikzpicture}[scale=0.9,yscale=0.8]
  \def\L#1#2#3{\draw[very thick] (#3-0.2,0) -- (#3,1) node[pos=0.15,left,inner sep=2pt,blue] {#2} node[pos=0.7] (#1) {};}
  \def\D#1#2#3{\draw ($(#1)+(-0.1,0.12)$) -- ($(#1)+(0.5,-0.6)$) node[above, inner sep=4pt,pos=0.7, blue] {#3} node [pos=0.85] (#2) {};}
  \def\U#1#2#3{\draw ($(#1)+(-0.1,-0.12)$) -- ($(#1)+(0.5,0.6)$) node[below, inner sep=4pt,pos=0.7, blue] {#3} node [pos=0.85] (#2) {};}
  \def\H#1#2#3{\draw ($(#1)+(-0.1,0)$) -- ($(#1)+(0.6,0)$) node[below, inner sep=2pt,pos=0.5, blue] {#3} node [pos=0.85] (#2) {};}
  \def\RD#1#2{\draw[very thick] ($(#1 |- 0,0)+(0.16,0)$) -- ($(#1 |- 0,1)+(-0.04,0)$) node[pos=0.15,right,inner sep=2pt,blue] {#2} node[pos=0.7] {};}
  \def\RU#1#2{\draw[very thick] ($(#1 |- 0,0)+(-0.04,0)$) -- ($(#1 |- 0,1)+(0.16,0)$) node[pos=0.85,right,inner sep=2pt,blue] {#2} node[pos=0.7] {};}
  \L{a0}{6}{-3}; \D{a0}{a1}{4}; \U{a1}{a2}{2};
  \L{v0}{3}{0}; \D{v0}{v1}{4}; \U{v1}{v2}{5}; \RD{v2}{6};
  \L{w0}{3}{3.2}; \D{w0}{w1}{1}; \U{w1}{w2}{2}; \D{w2}{w3}{3}; \U{w3}{w4}{4}; \D{w4}{w5}{5}; \RU{w5}{6};
  \L{z0}{3}{8}; \D{z0}{z1}{1}; \U{z1}{z2}{1}; \D{z2}{z3}{2}; \U{z3}{z4}{3}; \D{z4}{z5}{4}; \U{z5}{z6}{5}; \RD{z6}{6};
\end{tikzpicture}\\
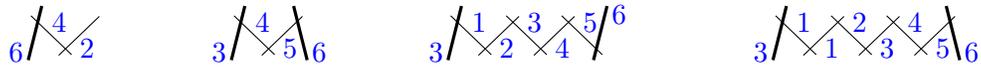
\figuretitle Outer chain of type $6\frac{4}{}$ and inner chains of type $3\edge{1}{5}{n}6$ for $n=-1,0,1$.
\end{center}

\begin{theoremABC}[Structure of chains]
\label{thmB}
Let $m,m'\!\ge\! 1$, $d\in\frac{\Z}{m\Z}, d'\in\frac{\Z}{m'\Z}$ and $n\in\Z$.
\begin{enumerate}
\item (Hirzebruch--Jung)
There are outer chains of type $m\frac{d}{}$ if and only $m\ge 2$ and $d\ne 0$.
\item
There are inner chains of type $m\edge{d}{d'}{n}m'$ if and only if
$$
  \gcd(d,m)=\gcd(d',m')
  \qquad\text{and}\qquad
   n + \tfrac{\inv(d,m)}{m} + \tfrac{\inv(d',m')}{m'} > 0,
$$
where $\inv(a,b)$ is the smallest integer $x\ge 0$ with $ax\equiv\gcd(a,b)\bmod b$.
\item
The type uniquely determines all multiplicities in the chain, both in outer and inner case.
\end{enumerate}
\end{theoremABC}

\subsubsection{Euler characteristic {\rm (\S\ref{stotalgenus})}}
\label{sssGenus}

Since chains of $\P^1$s are determined by their types, we get

\begin{theoremABC}[Decomposition]
\label{thmC}
A reduction family determines the following data, and any data satisfying these conditions
arises from a unique reduction family:
\begin{enumerate}
\item
invariants $\mgOL=\mgOLG$ for each principal component $\G\subset\cC_k$:

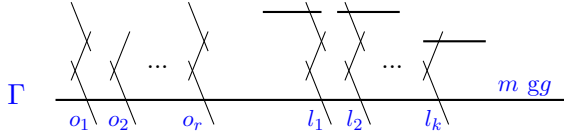
\begin{center}
\begin{tikzpicture}[TPMain,scale=1.3,edgemults/.style={above right=-0.9 and -0.14,scale=0.8,blue},
    edgemult/.style={scale=0.9},inv/.style={scale=0.9,anchor=west,align=left}]
  \draw[-] (0,0)
  edge[thickedge] ++(13,0)
  ++(0.5,1) edge[ethinedge] ++(0.3,-1) node[edgemults] {$\!o_1$} ++(0,0) edge[thinedge] ++(0.3,1) ++(0.3,1) edge[thinedge] ++(-0.3,1) ++(0.2,-2)
  ++(0.5,1) edge[ethinedge] ++(0.3,-1) node[edgemults] {$\!o_2$} ++(0,0) edge[thinedge] ++(0.3,1) ++(0.3,1) ++(0.2,-2)
  node[edgemult,above right=0.3 and 0.05] {$...$} ++(1,0)
  ++(0.5,1) edge[ethinedge] ++(0.3,-1) node[edgemults] {$\!\!o_r$} ++(0,0) edge[thinedge] ++(0.3,1) ++(0.3,1) edge[thinedge] ++(-0.3,1) ++(0.2,-2)
  ++(2,0)
  ++(0.5,1) node[edgemults] (CX) {$l_1$} ++(0,0) edge[ethinedge] ++(0.3,-1) edge[thinedge] ++(0.3,1) ++(0.3,1) edge[thinedge] ++(-0.3,1) ++(-1.5,1) edge[thickedge] ++(1.5,0) ++(1.5,-1) ++(0.2,-2)
  ++(0.5,1) node[edgemults] {$l_2$} ++(0,0) edge[ethinedge] ++(0.3,-1) edge[thinedge] ++(0.3,1) ++(0.3,1) edge[thinedge] ++(-0.3,1) ++(-0.6,1) edge[thickedge] ++(1.6,0) ++(0.6,-1) ++(0.2,-2)
  node[edgemult,above right=0.3 and 0.05] {$...$} ++(1,0)
  ++(0.5,1) node[edgemults] {$l_k$} ++(0,0) edge[ethinedge] ++(0.3,-1) edge[thinedge] ++(0.3,1) ++(-0.1,1) edge[thickedge] ++(1.6,0) ++(0.1,-1) ++(0.3,1) ++(0.2,-2)
  ++(2,0)
  node[mainmult,above=-1.5pt,scale=1.1,blue] {$m\ {\rm g}g$} 
  ++(-13.0,0.2) node[scale=1,blue] {$\G$}
  ;
  \node[inv] at (15,1.4) {
    $m\ge 1$ multiplicity of $\G$ in $\cC_k$\\
    $g\ge 0$ geometric genus of $\G$\\
    $\cO=\{o_1,...,o_r\}$ multiset of initial\\
    multiplicities of outer chains in $\smash{\frac{\Z}{m\Z}}\!\setminus\!\{0\}$\\
    $\cL=\{l_1,...,l_k\}$ multiset of initial\\
    multiplicities of inner chains in $\smash{\frac{\Z}{m\Z}}$
  };
\end{tikzpicture}
\end{center}
\noindent
such that the elements of $\cO\cup\cL$ sum to 0 mod $m$, and either $|\cO|+|\cL|\ge 3$ or $g>0$.

\item
inner chains between components, specified by a
decomposition of $\coprod_\G \cL_\G$ into pairs,
such that for every pair $(d\in\cL_\G,d'\in\cL_{\G'})$, we have $\gcd(d,m_\G)=\gcd(d',m_{\G'})$.
\end{enumerate}
Every reduction type in the family is obtained by allocating a depth for every inner chain,
subject to the minimum depth condition of Theorem \ref{thmB} (2).
\end{theoremABC}
It turns out that the total genus of the reduction type only depends on the invariants in (1):

\begin{theoremABC}
\label{thmD}
There is a unique invariant $\chi_\G = \text{some function}\,(m_\G,g_\G,\cO_\G,\cL_\G)\in\Z$
such that for every curve $C/K$, we have
$$
  \chi_{\G_1}+\ldots+\chi_{\G_n}
  = 2-2g(C),
  \eqno{\text{\rm(total genus formula)}}
$$
where $\G_1,...,\G_n$ are the principal components of the mrnc model of $C$. Explicitly,
$$
\chi_\G
 = m_\G\!\left( 2\!-\!2g_\G
   \!-\!\substack{\scriptstyle\text{\rm number of}\\[0pt]\scriptstyle\text{\rm chains from }\G}
   \right)
   + \substack{\scriptstyle\text{\rm sum of weights of}\\[0pt]\scriptstyle\text{\rm outer chains from }\G}
= m_\G\big( 2\!-\!2g_\G \!-\! |\cO_\G|\!-\!|\cL_\G|\big) \!+\!\! {\sum_{o\in\cO_\G}}\!\gcd(o,m).
$$
\end{theoremABC}

Uniqueness is easy (glue two identical copies of $\G$ to form a reduction type),
and additivity for this choice of $\chi_\G$ is proved in \ref{eulerformula}.
We call $\chi_\G$ the \textbf{Euler characteristic} of $\G$.
It might well have some geometric `stacky' interpretation, and it turns out
that $\chi_\G\le 0$. Apart from one exception called a D-tail, it is always
strictly negative, and there are finitely many tuples $\mgOL$ with given $\chi_\G<0$.
This already implies Artin-Winters' theorem
that there are finitely many reduction families in every genus $g(C)>1$.

\begin{theoremABC}[Invariants of principal components \S\ref{sprininv}]
\label{thmE}
\noindent\par\noindent
\begin{enumerate}
\item
$\chi_\G\le 0$ for every principal component $\G$. 
\item
$\chi_\G=0$ for a principal component $\G$ with invariants $\mgOL$ if and only if either
\begin{itemize}
\item
$g(C)=1$ and $\G$ is the unique principal component of $\cC_k$, which is a multiple of a Kodaira type
\redtype{I_{g1}}, \redtype{I^*_0}, \redtype{IV}, \redtype{IV^*}, \redtype{III}, \redtype{III^*}, \redtype{II}, \redtype{II^*};
or
\item
$\mgOL\!=\!(2n,0,\{n,n\},\{2n\})$ for some $n\!\ge\!1$, in which case we call $\G$ a \textbf{D-tail}:
\smallskip

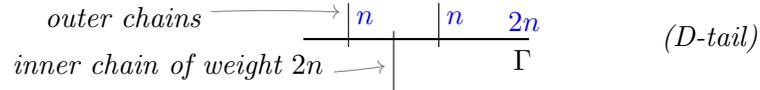
\begin{center}
\begin{tikzpicture}[scale=0.6]
  \draw[thick] (0,0) -- (5,0) node[below right=0pt and -10pt] {$\Gamma$} node[above right=-1pt and -12pt,blue] {$2n$};
  \draw[shorten <=-3pt] (1,0) -- (1,0.8) node[above right=-12pt and -1pt,blue] {$n$};
  \draw[shorten <=-3pt] (3,0) -- (3,0.8) node[above right=-12pt and -1pt,blue] {$n$};
  \draw[shorten <=-3pt] (2,0) -- (2,-1.2) node[blue,below right=-12pt and -1pt] {};
  \node[align=center] (inner) at (-3,-0.6) {inner chain of weight $2n$};
  \draw[->,black!50] (inner) edge[bend right=2] (1.8,-0.6);
  \node (outer) at (-4,0.5) {outer chains};
  \draw[->,black!50] (outer) edge[bend left=2] (0.8,0.6);
  \node at (9,0) {(D-tail)};
\end{tikzpicture}
\end{center}
\end{itemize}
\item
There are finitely many possible tuples $\mgOL$ of invariants of principal components with fixed $\chi\!<\!0$.
They satisfy the following bounds, sharp for all $\chi<0$:
\begin{center}
   $m\le -6\chi$,\ \ $g\le \lfloor 1\!-\!\chi/2\rfloor$,\ \ $|\cO|\le 4\!-\!\chi$,\ \  $|\cL|\le 2\!-\!\chi$.
\end{center}
\end{enumerate}
\end{theoremABC}

In addition to the bound $m\le -6\chi$ of (3),
we prove that $\frac{m}{-\chi}$ (for $\chi\!<\!0$) takes values (\ref{thmlargem})
\begin{center}
$\le 2$,\quad
or $(2\!+\!\frac 3k)_{k\ge 1}=\ldots,
2\frac{3}{4}, 3, 3\frac{1}{2}, 5
$, \quad
or $(2\!+\!\frac 4k)_{k\ge 1}=\ldots,
3, 3\frac{1}{3}, 4, 6
$.
\end{center}

\noindent
See Table \ref{tablargem} for the corresponding principal components.
They have `large' multiplicity $m=m_\G$ relative to $\chi=\chi_\G$.
This includes all those with $m>4g(C)\!-\!4$ in $\cC_k$ (when $g(C)>1$), that is
$m>4$ in genus 2, $m>8$ in genus 3, etc.

\subsubsection{Principal types {\rm (\S\ref{sprininv}, \S\ref{sshapes})}}
\label{sssPrinTypes}

There are two points of view what to consider the basic blocks that constitute reduction types.
For example, consider
D-tails $\G_1,\G_2,\G_3,\G_4$ linked as follows:

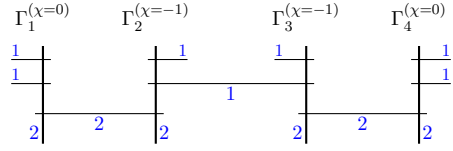
\begin{center}
\begin{tikzpicture}[scale=1]
  \def\vcomp#1#2#3{\draw[thick] (#1,0)-- node[blue,pos=0.11,inner sep=1.5pt,scale=0.7,#3] {2} (#1,1.3) node[pos=1,above=3pt,scale=0.7] {$#2$}}
  \def\hcompright#1#2{\draw[shorten <=-3pt] (#1,#2)--++(12pt,0) node[blue,pos=0.85,inner sep=2pt,above,scale=0.6] {1}}
  \def\hcompleft#1#2{\draw[shorten <=-3pt] (#1,#2)--++(-12pt,0) node[blue,pos=0.85,inner sep=2pt,above,scale=0.6] {1}}
  \def\edge#1#2#3#4{\draw[shorten <=-3pt,shorten >=-3pt] (#1,#2) --++ (#3,0) node[blue,pos=0.5,inner sep=1.5pt,below,scale=0.7] {#4}}
  \vcomp{0}{\G_1^{(\chi=0)}}{left};
    \hcompleft{0}{0.8};    \hcompleft{0}{1.12};    \edge{0}{0.4}{1.5}{2};
  \vcomp{1.5}{\G_2^{(\chi=-1)}}{right};
    \hcompright{1.5}{1.12};
    \edge{1.5}{0.8}{2}{1};
  \vcomp{3.5}{\G_3^{(\chi=-1)}}{left};
    \hcompleft{3.5}{1.12};   \edge{3.5}{0.4}{1.5}{2};
  \vcomp{5}{\G_4^{(\chi=0)}}{right};
    \hcompright{5}{0.8}; \hcompright{5}{1.12};
\end{tikzpicture}
\end{center}

\noindent
The $\G_i$ are the principal components, with $\chi=0,-1,-1,0$ as shown, and
$0\!-\!1\!-\!1\!+\!0=2\!-\!2\cdot 2$, so this is reduction type in genus 2.
We could

\begin{itemize}
\item[(a)] break it into its four principal components, and call it something like \redtype{D\e{2} D\e D\e{2} D}, or
\item[(b)] view D-tails as belonging to the two components with $\chi<0$, and call it, say, \redtype{I^*_2\e I^*_2}.
\end{itemize}

\noindent
We will go with option (b), and call this type \redtype{I^*_2\e I^*_2}. Generally, define

\begin{itemize}
\item
A \textbf{principal type} is a principal component $\G$ such that $\chi_\G<0$, taken together with
its \textbf{loops} (inner chains from $\G$ to itself) and \textbf{D-tails} as in Theorem \ref{thmE}(2).
\item
An \textbf{edge} is an inner chain between two principal types.
\end{itemize}

One reason to choose (b) over (a) is geometric: principal types seem to behave better in extensions than principal
components. Lorenzini \cite{Lo2} proved that an elliptic curve with potentially good reduction
may have type $\InS{n}$ when $\vchar k=2$. Similarly, genus 2 reduction types such as

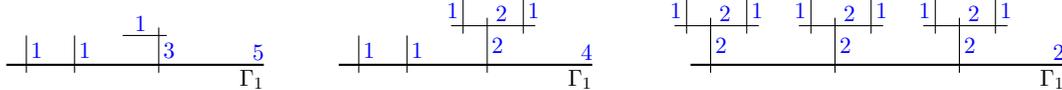
\begin{center}
\begin{tikzpicture}[xscale=0.8,yscale=0.65] \draw[thick,shorten >=2] (-1/3,0)--(121/30,0) node[inner sep=2pt,blue,above left,scale=0.8] {5} node[inner sep=2pt,below left,scale=0.8] {$\Gamma_1$}; \draw[shorten <=-3pt] (0,0)--node[inner sep=2pt,scale=0.8,blue,right] {1} (0,3/5); \draw[shorten <=-3pt] (4/5,0)--node[inner sep=2pt,scale=0.8,blue,right] {1} (4/5,3/5); \draw[shorten <=-3pt, shorten >=-3pt] (11/5,0)--node[inner sep=2pt,scale=0.8,blue,right] {3} (11/5,3/5); \draw[shorten <=-3pt] (11/5,3/5)--node[inner sep=2pt,scale=0.8,blue,above] {1} (8/5,3/5); \end{tikzpicture}  \qquad \begin{tikzpicture}[xscale=0.8,yscale=0.65] \draw[thick,shorten >=2] (-1/3,0)--(119/30,0) node[inner sep=2pt,blue,above left,scale=0.8] {4} node[inner sep=2pt,below left,scale=0.8] {$\Gamma_1$}; \draw[shorten <=-3pt] (0,0)--node[inner sep=2pt,scale=0.8,blue,right] {1} (0,3/5); \draw[shorten <=-3pt] (4/5,0)--node[inner sep=2pt,scale=0.8,blue,right] {1} (4/5,3/5); \draw[shorten <=-3pt, shorten >=-3pt] (32/15,0)--node[inner sep=2pt,scale=0.8,blue,right] {2} (32/15,4/5); \draw (23/15,4/5)--node[inner sep=2pt,scale=0.8,blue,above,pos=0.6] {2} (44/15,4/5); \draw[shorten <=-3pt] (26/15,4/5)--node[inner sep=2pt,scale=0.8,blue,left] {1} (26/15,7/5); \draw[shorten <=-3pt] (41/15,4/5)--node[inner sep=2pt,scale=0.8,blue,right] {1} (41/15,7/5); \end{tikzpicture}  \qquad \begin{tikzpicture}[xscale=0.8,yscale=0.65] \draw[thick,shorten >=2] (4/15,0)--(197/30,0) node[inner sep=2pt,blue,above left,scale=0.8] {2} node[inner sep=2pt,below left,scale=0.8] {$\Gamma_1$}; \draw[shorten <=-3pt, shorten >=-3pt] (3/5,0)--node[inner sep=2pt,scale=0.8,blue,right] {2} (3/5,4/5); \draw (0,4/5)--node[inner sep=2pt,scale=0.8,blue,above,pos=0.6] {2} (7/5,4/5); \draw[shorten <=-3pt] (1/5,4/5)--node[inner sep=2pt,scale=0.8,blue,left] {1} (1/5,7/5); \draw[shorten <=-3pt] (6/5,4/5)--node[inner sep=2pt,scale=0.8,blue,right] {1} (6/5,7/5); \draw[shorten <=-3pt, shorten >=-3pt] (8/3,0)--node[inner sep=2pt,scale=0.8,blue,right] {2} (8/3,4/5); \draw (31/15,4/5)--node[inner sep=2pt,scale=0.8,blue,above,pos=0.6] {2} (52/15,4/5); \draw[shorten <=-3pt] (34/15,4/5)--node[inner sep=2pt,scale=0.8,blue,left] {1} (34/15,7/5); \draw[shorten <=-3pt] (49/15,4/5)--node[inner sep=2pt,scale=0.8,blue,right] {1} (49/15,7/5); \draw[shorten <=-3pt, shorten >=-3pt] (71/15,0)--node[inner sep=2pt,scale=0.8,blue,right] {2} (71/15,4/5); \draw (62/15,4/5)--node[inner sep=2pt,scale=0.8,blue,above,pos=0.6] {2} (83/15,4/5); \draw[shorten <=-3pt] (13/3,4/5)--node[inner sep=2pt,scale=0.8,blue,left] {1} (13/3,7/5); \draw[shorten <=-3pt] (16/3,4/5)--node[inner sep=2pt,scale=0.8,blue,right] {1} (16/3,7/5); \end{tikzpicture}
\end{center}

\noindent
sometimes belong to potentially good curves when $\vchar k=2$. Their number of principal components fluctuates,
but they all have \textbf{only one} principal type, with $\chi=2-2g(C)=-2$ and a varying number of D-tails.
Generally, it seems that potentially good curves have one principal type in every genus.

In addition, (b) puts all Kodaira types on an equal footing, and
Ogg observed (cf. Kod($\G$) notation in \cite{Og}) that this makes the genus 2
classification cleaner. This is the convention adopted by Namikawa-Ueno as well, and it seems to work well in any~genus.

Here is an example that has three principal types, with outer chains, loops and D-tails:

\begin{center}
\tikzset{lab/.style={scale=0.9}}
\begin{tikzpicture}[xscale=0.9,yscale=0.7] \draw[thick,shorten >=2] (4/15,0)--(111/20,0) node[inner sep=2pt,blue,above left,scale=0.8] {4} node[inner sep=2pt,below left,scale=0.8] {$\Gamma_3$}; \draw[shorten <=-3pt, shorten >=-3pt] (3/5,0)--node[inner sep=2pt,scale=0.8,blue,right] {3} (3/5,3/5); \draw[shorten <=-3pt, shorten >=-3pt] (3/5,3/5)--node[inner sep=2pt,scale=0.8,blue,above] {2} (0,3/5); \draw[shorten <=-3pt] (0,3/5)--node[inner sep=2pt,scale=0.8,blue,left] {1} (0,6/5); \draw[shorten <=-3pt, shorten >=-3pt] (2,0)--node[inner sep=2pt,scale=0.8,blue,right] {2} (2,3/5); \draw[shorten <=-3pt, shorten >=-3pt] (2,3/5)--node[inner sep=2pt,scale=0.8,blue,above] {2} (7/5,3/5); \draw[shorten <=-3pt, shorten >=-3pt] (7/5,3/5)--node[inner sep=2pt,scale=0.8,blue,left] {2} (7/5,6/5); \draw (4/5,6/5)--node[inner sep=2pt,scale=0.8,blue,above,pos=0.6] {2} (11/5,6/5); \draw[shorten <=-3pt] (1,6/5)--node[inner sep=2pt,scale=0.8,blue,left] {1} (1,9/5); \draw[shorten <=-3pt] (2,6/5)--node[inner sep=2pt,scale=0.8,blue,right] {1} (2,9/5); \draw[shorten <=-3pt, shorten >=-3pt] (223/60,0)--node[inner sep=2pt,scale=0.8,blue,right] {3} (223/60,3/5); \draw[shorten <=-3pt, shorten >=-3pt] (223/60,3/5)--node[inner sep=2pt,scale=0.8,blue,below left=-1pt and -1pt] {2} (49/15,21/20); \draw[thick,shorten >=2] (44/15,21/20)--(57/10,21/20) node[inner sep=2pt,blue,above left,scale=0.8] {1 \smash{g}1} node[inner sep=2pt,below left,scale=0.8] {$\Gamma_2$}; \draw[shorten <=-3pt, shorten >=-3pt] (58/15,21/20)--node[inner sep=2pt,scale=0.8,blue,right] {2} (58/15,37/20); \draw[thick,shorten >=2] (53/15,37/20)--(69/10,37/20) node[inner sep=2pt,blue,above left,scale=0.8] {5} node[inner sep=2pt,below left,scale=0.8] {$\Gamma_1$}; \draw[shorten <=-3pt, shorten >=-3pt] (67/15,37/20)--node[inner sep=2pt,scale=0.8,blue,left] {1} (67/15,49/20); \draw[shorten <=-3pt, shorten >=-3pt] (67/15,49/20)--node[inner sep=2pt,scale=0.8,blue,above] {1} (76/15,49/20); \draw[shorten <=-3pt, shorten >=-3pt] (76/15,37/20)--node[inner sep=2pt,scale=0.8,blue,right] {2} (76/15,49/20);
\node[lab,inner sep=2pt,align=center] at (-2,0.5) (O) {outer chain};
\draw[->,black!50] (O) edge (-0.4,0.75);
\node[lab,inner sep=1pt,align=center] at (10,0.6) (I) {Principal types $\G_1$,$\G_2$,$\G_3$\\with edges between them\\of type $5\edge120 1$ and $1\edge130 4$};
\draw[->,black!50] (I) edge (7.2,1.8)      
                        (I) edge (5.9,1)        
                        (I) edge (6,0);         
\node[lab,inner sep=4pt,align=center] at (1,2.5) (D) {D-tail};
\draw[->,black!50] (D) edge[bend left=0] (1.5,1.7);
\node[lab,inner sep=4pt,align=center] at (2.8,2.5) (L) {Loop};
\draw[->,black!50] (L) edge[bend left=5] (4.1,2.4);
\end{tikzpicture}\\[4pt]
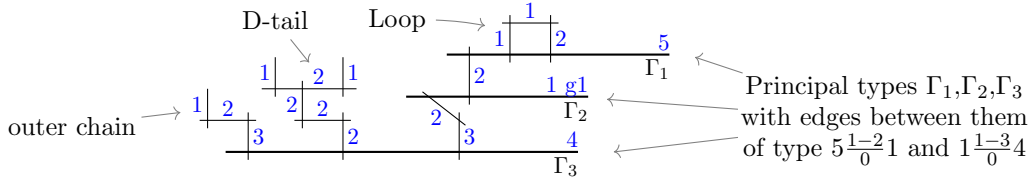
\figuretitle Reduction type
\redtype{5^{1,2,2}_1\e I_{g1}\e III^*_{3D}}
of genus
6.
\end{center}

\subsubsection{Canonical label {\rm (\S\ref{slabel})}}

\label{sssLabel}

We propose a naming convention,
attaching a canonical label to every reduction type that determines it uniquely,
and incorporates Kodaira's names.
In addition to
\redtype{5^{1,2,2}_1\e I_{g1}\e III^*_{3D}}
above, here are a few other examples, in genus
3, 3, 4,
6, 7 and 11.

\label{intropage}
\begin{center}
\begin{tabular}{c@{\qquad\quad}c@{\qquad\quad}c}
\cbox{\begin{tikzpicture}[xscale=0.8,yscale=0.6] \draw[thick,shorten >=2] (4/15,0)--(143/30,0) node[inner sep=2pt,blue,above left,scale=0.8] {6}; \draw[shorten <=-3pt, shorten >=-3pt] (3/5,0)--node[inner sep=2pt,scale=0.8,blue,right] {5} (3/5,3/5); \draw[shorten <=-3pt, shorten >=-3pt] (3/5,3/5)--node[inner sep=2pt,scale=0.8,blue,above] {4} (0,3/5); \draw[shorten <=-3pt, shorten >=-3pt] (0,3/5)--node[inner sep=2pt,scale=0.8,blue,left] {3} (0,6/5); \draw[shorten <=-3pt, shorten >=-3pt] (0,6/5)--node[inner sep=2pt,scale=0.8,blue,above] {2} (3/5,6/5); \draw[shorten <=-3pt] (3/5,6/5)--node[inner sep=2pt,scale=0.8,blue,right] {1} (3/5,9/5); \draw[shorten <=-3pt] (7/5,0)--node[inner sep=2pt,scale=0.8,blue,right] {3} (7/5,3/5); \draw[shorten <=-3pt, shorten >=-3pt] (47/15,0)--node[inner sep=2pt,scale=0.8,blue,right] {4} (47/15,3/5); \draw[shorten <=-3pt, shorten >=-3pt] (47/15,3/5)--node[inner sep=2pt,scale=0.8,blue,above] {2} (38/15,3/5); \draw[shorten <=-3pt, shorten >=-3pt] (38/15,3/5)--node[inner sep=2pt,scale=0.8,blue,left] {2} (38/15,6/5); \draw[thick,shorten >=2] (11/5,6/5)--(167/30,6/5) node[inner sep=2pt,blue,above left,scale=0.8] {4}; \draw[shorten <=-3pt] (47/15,6/5)--node[inner sep=2pt,scale=0.8,blue,right] {1} (47/15,9/5); \draw[shorten <=-3pt] (59/15,6/5)--node[inner sep=2pt,scale=0.8,blue,right] {1} (59/15,9/5); \end{tikzpicture} } & \cbox{\begin{tikzpicture}[xscale=0.8,yscale=0.6] \draw[thick,shorten >=2] (-1/3,0)--(127/30,0) node[inner sep=2pt,blue,above left,scale=0.8] {2 \smash{g}1}; \draw[shorten <=-3pt] (0,0)--node[inner sep=2pt,scale=0.8,blue,right] {1} (0,3/5); \draw[shorten <=-3pt] (4/5,0)--node[inner sep=2pt,scale=0.8,blue,right] {1} (4/5,3/5); \draw[shorten <=-3pt] (8/5,0)--node[inner sep=2pt,scale=0.8,blue,right] {1} (8/5,3/5); \draw[shorten <=-3pt] (12/5,0)--node[inner sep=2pt,scale=0.8,blue,right] {1} (12/5,3/5); \end{tikzpicture} } & \cbox{\begin{tikzpicture}[xscale=0.8,yscale=0.6] \draw[thick,shorten >=2] (-1/3,0)--(107/30,0) node[inner sep=2pt,blue,above left,scale=0.8] {2}; \draw[shorten <=-3pt] (0,0)--node[inner sep=2pt,scale=0.8,blue,right] {1} (0,3/5); \draw[shorten <=-3pt, shorten >=-3pt] (4/5,0)--node[inner sep=2pt,scale=0.8,blue,right] {1} (4/5,8/5); \draw[blue!80!white,densely dashed] (29/15,0)--(29/15,4/5);\draw[thick,shorten >=2] (8/5,4/5)--(131/30,4/5) node[inner sep=2pt,blue,above left,scale=0.8] {2 \smash{g}1}; \draw[blue!80!white,densely dashed] (38/15,4/5)--(38/15,8/5);\draw[thick,shorten >=2] (7/15,8/5)--(25/6,8/5) node[inner sep=2pt,blue,above left,scale=0.8] {2}; \draw[shorten <=-3pt] (7/5,8/5)--node[inner sep=2pt,scale=0.8,blue,right] {1} (7/5,11/5); \end{tikzpicture} } \Bcr
Type \redtype{III\e{2-4}(1)II^*} & Type \redtype{I^*_{0,g1}} & Type \redtype{[2]I_{g1}\e D\e D\e c_1} \Tcr
\cbox{\begin{tikzpicture}[xscale=0.8,yscale=0.6] \draw[thick,shorten >=2] (-1/3,0)--(14/3,0) node[inner sep=2pt,blue,above left,scale=0.8] {6}; \draw[shorten <=-3pt] (0,0)--node[inner sep=2pt,scale=0.8,blue,right] {2} (0,3/5); \draw[shorten <=-3pt, shorten >=-3pt] (4/5,0)--node[inner sep=2pt,scale=0.8,blue,left] {1} (4/5,3/5); \draw[shorten <=-3pt, shorten >=-3pt] (4/5,3/5)--node[inner sep=2pt,scale=0.8,blue,above left=0pt and -1pt] {1} (5/4,21/20); \draw[shorten <=-3pt, shorten >=-3pt] (17/10,3/5)--node[inner sep=2pt,scale=0.8,blue,above right=0pt and -1pt] {1} (5/4,21/20); \draw[shorten <=-3pt, shorten >=-3pt] (17/10,0)--node[inner sep=2pt,scale=0.8,blue,right] {1} (17/10,3/5); \draw[shorten <=-3pt, shorten >=-3pt] (91/30,0)--node[inner sep=2pt,scale=0.8,blue,right] {2} (91/30,4/5); \draw (73/30,4/5)--node[inner sep=2pt,scale=0.8,blue,above,pos=0.6] {2} (23/6,4/5); \draw[shorten <=-3pt] (79/30,4/5)--node[inner sep=2pt,scale=0.8,blue,left] {1} (79/30,7/5); \draw[shorten <=-3pt] (109/30,4/5)--node[inner sep=2pt,scale=0.8,blue,right] {1} (109/30,7/5); \end{tikzpicture} } & \cbox{\begin{tikzpicture}[xscale=0.8,yscale=0.6] \draw[thick,shorten >=2] (7/60,0)--(41/12,0) node[inner sep=2pt,blue,above left,scale=0.8] {1 \smash{g}3}; \draw[shorten <=-3pt, shorten >=-3pt] (9/20,0)--node[inner sep=2pt,scale=0.8,blue,right] {1} (9/20,32/35); \draw[shorten <=-3pt, shorten >=-3pt] (9/20,32/35)--node[inner sep=2pt,scale=0.8,blue,below left=-1pt and -1pt] {1} (0,8/5); \draw[shorten <=-3pt, shorten >=-3pt] (19/12,0)--node[inner sep=2pt,scale=0.8,blue,right] {1} (19/12,4/5); \draw[thick,shorten >=2] (5/4,4/5)--(241/60,4/5) node[inner sep=2pt,blue,above left,scale=0.8] {1 \smash{g}1}; \draw[shorten <=-3pt, shorten >=-3pt] (131/60,4/5)--node[inner sep=2pt,scale=0.8,blue,right] {1} (131/60,8/5); \draw[thick,shorten >=2] (-1/3,8/5)--(241/60,8/5) node[inner sep=2pt,blue,above left,scale=0.8] {1 \smash{g}2}; \end{tikzpicture} } & \cbox{\begin{tikzpicture}[xscale=0.8,yscale=0.6] \draw[thick,shorten >=2] (-0.33,0)--(5.08,0) node[inner sep=2pt,blue,above left,scale=0.8] {6}; \draw[shorten <=-3pt, shorten >=-3pt] (3/5,0)--node[inner sep=2pt,scale=0.8,blue,above right] {1} (3/10,3/5); \draw[shorten <=-3pt, shorten >=-3pt] (0,0)--node[inner sep=2pt,scale=0.8,blue,above left] {1} (3/10,3/5); \path[draw,thick] (4/3,0) edge[white,line width=2] ++(0.6,0) edge[out=0,in=-90] ++(0.6,0.4) ++(0.6,0.4) edge[out=90,in=0] ++(-0.3,0.3) ++(-0.3,0.3) edge[out=180,in=90] ++(-0.3,-0.3) ++(-0.3,-0.3) edge[out=-90,in=180] ++(0.6,-0.4); \draw[shorten <=-3pt, shorten >=-3pt] (3.45,0)--node[inner sep=2pt,scale=0.8,blue,right] {4} (3.45,3/5); \draw[shorten <=-3pt, shorten >=-3pt] (3.45,3/5)--node[inner sep=2pt,scale=0.8,blue,below left=-1pt and -1pt] {2} (3,21/20); \draw[thick,shorten >=2] (2.66,21/20)--(6.03,21/20) node[inner sep=2pt,blue,above left,scale=0.8] {8}; \draw[shorten <=-3pt] (3.6,21/20)--node[inner sep=2pt,scale=0.8,blue,right] {2} (3.6,33/20); \draw[shorten <=-3pt] (4.4,21/20)--node[inner sep=2pt,scale=0.8,blue,right] {4} (4.4,33/20); \end{tikzpicture} } \Bcr
Type \redtype{6^{1,1,2,2}_{3,1D}} & Type \redtype{I_{g3}\e(3)I_{g2}\e(2)I_{g1}\e(2)c_1} & Type \redtype{6^{1,1,4}_{1,1}\e [2]III} \Tcr
\end{tabular}\\[4pt]
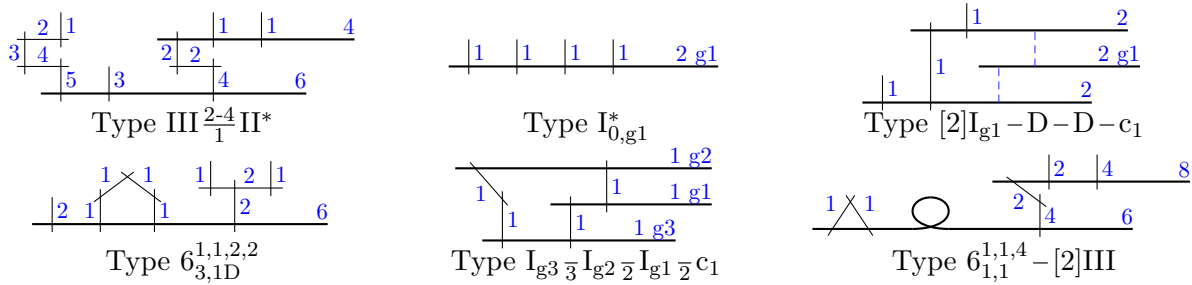
\figuretitle Some reduction types and their names.
\end{center}

Essentially, we specify a canonical path through the reduction type and the label
records principal types and types of edges along the way.
One can reconstruct the whole reduction type directly from the label.
It also gives an ordering of the principal components and chains,
which is useful if one needs a finer classification over a non-algebraically closed residue field
(e.g. split/non-split multiplicative reduction for elliptic curves), since a canonical ordering makes
it is easy to specify the Galois action.

\subsubsection{Shapes {\rm (\S\ref{sshapes})}}
\label{sssShapes}

Define the \textbf{shape} of a principal type $\G$ to be the tuple
\begin{center}
  ($\chi_\G$; weights of outgoing edges).
\end{center}
Principal types of the same shape can be thought of as interchangeable blocks: swapping one of them by another does
gives a different reduction family of the same genus. We say such reduction families have the same \textbf{shape},
and we think of the shape of a family as a graph whose vertices are $\chi$'s of principal types, and edges
are marked with weights.
This is how we will draw principal types (with inner chains, loops and D-tails
pointing up and outer chains pointing down, see \ref{notprintype}), their shapes, and the shape of
$R=
\redtype{5^{1,2,2}_1\e I_{g1}\e III^*_{3D}}
$ from above:

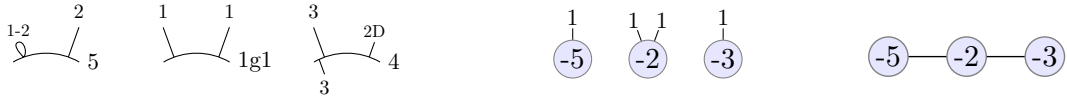
\begin{center}
\begin{tikzpicture}[scale=1,baseline=16pt]
\draw
\shapearc{25pt}{1-2,2}{}{5}
++(2,0)
\shapearc{25pt}{1,1}{}{1g1}
++(2,0)
\shapearc{25pt}{3,.2}{3}{4}
;
\node[pPShapeChi] (1) at (7,0.8) {-5}; \draw[shorten >= 4pt] (1) -- ++(90:0.55) node[pPShapeWeight] {1};
\node[pPShapeChi] (2) at (8,0.8) {-2}; \draw[shorten >= 4pt] (2) -- ++(110:0.55) node[pPShapeWeight] {1}; \draw[shorten >= 4pt] (2) -- ++(70:0.55) node[pPShapeWeight] {1};
\node[pPShapeChi] (3) at (9,0.8) {-3}; \draw[shorten >= 4pt] (3) -- ++(90:0.55) node[pPShapeWeight] {1};
\end{tikzpicture}
\qquad\qquad
\scalebox{1.3}{\begin{tikzpicture}[xscale=0.8,yscale=0.8,auto=left, v/.style={pShapeVertex}, l/.style={pShapeEdge}] \node[v] (1) at (1,1) {-5}; \node[v] (2) at (2,1) {-2}; \node[v] (3) at (3,1) {-3}; \draw[] (1)--node[l] {} (2); \draw[] (2)--node[l] {} (3); \end{tikzpicture}}
\end{center}

\noindent
In fact, there are 8 principal types of shape $(-2,[1,1])$, namely
$$
\begin{tikzpicture}[scale=0.6] \path[use as bounding box] (-1,0.5) rectangle (21,2.2); \draw(0,0) \shapearc{35pt}{1,1}{}{1g1} ++(2.75,0) \shapearc{35pt}{1-1,1,1}{}{1} ++(2.75,0) \shapearc{35pt}{1,1}{1,1}{2} ++(2.75,0) \shapearc{35pt}{1,1,.2}{}{2} ++(2.75,0) \shapearc{35pt}{1,1}{1}{3} ++(2.75,0) \shapearc{35pt}{2,2}{2}{3} ++(2.75,0) \shapearc{35pt}{1,1}{2}{4} ++(2.75,0) \shapearc{35pt}{3,3}{2}{4}; \end{tikzpicture}
$$
\begin{center}
\tabletitle\label{tabI2}%
Principal types of shape $(-2; 1,1)$.
\end{center}

\noindent
(essentially Kodaira types that are not $\II$ and $\IIS$, with two weight 1 edges)
and putting each one of them in the middle of $R$ gives a different genus 6 family.
Similarly, there are 108 principal types of shape $(-5,[1])$ (Table \ref{tabPS5}, top group) and
39 of shape $(-3,[1])$ (Table \ref{tabPS3}, top group),
and $108\cdot 8\cdot 39=33696$ genus 6 families of this shape in total.

In genus 6, the number 33696 is small, and there are shapes with 512000 families in them:

\begin{theoremABC}
The number $N(S)$ of families of a fixed shape $S$ in genus $g\ge 2$ satisfies
$$
N(S)
\le 10^g \cdot
\begin{cases}
0.55,  & \text{if } g=2,\\[-3pt]
0.44,  & \text{if } g=3,4,5,\\[-3pt]
0.512, & \text{if } g\ge 6.
\end{cases}
$$
There is a unique shape that attains these bounds when $g\le 7$, and $\lceil\frac{(g-5)^2}{6}\rceil$ of them
when $g\ge 6$.
\end{theoremABC}

Thus, there is a `small' (subexponential in~$g$) number of principal types,
`large' (exponential in~$g$) number of ways to use them as vertices to construct families
of a fixed shape, and `very large' (superexponential in~$g$) number of shapes. The latter is not very noticable
when $g$ is small, but dominates the growth when $g$ is large (even for semistable snc families, see \ref{thmgg}).

\subsection{Notation}

\label{intronotation}

All \textbf{graphs} are undirected, and multiple edges are allowed (but no loops).
In \S\S\ref{sintro}-\ref{sredtypes},
The \textbf{curves} $C/K$ we consider are smooth, complete and connected, but are \textbf{not} assumed to be
geometrically connected; otherwise there are subtle realisability issues for reduction types (see \S\ref{sredtypes}).
In \textbf{pictures of reduction types}, principal components are drawn as thick horizontal lines
with chains of $\P^1$s between them. Empty chains (two principal components meeting at a point) are
drawn as dashed lines. \textbf{Reduction families} are denoted by the reduction type in the family with minimal depths.
For example,
\redtype{I_2\e(3)I_4}
is a reduction type in the
\redtype{I_1\e I_1}-family. This family
consists of all reduction types
\redtype{I_{\text{$m$}}\e(\text{$n$})I_{\text{$r$}}}
with $m,n,r\ge 1$
(genus 2).

\begin{acknowledgements}
The author would like to thank
Raymond van Bommel,
Vladimir Dokchitser,
Qing Liu,
J\'ulia Martinez,
Allan Perez,
John Voight and
Art Waeterschoot for helpful discussions and comments on the paper.
\end{acknowledgements}

\section{Models}
\label{smodels}

\subsection{Regular and rnc models}

Let $K, O, k$ and $C/K$ be as above, and $\cC/O$ a model of $C/K$.
Recall that
$\cC$ is a \textbf{regular model} if $\cC$ is a regular scheme.
It is \textbf{regular with normal crossings (rnc)}
if $\cC$ is regular and the special fibre $\cC_k\subset\cC$ is a divisor with normal crossings, possibly non-reduced%
and with self-intersecting components.
Every curve has a \textbf{minimal regular model}, dominated by any other regular model, and a
\textbf{minimal rnc model (mrnc)}, dominated by any other rnc model.
Both are unique up to isomorphism; see \cite[\S3]{Saito} for an overview and~references.

The difference between the two is that a mrnc model can have genus 0 components of self-intersection $-1$
that meet the rest of the special fibre in $\ge 3$ points. Blowing such a component down gives another regular model
(Castelnuovo's criterion \cite[3.9]{Li}), with $\ge 3$ components meeting in a point, so it is not rnc anymore. Repeating this
process will eventually produce a minimal regular model, but it can have arbitrarily complicated
plane singularities \cite{Wi}.

\smallskip

\begin{example}
\label{exellmrnc}
For elliptic curves, below are three reduction types for which the mrnc model is different from
the minimal regular one.
Blowing down exceptional curves (dashed) from left to right gives a minimal regular model from the mrnc one,
see Morrow \cite[p. 208]{Mo}.

\smallskip
\begin{center}
\def\dashedtype#1#2#3{\draw[] (#1) \foreach\i in {#3} {   
    edge[dashedge,very thick] ++(1,0) ++(0.5,0.7) edge[thinedge] ++(0,-0.7) node[edgemult] {\i} ++(0.5,-0.7)
  } node[mainmult] (A) {#2};}
\def\threelines#1#2#3#4#5{\draw[thick] (#1)
    node[mainmult,left] {#2} to ++(2,2) ++(-1,0) to ++(0,-2) node[mainmult] {#3}
    ++(-1,2) edge[#5] ++(2,-2) ++(2,-2) node [mainmult,right] {#4};}
\def\twotouching#1#2#3#4{\draw[thick] (#1) node [mainmult,left] {#2}
    to [out=15,in=-15] ++(0,2) ++(1.146,0)
    edge[out=195,in=165,#4] ++(0,-2) ++(0,-2) node [mainmult] {#3};}
\def\cusp#1#2{\draw[thick] (#1) node [mainmult] {#2} to [out=100,in=0] ++(-1.5,1) to [out=0, in=-100] ++(1.5,1);}
\begin{tikzpicture}[scale=0.6,edgemult/.style={blue, right=-0.05, scale=0.85},mainmult/.style={blue, right, scale=0.85},
    thinedge/.style={shorten <= -15,shorten >= -8},thickedge/.style={thick},
    dashedge/.style={dash pattern=on 4pt off 1pt}
    ]
  \draw (-2,0.9) node {Type II};
  \dashedtype{2,0.3}{6}{1,2,3}         
  \threelines{6.6,0}123{dashedge}      
  \twotouching{10.1,0}12{dashedge}     
  \cusp{14,0}{1}                       
  \draw (-2,-1.8) node {Type III};
  \dashedtype{2,-2.4}{4}{2,1,1}        
  \threelines{7.9,-2.9}112{dashedge}
  \twotouching{12.6,-2.9}11{}
  \draw (-2,-4.5) node {Type IV};      
  \dashedtype{2,-5.1}{3}{1,1,1}
  \threelines{12.2,-5.6}111{}
  \draw (3.3,-6.7) node {Minimal regular with} ++(0,-0.8) node {normal crossings (mrnc)};
  \draw (13.3,-7.1) node {Minimal regular};
\end{tikzpicture}\\
\figuretitle
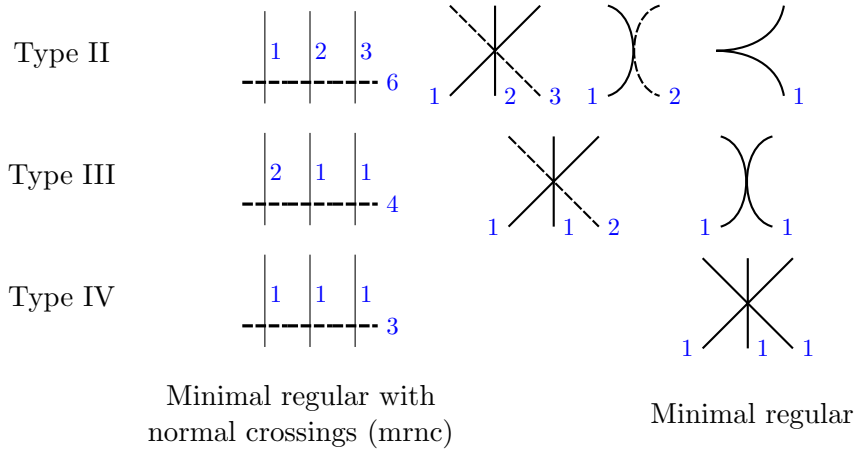
Special fibres of models for Type II, III, IV elliptic curves
\end{center}
\end{example}

\noindent
For the other families
\redtype{I_{g1}}, \In{n}, \redtype{I^*_0}, \InS{n}, \redtype{IV^*}, \redtype{III^*}, \redtype{II^*},
the mrnc and the minimal regular model agree.

\subsection{Model choice}

For classification purposes, mrnc models are nicer than minimal regular ones
in every aspect (except in their name).
They are combinatorial in nature, no analysis of complicated singularities is necessary, and they capture
arithmetic of the curve (e.g. the stable reduction field in the tame case \cite{Saito, Halle}) better.
Most importantly for us, they glue well, allowing for inductive constructions:

\begin{example}[Type \IImIVExlabel]
Here a Type II elliptic curve is linked to a Type IV one with a chain of three $\P^1$s of multiplicity 3,
and the mrnc and minimal models coincide.

\smallskip

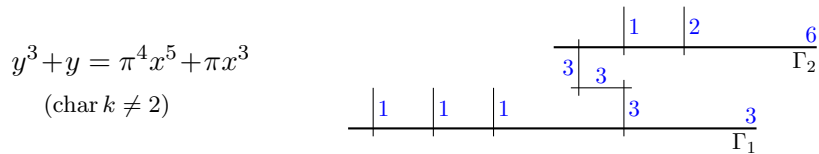
\begin{center}
\eqnbox{5cm}{\IImIVExeqn}{5pt}{1.1em}{($\vchar k\ne 2$)}
\qquad\quad
\pbox[c]{6cm}{\IImIVExfibre}
\end{center}

\noindent
Removing the bottom two $\P^1$s in the chain leaves mrnc models of the two elliptic curves, but they are not minimal regular any more.
Components $\G_1$ and $\G_2$ can now be blown down and the minimal regular models of the two elliptic
curves look very different (see Example \ref{exellmrnc}).
\end{example}

\section{Reduction types}
\label{sredtypes}

\subsection{Special fibre} \label{sspfibre}

The special fibre of any model $\cC/O$ of $C/K$ gives a divisor $[\cC_k]=\sum_\G m_\G\G$ on $\cC$.
The sum is finite, taken over irreducible curves $\G$ that we refer to as \textbf{components}.
Let

\begin{itemize}
\item
$m_\G$ = \textbf{multiplicity} of the component $\G$ in $\cC_k$,
\item
$g_\G$ = \textbf{geometric genus} of $\G$,
\item
$p_\G$ = \textbf{arithmetic genus} of $\G$.
\end{itemize}

\noindent
Suppose from now on that $\cC$ is regular with normal crossings (rnc).
Then $\Gamma$ can only have nodes (transversal self-intersection points) as singularities,
and $p_\G-g_\G\ge 0$ is the number of nodes~on~$\Gamma$.

\begin{definition}
\label{defincident}
The \textbf{incident components} to a component $\G$ is a multiset $I(\G)$ with
\begin{itemize}
\item
every component $\G'\ne\G$ taken $|\G\cap\G'|$ times,
\item
$\G$ itself taken $2(p_\G-g_\G)$ times (so every node on $\G$ contributes $\G$ to $I(\G)$ twice).
\end{itemize}
\end{definition}

\begin{center}
\begin{tikzpicture}[xscale=1,yscale=0.9]
\draw[thick,shorten >=2] (-0.33,0)--(5.46,0) node[inner sep=2pt,below left,scale=0.8] {$\Gamma_1$};
\draw[shorten <=-3pt] (0,0)--node[inner sep=2pt,scale=0.8,right] {$\Gamma_2$} (0,3/5);
\draw[shorten <=-3pt] (4/5,0)--node[inner sep=2pt,scale=0.8,right] {$\Gamma_3$} (4/5,3/5);
\draw[shorten <=-3pt, shorten >=-3pt] (29/20,0) ++(1,0) arc (0:180:0.5) node[inner sep=2pt,scale=0.8,midway,above] {$\Gamma_4$};\path[draw,thick] (91/30,0) edge[white,line width=2] ++(0.6,0) edge[out=0,in=-90] ++(0.6,0.4) ++(0.6,0.4) edge[out=90,in=0] ++(-0.3,0.3) ++(-0.3,0.3) edge[out=180,in=90] ++(-0.3,-0.3) ++(-0.3,-0.3) edge[out=-90,in=180] ++(0.6,-0.4);
\end{tikzpicture}\\
\figuretitle
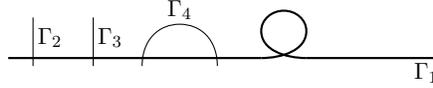
Example: $I(\G_1)=\{\G_2,\G_3,\G_4,\G_4,\G_1,\G_1\}$
\end{center}

\begin{remark}
Equivalently, take all ordinary double points $P\in\cC_s^{\red}(k)$. At such a point $P$, two components meet,
say $\G$ and $\G'$.
The component $\G$ is incident to $\G'$ and $\G'$ is incident to $\G$, so we count
one contribution of $\G'$ to $I(\G)$ and of $\G$ to~$I(\G')$ from $P$.
In particular, if $\G=\G'$, then $P$ contributes
$\G$ to $I(\G)$ twice, one for each tangent line.
\end{remark}

\subsection{Intersection pairing and adjunction}
\label{sspairing}

There is a symmetric bilinear $\Z$-valued intersection pairing on the free abelian group on the components.
It is particularly simple for rnc models because all components intersect transversally:
$$
  \G\cdot\G'=|\G\cap\G'|\quad\in\Z_{\ge 0} \rlap{\qquad for $\G\ne\G'$.}
$$
The left kernel of the pairing is of rank 1 and contains the divisor
$[\cC_k]=\sum_\G m_\G \G$. It follows that for a component $\G$,
\begin{equation}\label{selfintformula}
  \G\cdot \G \>=\> - \frac{1}{m_\G} \sum_{\G'\ne \G} m_{\G'} |\G\cap\G'|  \>=\>
  2(p_\G \!-\! g_\G) - \frac{1}{m_\G} \sum_{\G'\in I(\G)} m_{\G'}  \quad\in\Z_{\le 0}.
\end{equation}
This is strictly negative unless $\G$ is the only component of $\cC_k$.

There is a canonical divisor $K$ on $\cC$ whose intersection with every component is determined
by the adjunction formula (\cite[\S9 Thm 1.37]{Liu} or \cite[1.1(c)]{AW})
$$
  2 p_\G - 2 = \G\cdot\G + K\cdot \G.
$$
The genus of the generic fibre is given by
\cite[1.3]{AW},
\begin{equation}\label{adjunction}
  2 g(C) \!-\! 2 = \cC_k\cdot K = \sum_\G m_\G (2 p_\G \!-\! 2 \!-\! \G\cdot\G)
  = \sum_\G \Bigl(m_\G (2 g_\G \!-\!2) + \sum_{\G'\in I(\G)} m_{\G'}\Bigr).
\end{equation}
We will only need these expressions for $g(C)$, not $K$ itself or the adjunction formula.

\subsection{Reduction types}
\label{ssgeotypes}

The discrete invariants above constitute a reduction type:

\begin{definition}
\label{defredtype}
A (mrnc) \textbf{reduction type} is a tuple $R=(S; m; p; g; \cdot)$ with
\begin{itemize}
\item
$S$ a non-empty finite set, whose elements $\G$ are called components.
\item
$m\!=\!(m_\G)$, $p\!=\!(p_\G)$, $g\!=\!(g_\G)$ lists of integers indexed by $\G\in S$, with all
$m_\G\ge 1$ (multiplicity) and $p_\G\!\ge\! g_\G\!\ge\! 0$ (arithmetic/geometric genus).
\item
$\cdot: \Z S\times \Z S\to \Z$ is a symmetric bilinear pairing with $\G\cdot \G'\ge 0$ for $\G\ne \G'$, and such that
the left kernel of $\cdot$ is of rank one and contains $\sum_{\G\in S} m_\G \G$.
\item
Every $\G\in S$ with $p_\G\!=\!0$ and $\G\cdot\G\!=\!-1$ has $|I(\G)|\ge 3$, with $I(\G)$ defined by~\ref{defincident}.%
\footnote{This is the condition that the rnc model is minimal, see Proposition \ref{lemmrnc}.}
\end{itemize}
The \textbf{(total) genus} $g(R)$ is defined by \eqref{adjunction},
and \textbf{Euler chararacteristic} by $\chi(R)\!=\!2\!-\!2g(R)$.
\end{definition}

The special fibre of a mrnc model naturally gives a reduction type, and they all arise:

\begin{theorem}[Winters \cite{Wi}]
\label{winmain}
Let $R=(S; m; p; g; \cdot)$ be a reduction type, and $k$ an algebraically closed field of characteristic 0.
There is a fibration $X\to Y$ of a smooth projective surface $X/k$ over a non-singular curve $Y/k$ with
a fibre of type $R$.
\end{theorem}

Winters' theorem is more general, and allows to specify the actual curves for $\G\in S$, as long
as they are locally planar. The generic fibre of the fibration $X\to Y$ is a non-singular irreducible curve
$C/k(Y)$. Note, however, that it may not be geometrically connected! We can ask:

\begin{question}
Given a reduction type $R$ and $O, K, k\!=\!\bar k$,
is there a non-singular projective \emph{geometrically integral} curve $C/K$ whose mrnc model $\cC/O$ has special fibre of type $R$?
\end{question}

When $\vchar k=0$ and $\gcd(\{m_\G\}_\G)=1$, this answer is `yes' by Winters \cite[3.7]{Wi}.
When the $\gcd$ condition does not hold, this is much more subtle.
Liu, Lorenzini and Raynaud \cite[7.4]{LLR} showed that for `totally unipotent' types the answer is `no' when $\vchar k=0$,
and they give a complete classification in genus 1.
The question is in general open in genus $>1$, especially in mixed or positive characteristic.

\section{Principal components and chains of $\P^1$s}
\label{sprinchains}

We can now analyse reduction types purely combinatorially. We will still draw them as special fibres,
refer to components with $p_\G=0$ as $\P^1$s, and use other geometric terminology.
The term principal component is due to Xiao \cite{Xi} and inner/outer chains due to Faraggi--Nowell \cite{FN}.

\begin{definition}
\label{defprincipal}
A component $\G$ in a reduction type is \textbf{principal} if either the geometric genus $g_\G>0$ or
the number of incident components $|I(\G)|\ge 3$.
\end{definition}

\begin{definition}[Chains]
\label{defchains}
Non-principal components $\G$ have genus zero,
and meet the rest of the special fibre in 1 or 2 points. They form
\textbf{chains} of $\P^1$s of the form
$$
  \G_1-\G_2-...-\G_r,
$$
with $\G_i\cdot \G_{i+1}=1$. With the exception of one single looping chain of some multiplicity $m$ and length $n$
(Type $[m]\In{n}$), ignored from now on,
every reduction type in genus $\ge 1$ has at least one principal component (see e.g. \cite[\S3]{Ue}), and every chain is of one of
the following two types:

An \textbf{outer chain} $\G_1,..,\G_r$ is one which meets one principal component $\G_0$, and a \textbf{inner chain}
is one that meets two $\G_0$, $\G_{r+1}$ (with $\G_0=\G_{r+1}$ allowed), as follows:
$$
  \chainspicture
$$
The multiplicities of $\G_0$, $\G_{r+1}$ are part of the chain data,
and we view an intersection point of two principal components as an inner chain from $\G_0$ to $\G_1$ with~$r\!=\!0$.

We call $m_{\G_1}$ the \textbf{initial multiplicity} of the chain at $\G_0$.
In the inner case, we can also view the chain as going in the opposite direction, and call
$m_{\G_r}$ its initial multiplicity at~$\G_{r+1}$.

\end{definition}

\begin{definition}[Type]
\label{defchainstype}
Write $d_i=m_{\G_i}$ for the multiplicities of chain components.
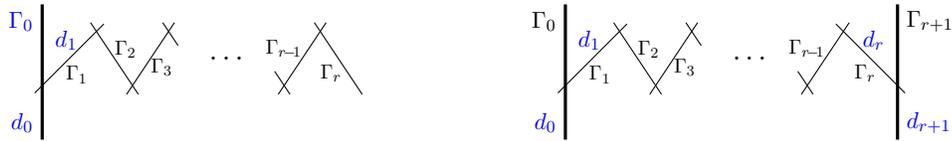
We call the chains above an \textbf{outer chain of type $d_0\frac{d_1}{}$} and
an \textbf{inner chain of type $d_0\edge{d_1}{d_r}{n}d_{r+1}$}, where:

\begin{itemize}
\item
the numbers above the line are viewed as residue classes $d_1\in\Z/d_0\Z$, $d_r\in\Z/d_{r+1}\Z$, and\\
wrtten as $1,2,...,d_0$, respectively $1,2,...,d_{r+1}$.
\item
the integer $n$ is the \textbf{depth} of the sequence $(d_i)$, defined as
\begin{center}
$n=-1\>+$ number of times $\gcd(d_0,...,d_{r+1})$ occurs in the sequence $d_0,...,d_{r+1}$.
\end{center}
\end{itemize}
\end{definition}

\begin{example}
For example, an inner chain of type $3\edge{1}{5}{n}6$ has depth $n$ and goes
from a principal component of multiplicity $m=3$ and initial multiplicity $d\equiv 1\mod 3$
to a principal component of multiplicity $m'=6$ and initial multiplicity $d'\equiv 5\mod 6$:
\begin{center}
\begin{tikzpicture}[scale=0.9]
  \def\L#1#2#3{\draw[very thick] (#3-0.2,0) -- (#3,1) node[pos=0.15,left,inner sep=2pt,blue] {#2} node[pos=0.7] (#1) {};}
  \def\D#1#2#3{\draw ($(#1)+(-0.1,0.12)$) -- ($(#1)+(0.5,-0.6)$) node[above, inner sep=4pt,pos=0.7, blue] {#3} node [pos=0.85] (#2) {};}
  \def\U#1#2#3{\draw ($(#1)+(-0.1,-0.12)$) -- ($(#1)+(0.5,0.6)$) node[below, inner sep=4pt,pos=0.7, blue] {#3} node [pos=0.85] (#2) {};}
  \def\H#1#2#3{\draw ($(#1)+(-0.1,0)$) -- ($(#1)+(0.6,0)$) node[below, inner sep=2pt,pos=0.5, blue] {#3} node [pos=0.85] (#2) {};}
  \def\RD#1#2{\draw[very thick] ($(#1 |- 0,0)+(0.16,0)$) -- ($(#1 |- 0,1)+(-0.04,0)$) node[pos=0.15,right,inner sep=2pt,blue] {#2} node[pos=0.7] {};}
  \def\RU#1#2{\draw[very thick] ($(#1 |- 0,0)+(-0.04,0)$) -- ($(#1 |- 0,1)+(0.16,0)$) node[pos=0.85,right,inner sep=2pt,blue] {#2} node[pos=0.7] {};}
  \L{v0}{3}{-3}; \D{v0}{v1}{1};
  \L{v0}{3}{0}; \D{v0}{v1}{4}; \U{v1}{v2}{5}; \RD{v2}{6};
  \L{w0}{3}{3.2}; \D{w0}{w1}{1}; \U{w1}{w2}{2}; \D{w2}{w3}{3}; \U{w3}{w4}{4}; \D{w4}{w5}{5}; \RU{w5}{6};
  \L{z0}{3}{8}; \D{z0}{z1}{1}; \U{z1}{z2}{1}; \D{z2}{z3}{2}; \U{z3}{z4}{3}; \D{z4}{z5}{4}; \U{z5}{z6}{5}; \RD{z6}{6};
\end{tikzpicture}\\
\figuretitle
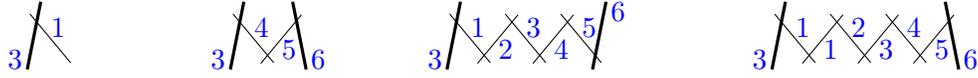
Outer chain of type $3\frac{1}{}$ and inner chains of type $3\edge{1}{5}{n}6$ with depth $n=-1,0,1$
\end{center}
\end{example}

\begin{remark}
The `$-1$' in the definition of depth is there for historical reasons, and agrees
with the Kodaira and Namikawa--Ueno conventions in genus 1 and 2.
In \cite{NU},
genus 2 reduction types with chains of depth $-1$ are called $\alpha$-types,
though Liu's algorithm \cite{Liu} uses $-1$ instead of $\alpha$.
See Table \ref{tabNU} (p. \pageref{tabNU}), `when $n=-1$' entries.
\end{remark}

In \S\ref{schains} we will see that the type determines the chain uniquely.

\begin{example}[Elliptic curves]
\label{ellcurvetableex}
For Kodaira types (p. \pageref{kodairatypes})
with the exception of $\In{n}$ (ignored, see \ref{defchains}), the principal components and chains are:

\smallskip

\begin{center}
\def\outer#1#2{#1\frac{#2}{}}
\begin{tabular}{|c|c|c|c|}
\hline
Kodaira       & Multiplicites &  Types of & Types of\cr
Type          & of principal  &  inner    & outer   \cr
              & components    &  chains   & chains  \cr
\hline
\redtype{I_{g1}} (\In{0}) & 1   &  & \cr
\IZS          & 2   &  & $\outer 21, \outer 21, \outer 21, \outer 21$ \cr
\InS{n}       & 2,\>2 & $2\edge22n2$ & $\outer 21, \outer 21, \outer 21, \outer 21$ \cr
\IV           & 3   &  & $\outer 31, \outer 31, \outer 31$ \cr
\IVS          & 3   &  & $\outer 32, \outer 32, \outer 32$ \cr
\III          & 4   &  & $\outer 41, \outer 41, \outer 42$ \cr
\IIIS         & 4   &  & $\outer 43, \outer 43, \outer 42$ \cr
\II           & 6   &  & $\outer 61, \outer 61, \outer 63$ \cr
\IIS          & 6   &  & $\outer 65, \outer 64, \outer 63$ \cr
\hline
\end{tabular}
\end{center}
\end{example}

\section{Classification of chains}  \label{schains}

\def\PP(#1,#2,#3,#4,#5){P_{#1,#3,#5}^{#2,#4}}

We now classify outer and inner chains by their types (Definition \ref{defchains}). We show that

\begin{itemize}
\item
Outer chains have type $m\frac{d}{}$ with $m\ge 2$ and $1\le d<m$. Conversely, for such $m$, $d$,
there are chains of this type, and they are unique in the sense that the number of components and
their multiplicities are uniquely determined by $m$ and $d$ (\ref{definoutseq}, \ref{outermain}).

\item
Inner chains have type $m\edge{d}{d'}{n}m'$ with $m,m'\ge 1$, $d\in\Z/m\Z$, $d'\in\Z/m'\Z$ and
$$
  \gcd(m,d)=\gcd(m',d'). \eqno{\text{(weight condition)}}
$$
Conversely, under these conditions it is possible to construct such chains, of any depth
$$
  n\ge n_{\min},     \eqno{\text{(minimal depth condition)}}
$$
with explicit $n_{\min}\in\{-1,0,1\}$ (\ref{definoutseq}, \ref{definnertype}, \ref{thmchain}).
Again, $m$, $d$, $m'$, $d'$ and $n$ determine all the
multiplicities in the chain uniquely (and explicitly, see Algorithm \ref{inneralg}).
\end{itemize}

The conditions that control multiplicities in a chain are
$\G\cdot\G\in\Z,\>\G\cdot\cC_k=0,\>\G\cdot\G\ne -1$ for every component $\G$ in the chain (\ref{defredtype}),
and translate to the following:

\begin{definition}
\label{definoutseq}
An \textbf{inner sequence} is a sequence of positive integers $d_0, d_1, \ldots, d_{r+1}$ for some $r\ge 0$, such that
$$
  \frac{d_{i-1}+d_{i+1}}{d_i}\in\Z_{\ge 2} \quad\text{for}\quad i=1,\ldots,r.
$$
An \textbf{outer sequence} is a sequence $d_0, d_1, \ldots, d_{r}$ satisfying this with $d_{r+1}=0$.
In both cases, we refer to $c=\gcd(d_0,d_1,...,d_{r+1})$ as the \textbf{weight} of the sequence, written $\weight(d_i)$.

As in Definition \ref{defchainstype}, we call $d_0\frac{d_1}{}$ the \textbf{type} of an outer sequence,
and $d_0\edge{d_1}{d_r}{n}d_{r+1}$ the \textbf{type} of an inner sequence, where
we view $d_1\in \Z/d_0\Z$, $d_r\in\Z/d_{r+1}\Z$,
and the \textbf{depth} $n\ge -1$ is
$$
  n = -1 + \text{number of occurrences of $c$ in the sequence $d_0,...,d_{r+1}$.}
$$
\end{definition}

\begin{example}
\label{exouterinner}
Here are some outer and inner sequences, all of weight 1:

\smallskip

\begin{tabular}{l@{ }l@{\qquad}l@{ }l@{\quad}l}
(a)&Outer: 8, 5, 2, 1. & (c)&Inner: 8, 5, 2, 1, 1, 1, 2, 3, 4 & (depth $2$).\cr
(b)&Outer: 4, 3, 2, 1. & (d)&Inner: 8, 5, 2, 3, 4 & (depth $-1$).\cr
\end{tabular}
\end{example}

\begin{proposition}
\label{lemmrnc}
A rnc model is mrnc if and only if the multiplicities in every inner chain form
an inner sequence,
and multiplicities in every outer chain form an outer sequence.
\end{proposition}

\begin{proof}
Recall that by the Castelnuovo's criterion, a regular model is not minimal if and only if it has a component
of genus $p_\G=g_\G=0$ with $\G\cdot\G=-1$. Such a component $\G$ can be blown down, resulting in another regular model.
If an rnc model is not mrnc, such a $\G$ meets the rest of the special fibre
in 1 or 2 points (otherwise blowing it down gives a singularity with $\ge 3$ branches, so the model is no longer rnc).
Thus $\G$ is necessarily non-principal, and is part of a chain with multiplicities as follows:

\bigskip

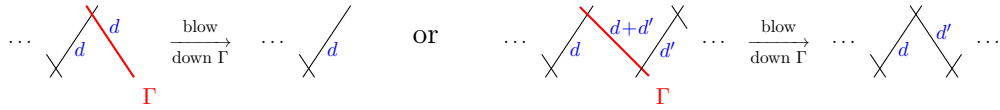
\begin{center}
\def\skewlong#1#2#3{edge node[#2,inner sep=0pt,scale=0.7,blue]{$#1$} ++(#3) ++(#3)}
\def\skewlongpar#1#2#3#4{edge[#4] node[#2,inner sep=0pt,scale=0.7,blue]{$#1$} ++(#3) ++(#3)}
\def\skewshort{edge ++(0.1,-0.2) ++(0.1,-0.2)}
\def\skewshortup{edge ++(0.1,0.2) ++(0.1,0.2)}
\pbox[c]{10cm}{\begin{tikzpicture}[TPMain,scale=1.2]
  \draw[thinedge] (0,0)
  node[scale=0.7] {$\ldots$} ++(0.9,-0.8)
  \skewshort
  \skewlong{d}{midway,below right}{1,2}
  \skewlongpar{d}{pos=0.4,above right=1pt}{1,-2}{red,thick}
  ++(0.6,-1) node[scale=0.8,red] {$\Gamma$} ++(-0.6,1)
  ++(2.0,1.0) node[scale=0.6,align=center] {blow\\[-4pt]$\xrightarrow{\qquad\quad}$\\[-3pt]down $\G$} ++(2.0,0)
  node[scale=0.7] {$\ldots$} ++(0.9,-0.8)
  \skewshort
  \skewlong{d}{midway,below right}{1,2}
  ;
\end{tikzpicture}}
\qquad\raise10pt\hbox{or}\qquad
\pbox[c]{10cm}{\begin{tikzpicture}[TPMain,scale=1.2]
  \draw[thinedge] (0,0)
  node[scale=0.7] {$\ldots$} ++(0.9,-0.8)
  \skewshort
  \skewlong{d}{midway,below right}{1,2}
  \skewlongpar{d\!+\!d'}{pos=0.4,above right}{1.5,-2}{red,thick}
  ++(0.6,-1) node[scale=0.8,red] {$\Gamma$} ++(-0.6,1)
  \skewlong{d'}{midway,below right}{1,2}
  \skewshort
  ++(0.9,-0.8) node[scale=0.7] {$\ldots$}
  ++(1.8,0.0) node[scale=0.6,align=center] {blow\\[-4pt]$\xrightarrow{\qquad\quad}$\\[-3pt]down $\G$} ++(1.8,0)
  node[scale=0.7] {$\ldots$} ++(0.9,-0.8)
  \skewshort
  \skewlong{d}{midway,below right}{1,2}
  \skewlong{d'}{midway,above right=1pt}{1,-2}
  \skewshortup
  ++(0.9,0.8) node[scale=0.7] {$\ldots$}
  ;
\end{tikzpicture}}
\end{center}

Equivalently, $\frac{d_{i-1}+d_{i+1}}{d_i}=1$, so the $d_i$ do not form an inner/outer sequence.
\end{proof}

The main result of this section is that types determine the sequences.
There are bijections (Definition \ref{definnertype}, Corollary \ref{outermain} and Theorem \ref{thmchain}(3))
$$
\begin{array}{ccc}
  \text{\{\rm inner sequences\}}&\lar&
  \text{\{\rm inner types\}}\cr
  (d_0,d_1,\ldots,d_r)
  &\longmapsto&
  d_0\edge{d_1}{d_r}{n}d_{r+1}\cr
\end{array}\quad
\begin{array}{ccc}
  \text{\{\rm outer sequences\}}&\lar&
  \text{\{\rm outer types\}}\cr
  (d_0,d_1,\ldots,d_r)
  &\longmapsto&
  d_0\frac{d_1}{}.\cr
\end{array}
$$

First, we show that the gcd of two consecutive terms is constant in the sequence $(d_i)$,
and that it is convex. In particular, it is first decreasing,
then constant, then increasing, with some of these three steps possibly omitted,
like in Example \ref{exouterinner}.

\begin{lemma}[Weight]
\label{gcdcondition}
Let $d_0,...,d_{r+1}$ be an inner or an outer sequence. Its weight
$$
  \weight(d_i) = \gcd(d_0,d_1,...,d_{r+1})
$$
is the gcd of any two consecutive terms in the sequence,
$$
  \weight(d_i) = \gcd(d_0,d_1)=\gcd(d_0,d_1)=...=\gcd(d_r,d_{r+1}).
$$
\end{lemma}

\begin{proof}
The condition $\frac{d_{i-1}+d_{i+1}}{d_i}\in\Z$ gives $d_{i-1}\equiv -d_{i+1}\bmod d_i$ for $i=1,...,r$, and in
particular $\gcd(d_{i-1},d_i)=\gcd(d_i,d_{i+1})$.
\end{proof}

\begin{lemma}[Convexity]
\label{convexity}
An inner or an outer sequence $d_0,...,d_{r+1}$ is convex.
\end{lemma}

\begin{proof}
By definition, $\frac{d_{i-1}+d_{i+1}}{d_i}\ge 2$, so
$d_i\le \frac{d_{i-1}+d_{i+1}}2$.
\end{proof}

Outer sequences are essentially the terms in Hirzebruch--Jung continued fractions (cf. \cite{Hi}),
and are determined by their two initial terms. This is well-known:

\begin{corollary}[Outer sequences]
\label{outermain}
An outer sequence $m,d,\ldots$ is strictly decreasing, has $1\le d<m$ and ends with $d_r=\gcd(m,d)$,
the weight of the sequence.
Conversely, for $1\le d<m$ integers, there is a unique outer sequence of type $m\frac{d}{}$.
\end{corollary}

\begin{proof}
In an outer sequence, $d_r|d_{r-1}$. Together with \ref{gcdcondition} and \ref{convexity}
this shows that $d_r=\gcd(d_0,d_1,...,d_r)$ and that $(d_i)$ is strictly decreasing:
$$
  d_0=m, \>\>d_1=d, \>\>d_2,...,d_{r-1}, \>\>d_r=\gcd(d_0,...,d_r).
$$
Inductively, the $d_i$ satisfy $1\le d_i<d_{i-1}$ and
$d_i\equiv-d_{i-2}\bmod d_{i-1}$ for $i>1$. This clearly has a unique solution on every step,
terminating at $\gcd(d_0,...,d_r)$, as in Euclid's algorithm.
\end{proof}

Now we turn to inner sequences and their types.

\begin{remark}
\label{reminner}
Any inner sequence, say,
$$
  m,d,...,d',m',
$$
satisfies $\gcd(m,d)=\gcd(m',d')$ by Lemma \ref{gcdcondition}.
Conversely, under this condition, let $c$ be this gcd (the weight of the sequence), and suppose
for now that $d<m$, $d'<m'$.
Take outer sequences $m,d,\ldots,c$ and $m',d',\ldots,c$, and glue them to
$$
  m,d,\ldots,c,c,\ldots,c,\ldots,d',m'
$$
with any number $n+1\ge 1$ of $c$'s. This is clearly an inner sequence of
arbitrary depth $n\ge 0$, and we will show that it is unique.

Sequences of depth $-1$ are slightly more subtle.
For example, consider
$$
  \text{(i) }8,5,\ldots,3,4  \qquad\text{and}\qquad  \text{(ii) }5,3,\ldots,3,5.
$$
An inner sequence as in (i) of depth $-1$ exists (Example \ref{exouterinner}(d)),
and one as in (ii) does not. It turns out that the precise condition is this (Theorem \ref{thmchain} below):
\end{remark}

\begin{definition}[Inner type]
\label{definnertype}
Let $m,m'\!\ge\! 1$, $d\in\frac{\Z}{m\Z}, d'\in\frac{\Z}{m'\Z}$ and $n\in\Z$ such that
$$
  \gcd(d,m)=\gcd(d',m')
  \qquad\text{and}\qquad
   n + \frac{\inv(d,m)}{m} + \frac{\inv(d',m')}{m'} > 0,
$$
where $\inv(d,m)$ is the smallest integer $x\ge 0$ with $dx\equiv\gcd(d,m)\bmod m$.
Then $m\edge{d}{d'}{n}m'$ is called an \textbf{inner type}, and $\gcd(d,m)$ its \textbf{weight}.
\end{definition}

This condition can be understood through the connection to 1-paths between rational numbers.
These are closely related to the Hirzebruch-Jung continued fractions and
resolutions of toric singularities.
Existence and uniqueness of these 1-paths is due to Obus-Wewers \cite[App. A]{OW}.

\begin{definition}\cite[A.4]{OW}.
For rational numbers $a\ge a'$,
a \textbf{1-path} is a sequence of fractions in lowest terms for some $r\ge 0$,
$$
  a = \frac{n_0}{d_0} > \frac{n_1}{d_1} > \ldots > \frac{n_{r+1}}{d_{r+1}} = a' \qquad \text{with}
    \quad n_i d_{i+1} - d_i n_{i+1} = 1, \>\>0\le i\le r.
$$
It is \textbf{shortest} if $d_i\ne d_{i-1}+d_{i+1}$ for all $i=1,...,r$.
\end{definition}

\begin{example}
\label{exHJ1}
This is an example of a shortest 1-path:
$$
  \frac{19}{23} > \frac{14}{17} > \frac{9}{11} > \frac 45 > \frac 34 > \frac 57 > \frac{7}{10} > \frac{16}{23}.
$$
\end{example}

\begin{proposition}
\label{fracmain}
\noindent\par\noindent
\begin{enumerate}
\item
If $\frac{n}{d} > \frac ab > \frac{n'}{d'}$
with $n d' - d n' = 1$, then $b\ge d+d'$.
\item
Conversely, if $n d' - d n' > 1$ with $(n,d)=(n',d')=1$, then there is $a/b$ with
$$
  \frac{n}{d} > \frac ab > \frac{n'}{d'}, \qquad b\le \max(d,d').
$$
\item \cite{OW}
For every $\frac nd>\frac{n'}{d'}$ there is a unique shortest 1-path $P$ from $\frac nd$ to $\frac{n'}{d'}$.
\item
This path $P$ is computed inductively as follows:
\begin{itemize}
\item
If $n d' - d n' = 1$ then $P=\{\frac nd,\frac{n'}{d'}\}$. Otherwise,
\item
If
$d\ge d'$, let $a,b$ be unique integers with $b n - a d = 1$ and $1\le b\le d$. Then
$$
  P=\{\frac nd\}\>\cup\>\bigl(\text{shortest 1-path from }\frac ab\text{ to }\frac{n'}{d'}\bigr).
$$
\item
If $d\le d'$, let $a,b$ be unique integers with $a d' \!-\! b n' = 1$ and $1\le b\le d'$.~Then
$$
  P=\bigl(\text{shortest 1-path from }\frac{n}{d}\text{ to }\frac ab\bigr)\>\cup\>\{\frac {n'}{d'}\}.
$$
\end{itemize}
\end{enumerate}
\end{proposition}

\begin{proof}
(1)
Indeed,
$$
  \frac{1}{dd'} = \Bigl(\frac{n}{d} - \frac{a}{b}\Bigr) + \Bigl(\frac{a}{b} - \frac{n'}{d'}\Bigr)
    \ge \frac1{bd} + \frac{1}{bd'}
    = \frac{d+d'}{b}\frac{1}{dd'}.
$$

(2)
Assume $d\ge d'$ (else change signs of the $n, n'$ and swap the two fractions).
Let $a/b$ be as in (4), so that
$$
  \frac{n}{d} > \frac ab, \qquad b n - a d = 1.
$$
In addition, $\frac{n'}{d'}>\frac ab$, for otherwise $\frac{n'}{d'}\notin [\frac ab,\frac{n}{d}]$,
contradicting (1) as $d'\le d< d+b$.

(3 - Existence, 4): Applying (2) repeatedly constructs a path as in (4). The process stops since there are only finitely
many fractions in $[\frac{n'}{d'},\frac nd]$ of denominator $\le\max(d,d')$.

(3 - Uniqueness):
Suppose there is a shortest path $P$ from $\frac{n}{d}$ to $\frac{n'}{d'}$ without $a/b$ as in (4), and say $d\ge d'$
(handling $d\le d'$ similarly, as in (2)). Thus,
$$
  \frac{n}{d} > \frac uv > \ldots > \frac{n'}{d'}.
$$
If $\frac{n}{d} > \frac ab>\frac uv$, this contradicts (1), because $b\le d<d+v$. If
$\frac{n}{d} > \frac uv > \frac ab$, then $v>d+b>d\ge d'$, so the denominators in $P$ have a local maximum
not on the boundary. This contradicts the fact that
the denominators of a shortest 1-path form an inner sequence (see proof of Theorem \ref{thmchain} Step 1 below).
\end{proof}

\begin{corollary}
The shortest 1-path from $\frac{n}{d}$ to $\frac{n'}{d'}$ contains the fractions in the interval
$[\frac{n'}{d'},\frac{n}{d}]$ of the minimal denominator.
\end{corollary}

\begin{proof}
This follows from Proposition \ref{fracmain} (1).
\end{proof}

\begin{example}
\label{exHJ2}
In Example \ref{exHJ1}, the fraction of minimal denominator in the interval $[\frac{16}{23},\frac{19}{23}]$
is $\frac 34$,
and it is a member of the shortest 1-path from $\frac{19}{23}$ to $\frac{16}{23}$.
\end{example}

Now we can prove the main result of this section:

\begin{theorem}[Inner chains]
\label{thmchain}
Consider the maps
$$
\begin{array}{ccccc}
  \Bigl\{\pbox{8em}{\ \,Shortest 1-paths\\up to integer shifts}\Bigr\}
  &\lar&
  \Bigl\{\pbox{8em}{\rm \ inner sequences\\{}\hbox to 15pt{\hfill}of weight 1}\Bigr\}
  &\lar&
  \Bigl\{\pbox{8em}{\rm inner types\\{}\hbox to 2pt{\hfill}of weight 1}\Bigr\}
  \cr
\frac{n_0}{d_0} > \frac{n_1}{d_1} > \ldots > \frac{n_{r+1}}{d_{r+1}}
  &{\buildrel \alpha \over \longmapsto}&
  (d_0,d_1,\ldots,d_r,d_{r+1})
  &{\buildrel \beta \over \longmapsto}&
  d_0\edge{d_1}{d_r}{n}d_{r+1},\cr
\end{array}
$$
where $n$ in $\beta$ is $-1\>\>+$ the number of 1's among $d_0,...,d_{r+1}$. Then
\begin{enumerate}
\item
Both $\alpha$ and $\beta$ are bijections.
\item
To compute the inverse $(\beta\circ\alpha)^{-1}$,
let $m\edge{l}{l'}{n}m'$ be an inner type of gcd~1.
Solve $a'l'\equiv -1\mod m'$ with $a'\in\Z$ and set $\frac{a}{m}=$ unique rational
$>\frac{a'}{m'}$
with denominator~$m'$, numerator $\equiv l^{-1}\bmod m$ and
$|[\frac{a'}{m'},\frac{a}{m}]\cap\Z|=n+1$.~Then
$$
  (\beta\circ\alpha)^{-1}\colon\quad
  \textstyle m\edge{l}{l'}{n}m' \quad\longmapsto\quad P_{\frac{a}{m},\frac{a'}{m'}}.
$$
\item
Dropping the weight 1 condition, we still have a bijection induced by $\beta$,
$$
\begin{array}{ccc}
  \text{\{\rm inner sequences\}}&\lar&
  \text{\{\rm inner types\}}\cr
  (d_0,d_1,\ldots,d_r)
  &\longmapsto&
  d_0\edge{d_1}{d_r}{n}d_{r+1}\cr
\end{array}
$$
\end{enumerate}
\end{theorem}

\begin{proof}
\emph{Step 1: $\alpha$ is a bijection.}
Say $(\frac{n_i}{d_i})_{i=0}^{i=r+1}$ is a 1-path, so
$n_i d_{i+1} - d_i n_{i+1} = 1$ for $i=0,...,r$.
In particular, $\gcd(n_i,d_i)=1$ and $\gcd(d_i,d_{i+1})=1$.
Combining the 1-path condition for $i$ and $i-1$, we get
\begin{equation}\label{nini}
  n_{i+1} = -n_{i-1} + n_i\frac{d_{i-1}+d_{i+1}}{d_i}.
\end{equation}
As $\gcd(n_i,d_i)=1$, this shows that $\frac{d_{i-1}+d_{i+1}}{d_i}\in \N$.
By the shortest path condition, the quotient is $\ne 1$, so that $(d_i)$ is an inner sequence,
and as noted above it has gcd 1.

Conversely, given an inner sequence $d_0,...,d_{r+1}$ of gcd 1, we can solve the equation
$$
  n_0 d_1 - d_1 n_0 = 1
$$
in integers $n_0, n_1$. Moreover, $n_0$ is unique modulo $d_0$, equivalently
$\frac{n_0}{d_0}$ is unique mod~$\Z$, and $n_1$ is unique once $n_0$ is fixed.
Then \eqref{nini} determines $n_i$ for $i>1$ inductively and uniquely,
and we find a 1-path with the given $d_i$, unique up to integer shifts.
The condition $d_i\ne d_{i-1}+d_{i+1}$ guarantees that
$(\frac{n_i}{d_i})$ is a shortest 1-path.

\emph{Step 2: $\beta\circ\alpha$ is a bijection.}
This follows from Proposition \ref{fracmain} and the explicit form of $(\beta\circ\alpha)^{-1}$. In particular $\beta$
is bijective as well, so (1) and (2) hold.

\emph{Step 3: (3) holds}. Since $\beta$ is bijective, (3) holds for weight 1 inner sequences.
Multiplying all the $d_i$ by the weight, we find that it holds for arbitrary inner sequences as well.
\end{proof}

Finally, the proof of the theorem combined with Proposition \ref{fracmain}(4) gives an algorithm
to construct the unique inner sequence of a given type:

\begin{algorithm}
\label{inneralg}
Suppose $m\edge{d}{d'}{n}m'$ is an inner type (Definition \ref{definnertype}), and
$P=\PP(m,d,m',d',n)$ is the unique inner sequence of this type, of weight $c$.
Write $\tilde x$ for the lift from $\Z/m\Z$ to $\{1,...,m\}\subset\Z$, and
from $\Z/m'\Z$ to $\{1,...,m'\}\subset\Z$.

If $n\!\ne\!-1$,~then
$$
  P = (m,\tilde d,\ldots,c,c,\ldots,c,\tilde d',\ldots m')
$$
where
\begin{enumerate}
\item
$m,\tilde d,\ldots,c$ is the (unique) outer sequence of type $m\frac{d}{}$ if $c<m$,
and $\emptyset$ if $c=m$.
\item
$c,...,\tilde d',m'$ is the reversed outer sequence of type $m'\frac{d'}{}$ if $c\!<\!m'$,
and $\emptyset$ if $c\!=\!m'$.
\item
the total number of $c$'s in $P$ is $n+1$.
\end{enumerate}
If $n=-1$, then $P$ is constructed inductively as follows:
\begin{enumerate}
\setcounter{enumi}{3}
\item
If $d'=m$ and $d=m'$, then $P=(m,m')$. Otherwise:
\item
If $m\!\ge\! m'$ then $P \!=\! (m)\cup \PP(\tilde d,-m,m',d',-1)$.
\item
If $m\!\le\! m'$ then $P \!=\! \PP(m,d,\tilde d',-m',-1)\cup (m')$.
\end{enumerate}
\end{algorithm}

\begin{example}[$\PP(8,5,4,3,-1)$ and $\PP(5,3,5,3,-1)$]
To follow up the examples in Remark \ref{reminner},
$$
  \PP(8,5,4,3,-1)=(8)\cup \PP(5,2,4,3,-1) = (8,5) \cup \PP(2,1,4,3,-1) = (8,5) \cup \PP(2,1,3,2,-1) \cup (4) =
    (8,5,2,3,4)
$$
using Step (5), (5), (6), (4) of the algorithm. Similarly, using Step (1)--(3),
$$
  \PP(8,5,4,3,0) = (8,5,2,1,2,3,4).
$$
By Theorem \ref{thmchain}, $\PP(5,3,5,3,-1)$ does not exist, since
$$
  -1 + \frac{\inv(3,5)}{5} + \frac{\inv(3,5)}{5} = -1 + \frac 25 + \frac 25 < 0,
$$
so the inner type condition is not satisfied (Definition \ref{definnertype}).
\end{example}

\begin{example}
\label{ex23chain}
To follow up Example \ref{exHJ1}, denominators in the shortest path
\begin{equation}\label{exHJpath}
  \frac{19}{23} > \frac{14}{17} > \frac{9}{11} > \frac 45 > \frac 34 > \frac 57 > \frac{7}{10} > \frac{16}{23}
\end{equation}
form an inner sequence of weight 1
$$
  23,17,11,5,4,7,10,23 \quad \text{of inner type} \quad \textstyle 23\,\edge{17}{10}{-1}\,23.
$$
Conversely, from this inner type, Theorem \ref{thmchain}(2) constructs the shortest 1-path:
take a solution to $10 a'\equiv -1\bmod 23$, say $a'=16$. (Any other choice will shift the resulting 1-path by an integer.)
There is unique fraction $>\frac{16}{23}$ of the form $\frac{a}{23}$ with $a\equiv 17^{-1} \bmod 23$ and
$|[\frac{a'}{m'},\frac{a}{m}]\cap\Z|=0$, namely $\frac{19}{23}$. Then
Proposition \ref{fracmain} (4) reconstructs the sought 1-path it as \eqref{exHJpath}.
This is what Algorithm \ref{inneralg} does, phrased in terms of inner sequences.

Here is an example of a reduction type with
an inner chain between two principal components $\G_1, \G_2$
of type $23\,\edge{17}{10}{-1}\,23$, and outer chains
of type $23\frac{9}{}$, $23\frac{4}{}$ and $23\frac{3}{}$.
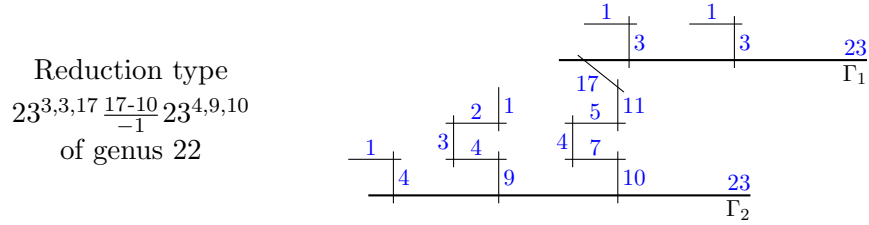
$$
\cbox{\begin{tabular}{cc}
Reduction type \\[3pt]
\redtype{23^{3,3,17}\e{17-10}(-1)23^{4,9,10}}
&\quad
\cr
of genus
22
\end{tabular}}
\qquad
\cbox{
\begin{tikzpicture}[xscale=1,yscale=0.8] \draw[thick,shorten >=2] (4/15,0)--(65/12,0) node[inner sep=2pt,blue,above left,scale=0.8] {23} node[inner sep=2pt,below left,scale=0.8] {$\Gamma_2$}; \draw[shorten <=-3pt, shorten >=-3pt] (3/5,0)--node[inner sep=2pt,scale=0.8,blue,right] {4} (3/5,3/5); \draw[shorten <=-3pt] (3/5,3/5)--node[inner sep=2pt,scale=0.8,blue,above] {1} (0,3/5); \draw[shorten <=-3pt, shorten >=-3pt] (2,0)--node[inner sep=2pt,scale=0.8,blue,right] {9} (2,3/5); \draw[shorten <=-3pt, shorten >=-3pt] (2,3/5)--node[inner sep=2pt,scale=0.8,blue,above] {4} (7/5,3/5); \draw[shorten <=-3pt, shorten >=-3pt] (7/5,3/5)--node[inner sep=2pt,scale=0.8,blue,left] {3} (7/5,6/5); \draw[shorten <=-3pt, shorten >=-3pt] (7/5,6/5)--node[inner sep=2pt,scale=0.8,blue,above] {2} (2,6/5); \draw[shorten <=-3pt] (2,6/5)--node[inner sep=2pt,scale=0.8,blue,right] {1} (2,9/5); \draw[shorten <=-3pt, shorten >=-3pt] (43/12,0)--node[inner sep=2pt,scale=0.8,blue,right] {10} (43/12,3/5); \draw[shorten <=-3pt, shorten >=-3pt] (43/12,3/5)--node[inner sep=2pt,scale=0.8,blue,above] {7} (179/60,3/5); \draw[shorten <=-3pt, shorten >=-3pt] (179/60,3/5)--node[inner sep=2pt,scale=0.8,blue,left] {4} (179/60,6/5); \draw[shorten <=-3pt, shorten >=-3pt] (179/60,6/5)--node[inner sep=2pt,scale=0.8,blue,above] {5} (43/12,6/5); \draw[shorten <=-3pt, shorten >=-3pt] (43/12,6/5)--node[inner sep=2pt,scale=0.8,blue,right] {11} (43/12,9/5); \draw[shorten <=-3pt, shorten >=-3pt] (43/12,9/5)--node[inner sep=2pt,scale=0.8,blue,below left=-1pt and -1pt] {17} (47/15,9/4); \draw[thick,shorten >=2] (14/5,9/4)--(209/30,9/4) node[inner sep=2pt,blue,above left,scale=0.8] {23} node[inner sep=2pt,below left,scale=0.8] {$\Gamma_1$}; \draw[shorten <=-3pt, shorten >=-3pt] (56/15,9/4)--node[inner sep=2pt,scale=0.8,blue,right] {3} (56/15,57/20); \draw[shorten <=-3pt] (56/15,57/20)--node[inner sep=2pt,scale=0.8,blue,above] {1} (47/15,57/20); \draw[shorten <=-3pt, shorten >=-3pt] (77/15,9/4)--node[inner sep=2pt,scale=0.8,blue,right] {3} (77/15,57/20); \draw[shorten <=-3pt] (77/15,57/20)--node[inner sep=2pt,scale=0.8,blue,above] {1} (68/15,57/20); \end{tikzpicture}
}
$$
\end{example}

\section{Total genus formula}
\label{stotalgenus}

\begin{notation}
\label{notmgOL}
Let $\G$ be a principal component of multiplicity $m=m_\G$ and geometric genus $g=g_\G$.
Decompose multiplicities of incident components $\{m_{\G'}\}_{\G'\in I(\G)}$ taken mod $m$
into two multisubsets of $\Z/m\Z$:

\begin{tabular}{lll}
  $\cO = \cO_\G$ & = & the multiset of initial multiplicities of outer chains from $\G$ taken mod $m$, \cr
  $\cL = \cL_\G$ & = & the multiset of initial multiplicities of inner chains from $\G$ taken mod $m$. \cr
\end{tabular}

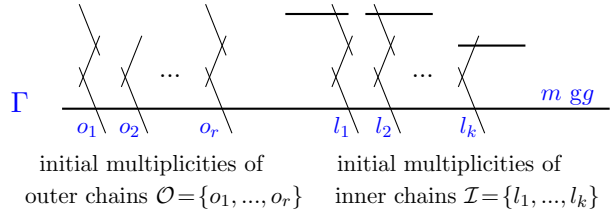
\bigskip
\begin{center}
\princomppicture
\end{center}

\end{notation}

We showed that chains of $\P^1$s are determined by their initial multiplicities and depths.
Therefore a reduction type can be reconstructed from

\begin{enumerate}
\item
Tuple of invariants $\mgOLG$ for each principal component $\G$;
\item
Inner chain endpoints ($=$ decomposition of $\coprod_\G \cL_\G$ into disjoint pairs);
\item
Depths of inner chains (one for every pair in (2)).
\end{enumerate}

\noindent
Varying depths we get what we call a reduction family, determined by (1)+(2):

\begin{definition}
Two reduction types belong to the same \textbf{reduction family} if one can be obtained from the other
by repeatedly replacing an inner chain by one of the same type except for its depth
($m\edge{l}{l'}{n}m'$ by $m\edge{l}{l'}{n'}m'$).
\end{definition}

We now show that the total Euler characteristic of the reduction type is the sum
$\sum_\G\chi_\G$, with a term $\chi_\G=\chi\mgOLG$ from each principal component $\G$. In particular,
it is completely independent of (2) and (3), so two reduction types in the same family have the same total genus.
This $\chi_\G$
will be the primary invariant for our classification.

\begin{definition}[Euler characteristic]
\label{defprineuler}
Let $\G$ be a principal component of multiplicity $m$, geometric genus $g$,
initial multiplicities of outer chains $\cO$ and initial multiplicities of outer chains $\cL$ (cf. \ref{notmgOL}).
We define the \textbf{Euler characteristic} of $\G$ by
$$
  \chi_\G = \chi\mgOL = (2\!-\!2g\!-\!|\cO|-\!|\cL|)m + \sum_{o\in\cO} \gcd(o,m).
$$
\end{definition}

\begin{theorem}[Total genus formula]
\label{eulerformula}
Let $R$ be a reduction type. Then
$$
  2-2g(R) = \sum_\G \chi_\G,
$$
the sum taken over all principal components of $R$.
\end{theorem}

\begin{proof}
Let us compute the right-hand side of the genus formula \eqref{adjunction}
$$
  2 g(R) - 2 = \sum_\G \Bigl(m_\G (2 g_\G -2) + \sum_{\G'\in I(\G)} m_{\G'}\Bigr).\eqno{(*)}
$$
Take a principal component $\G$ of multiplicity $m$, and let $\G=\G_0,\G_1,...,\G_r$ be an outer chain from~$\G$.
Write $d_i=m_{\G_i}$ (and $d_{r+1}=0$), so that $d_0=m$ and $d_1=o$ is the initial multiplicity of the chain.
The contribution from $\G_1,...,\G_r$ to $(*)$ is
$$
  \sum_{i=1}^r (d_{i-1}-2d_i+d_{i+1}).
$$
This is a telescopic series, adding up to
$$
  d_0-d_1-d_r = m-o-d_r,
$$
and $d_r=\gcd(o,m)$ by Corollary \ref{outermain}.
Similarly for an inner chain $\G=\G_0$, $\G_1$, ..., $\G_r$, $\G_{r+1}=\G'$ between $\G$ and $\G'$,
the contribution from $\G_1,...,\G_r$ to $(*)$ is
$$
  m-l + m'-l',
$$
where $l=d_1$ is the initial multiplicity of the chain, and $l'=d_r$.
Leaving $m'-l'$ to be accounted for by $\G'$, we find

\begin{center}
\begin{tabular}{ccccc}
contribution to $(*)$ from  & $+$ & contribution to $(*)$ from the $\G$-side & $+$ & contribution to $(*)$ \cr
outer chains from $\G$      &     & of the inner chains from $\G$            &     & from $\G$ itself\cr
\end{tabular}
\end{center}

$$
  = \Bigl(\sum_{o\in \cO(\G)}\!m\!-\!o\!-\!\gcd(o,m)\Bigr) + \Bigl(\sum_{l\in \cL(\G)}\!m\!-\!l\Bigr) +
    \Bigl(m (2 g_\G\!-\!2) + \sum_{\G'\in I(\G)} m_{\G'}\Bigr) = -\chi_\G,
$$
because $-o$'s and $-l$'s cancel with $\sum_{\G'\in I(\G)} m_{\G'}$, and the total number of $m$'s
is $|\cO(\G)|+|\cL(\G)|+2g_\G-2$.
Summing up over all principal components $\G$ gives $-\sum_\G \chi_\G$, as claimed.
\end{proof}

\section{Classification of cores}
\label{score}

We now prove that there are only finitely many tuples $\mgOL$ coming from principal components $\G$ with a given Euler
characteristic $\chi_\G$.
First, we look at tuples with $g=0$, $\cL=\emptyset$ and $\gcd(m,\cO)=1$, called `cores'.
We classify them, give bounds on $m$ and $|\cO|$ in terms of $\chi$ (\ref{coremain}),
and then deduce finiteness for principal components in general (\ref{cormgol}).

\begin{definition}
A \textbf{core} $\Psi=m^{o_1,o_2,...,o_k}$ is an integer $m\ge 1$ and unordered
elements $o_i\in \Z/m\Z$ with $o_i\ne m$, $\sum o_i\equiv 0\bmod m$ and $\gcd(m,o_1,...,o_k)=1$.
Its \textbf{Euler characteristic} is
$$
  \chi(\Psi) = m(2-k) + \sum_{i=1}^k \gcd (m,o_i).
$$
\end{definition}

\begin{notation}
The ten most frequently occurring cores will be written as follows:

\smallskip
\begin{center}
\begin{tabular}{l@{\>=\>}l@{\qquad}l@{\>=\>}l@{\qquad}l@{\>=\>}l@{\qquad}l@{\>=\>}l@{\qquad}l@{\>=\>}l@{\qquad}l@{\>=\>}l}
I     & $1^{\emptyset}$ &
D     & $2^{1,1}$     &
\IV   & $3^{1,1,1}$   &
\III  & $4^{1,1,2}$   &
\II   & $6^{1,2,3}$   \cr
\IZS  & $2^{1,1,1,1}$ &
T     & $3^{1,2}$     &
\IVS   & $3^{2,2,2}$   &
\IIIS  & $4^{3,3,2}$   &
\IIS   & $6^{5,4,3}$.  \cr
\end{tabular}
\end{center}
\end{notation}

\noindent
The cores I, D, T have $\chi=2$ and start an infinite family (\ref{coremain}(2), Table \ref{tabC2}).
The other seven are the cores with $\chi=0$ (\ref{coremain}(3), Table \ref{tabC0}), and their names
come from Kodaira types.

\smallskip

Note that $\chi_C\!=\!\chi(m,0,\cO,\emptyset)$ (cf. \ref{defprineuler}),
so a core is basically a component of genus $g=0$, no inner chains and $\gcd(m,\cO)=1$.
The tuple $(m,0,\cO,\emptyset)$ might not come from a principal component itself, but
every principal component has a core naturally associated to it:

\begin{definition}
\label{princore}
Let $\G$ be a principal component of multiplicity $m$,
and $o_1,...,o_k\in\Z/m\Z$ the multiplicities of incident components $\G'\in I(\G)$
that are non-zero in $\Z/m\Z$.
Extract $c\!=\!\gcd(m,o_1,...,o_k)$ from $m$ and the $o_i$, writing
$$
  m=c\bar m, \quad o_1=c\bar o_1, \>\>\ldots,\>\> o_k=c\bar o_k, \qquad \gcd(\bar m,\bar o_1,...,\bar o_k)=1.
$$
We call $c$ the \textbf{weight of $\G$} and ${\bar m}^{\bar o_1,\bar o_2,...,\bar o_k}$
the \textbf{core of $\G$}.
\end{definition}

\begin{example} [cf. Example \ref{ellcurvetableex}]
Kodaira types $\IZS, \IV, \IVS, \III, \IIIS, \II, \IIS$ each have a unique principal component $\G$,
whose core is denoted by the same Kodaira symbol. Type $\InS{n}$ has two principal components with core D,
type $\IZ=\redtype{I_{g1}}$ has one principal component with core I,
and type $\In{n}$ has no principal components.
\end{example}

\begin{theorem}
\label{coremain}
Let $\Psi=m^{o_1,o_2,...,o_k}$ be a core. Then
\begin{enumerate}
\item
$\chi(\Psi)\in \{2,0,-2,-4,-6,-8,...\}$.
\item
$\chi(\Psi)=2 \>\>\iff\>\> k\le 2\>\>\iff\>\> \Psi=1^\emptyset$ or $\Psi=m^{a,m-a}$ for $m\ge 2$ and $(a,m)=1$.
\item
$\chi(\Psi)=0 \>\>\iff\>\> \Psi\in\{ \IZS, \IV, \IVS, \III, \IIIS, \II{}, \IIS\}$.
\item
If $\chi(\Psi)\le 0$, then $2\le m\le 6-2\chi(\Psi)$. In fact,
$$
  m \le \frac{\alpha-2\chi(\Psi)}{k-2} \qquad\text{with}\ \
  \alpha=\Bigl\{
  \begin{tabular}{@{\>}l@{\ \>}l}
    8 & if $k$ is even,\\
    6 & if $k$ is odd.
  \end{tabular}
$$
with equality if and only if $\Psi$ is one of
$$
\begin{array}{l@{\quad}l}
6^{2,2,2,3,\ldots,3}, 6^{4,4,4,3,\ldots,3}  & \text{\rm ($k\ge 5$ odd)},\cr
30^{4,6,20,15,\ldots,15}, 30^{14,6,10,15,\ldots,15},
30^{16,24,20,15,\ldots,15}, 30^{26,24,10,15,\ldots,15} & \text{\rm ($k\ge 5$ odd)},\cr
\smash{m^{\frac{mk}2-a,a,\frac m2,\ldots,\frac m2}}\quad\text{with } m\equiv 2\bmod 4,\>\gcd(a,m)=2 & \text{\rm ($k\ge 3$ arbitrary).}\cr
\end{array}
$$
\item
$k\le 4-\chi(\Psi)$, with equality if and only if $k=2$ or $\Psi=2^{1,1,...,1}$ with $k$ even.
\item
For every $B\le 0$, there are finitely many cores $\Psi$ with $\chi(\Psi)=B$.
\end{enumerate}
\end{theorem}

\begin{proof}
Let us label the two conditions from the definition of the core:
$$
  \sum_{i=1}^k o_i=0\in\Z/m\Z\quad(\dagger), \qquad\qquad \gcd(m,o_1,...,o_k)=1\quad(\ast).
$$
Observe also that $o_i\ne 0$ in $\Z/m\Z$, so $\gcd(o_i,m)\le \frac m2$ for all $i$, and
\begin{equation}\label{chimkrel}
   \chi(\Psi) = m(2\!-\!k) + \sum_{i=1}^k \gcd(o_i,m) \le m(2\!-\!k) +
     \frac{km}2 = \frac m2(4\!-\!k).
\end{equation}

($\chi(\Psi)\in 2\Z$) Indeed, if $m$ is odd, then
$\chi(\Psi)\equiv - k + \sum_{i=1}^k 1\equiv 0\mod 2$.
If $m$ is even, then $m(2-k)$ is even. By $(\dagger)$, the number of odd $o_i$ is even,
whence the number of odd $\gcd(m,o_i)$ is even as well, so $\chi(\Psi)$ is again even.

(1,2,3,4,5 when $k\le 2$) The conditions $(\ast)$ and $(\dagger)$ imply that
$$
  \Psi=1^\emptyset\qquad\text{and}\qquad \Psi=m^{a,m-a}\text{\ \ for \ }m\ge 2,\> (a,m)=1
$$
are the only possible cores with $k\le 2$. They have $\chi=2$ and satisfy all the claims,
so from now on we assume $k\ge 3$ and $m\ge 2$.

(1,2,3 when $k\ge 4$) When $k\ge 4$, from \eqref{chimkrel} we get $\chi(\Psi)\le 0$, with equality
if and only if $\Psi=2^{1,1,1,1}=\IZS$. This proves (1), (2) and (3) when $k\ge4$.

(1,2,3 when $k=3$) Suppose $k=3$. Write
$$
  \chi(\Psi) = m (\frac{1}{m_1} \!+\! \frac{1}{m_2} \!+\! \frac{1}{m_3} \!-\! 1), \qquad m_i=\frac{m}{\gcd(o_i,m)}\ge 2,
$$
and observe that $m=\lcm(m_1,m_2,m_3)$ by $(\ast)$.
If some prime power $p^a$ divides $m_1$ (say),
equivalently the $p$-adic valuation $v_p(o_1)\le v_p(m)-a$,
then one of the $o_2,o_3$ must also satisfy $v_p(o_i)\le v_p(m)-a$ by the strong triangle
inequality, and so $p^a|m_i$. So

\smallskip
\hbox to\textwidth{\hfill every prime power that divides one of the $m_i$ divides at least two of them.\hfill($\ddagger$)}
\smallskip

It is well-known (see e.g. \cite[\S 4.2]{BMS}) that triples $m_1,m_2,m_3\ge 2$ with $\sum \frac{1}{m_i}\ge 1$ are,
up to reordering,
$$
\begin{array}{ll}
  (2,2,n),\>\>(2,3,3),\>\>(2,3,4),\>\>(2,3,5) & (\text{`Spherical case' }\sum \frac{1}{m_i}>1 \>\iff\> \chi(\Psi)>0) \cr
  (2,3,6),\>\>(2,4,4),\>\>(3,3,3)             & (\text{`Parabolic case' }\sum \frac{1}{m_i}=1 \,\iff\> \chi(\Psi)=0) \cr
\end{array}
$$
The only spherical one satisfying ($\ddagger$) is (2,2,2), but in this case $m=2$, all $o_i=1$ and
$o_1+o_2+o_3\not\equiv 0\mod m$. So $\chi(\Psi)\le 0$. The parabolic ones give two types each with $\chi(\Psi)=0$,
the last six in (3).

(4 when $k\ge 3$) We now assume $k\ge 3$, $m\ge 2$ and $\chi(\Psi)< 0$.
As above, write
\begin{equation}
\label{chireciprocals}
  \chi(\Psi) = m \Bigl(\frac{1}{m_1} + ... + \frac{1}{m_k} - (k\!-\!2)\Bigr), \qquad m_i=\frac{m}{\gcd(o_i,m)}\ge 2,
\end{equation}
with $m=\lcm(m_1,...,m_k)$, and observe that ($\ddagger$) holds.
We need to show that
$$
  m \le \frac{\alpha-2\chi(\Psi)}{k-2} \qquad\text{with}\ \
  \alpha=\Bigl\{
  \begin{tabular}{@{\>}l@{\ \>}l}
    8 & if $k$ is even,\\
    6 & if $k$ is odd.
  \end{tabular}
$$
First note that
$$
  m \le \frac{\alpha-2\chi(\Psi)}{k-2} \qquad\liff\qquad
  m(k-2)+2\chi(\Psi)\le \alpha,
$$
and consider a stronger condition (without $\alpha$)
$$
  m(k-2)+2\chi(\Psi)\le 0
    \>\>\>\iff\>\>\>
  -\chi(\Psi) \ge m\,\frac{k-2}2
    \>\>\>{\buildrel\eqref{chireciprocals}\over\iff}\>\>\>
  \sum_{i=1}^k \frac{1}{m_i} \le \frac{k-2}2
    \>\>\>\iff\>\>\>
  \sum_{i=1}^k \bigl(\frac12-\frac{1}{m_i}) \ge 1.
$$

$$
\begin{array}{lcccl}
(m_i) & k & m & m(k-2)+2\chi(\Psi) & m(k-2)+2\chi(\Psi)\le\alpha\cr
\hline
2,...,2,n,n  \text{ ($n$ odd)}   & \ge 3  & 2n   & 8 & \text{sharp for $k$ even, fails for $k$ odd}\cr
2,...,2,n,2n \text{ ($n$ odd)}   & \ge 3  & 2n   & 6 & \text{sharp for $k$ odd}\cr
2,...,2,n,n  \text{ ($n$ even)}  & \ge 3  & n    & 4 \cr
3, 3, 3                   & 3       & 3   & 3 \cr
2,...,2, 4, 4, 4          & \ge 3   & 4   & 2 \cr
5, 5, 5                   & 3       & 5   & 1 \cr
2,...,2, 3, 3, 6          & \ge 3   & 6   & 4 \cr
2,...,2, 3, 6, 6          & \ge 3   & 6   & 2 \cr
3, 9, 9                   & 3       & 9   & 1 \cr
2,...,2, 3, 4, 12         & \ge 3   & 12  & 4 \cr
3, 5, 15                  & 3       & 15  & 3 \cr
3, 7, 21                  & 3       & 21  & 1 \cr
2,...,2, 3, 5, 30         & \ge 3   & 30  & 4 \cr
3, 3, 3, 3                & 4       & 3   & 2 \cr
2,...,2, 3, 3, 3          & \ge 4   & 6   & 6 & \text{sharp for $k$ odd}\cr
2,...,2, 3, 3, 3, 6       & \ge 4   & 6   & 2 \cr
2,...,2, 5, 5, 5          & \ge 4   & 10  & 2 \cr
2,...,2, 3, 3, 4, 4       & \ge 4   & 12  & 4 \cr
3, 3, 5, 5                & 4       & 15  & 2 \cr
2,...,2, 3, 9, 9          & \ge 4   & 18  & 2 \cr
2,...,2, 3, 5, 15         & \ge 4   & 30  & 6 & \text{sharp for $k$ odd} \cr
2,...,2, 3, 7, 21         & \ge 4   & 42  & 2 \cr
3, 3, 3, 3, 3             & 5       & 3   & 1 \cr
2,...,2, 3, 3, 3, 3       & \ge 5   & 6   & 4 \cr
2,...,2, 3, 3, 5, 5       & \ge 5   & 30  & 4 \cr
2,...,2, 3, 3, 3, 3, 3    & \ge 6   & 6   & 2 \cr
\end{array}
$$
\begin{center}
  \tabletitle\label{tabMChi}
  Invariants of cores $\Psi=m^{o_1,...,o_k}$ with $m(k-2)+2\chi(\Psi)>0$.
\end{center}

\noindent
Integers $2\le m_1\le m_2\le\ldots\le m_k$ ($k\ge 0$) such that
\begin{itemize}
\item[(a)] every prime power $p^j>2$ that divides one of the $m_i$ divides at least two of them, and
\item[(b)] $\sum_{i=1}^k \bigl(\frac12-\frac{1}{m_i}) < 1$
\end{itemize}

\noindent
fall into finitely many families. Indeed,
\begin{itemize}
\item
If $k\le 2$, the tuples $(m_i)$ satisfying (a)+(b) are
$$
  \emptyset, \quad (2), \quad (n,n) \text{ ($n$ arbitrary)} \quad (n,2n) \text{ ($n$ odd)}.
$$
\item
Tuple $(m_1,...,m_k)$ satisfies (a)+(b) if and only $(2,m_1,...,m_k)$ does.
\item
At most 5 of the $m_i$ are $\ge 3$, since $\frac12-\frac{1}{m_i}\ge 1/6$ when $m_i\ge 3$.
\item
If at least 3 of the $m_i$ are $\ge 3$ then (b) together with
$(\frac12-\frac{1}{3})+(\frac12-\frac{1}{12})+(\frac12-\frac{1}{12})=1$ imply that if a prime power $p^a$
divides some $m_i$ then $p^a\le 11$. That is, all $m_i\,|\,2^3\cdot 3^2\cdot 5\cdot 7\cdot 11$.
\end{itemize}

This is a finite search, and the families it yields are listed in Table \ref{tabMChi}.

The claimed bound holds for all of them, and is sharp where specified, with
one exception: when $(m_i)=(2,...,2,n,n)$ with $n$ odd (first row) and $k$ is odd, $m(k-2)+2\chi(\Psi)=8$
is not $\le 6$ as claimed. However, a tuple $(o_i)$ with such $m_i=\frac{m}{\gcd(o_i,m)}$
is of the form ($n$,$n$,...,$n$,even,even) and cannot add up to 0 mod $2n$ when $k$ is odd.
So such $(o_i)$ simply do not exist in this case.

In all other cases where the bound is sharp, the $(o_i)$ do exist, and are as stated in the theorem.

(5) We have $k\ge 3$, $m\ge 2$, and therefore $\chi(\Psi)\le 0$ by (1) and (2). Then \eqref{chimkrel} gives
$$
  k \le 4-\frac{2\chi(\Psi)}{m} \le 4-\chi(\Psi).
$$
Equalities hold if and only if all $m=2$ and all $\gcd(o_i,m)=\frac m2$ in \eqref{chimkrel}.
In other words, $\Psi=2^{1,1,...,1}$. In that case, $k$ is even by
$(\dagger)$, and the asserted equality holds.

(6) If $\chi(\Psi)=B\le 0$ is fixed, $m$ and $k$ are bounded by (4) and (5).
As all $o_i\in\Z/m\Z$, there are finitely many choices for them, so the total number
of possible cores with $\chi(\Psi)=B$ is finite.
\end{proof}

\begin{remark}
Cores $m^{a,\frac m2-a,\frac m2}$ with maximal $m=6-2\chi$ in (4)
are realised by hyperelliptic curves $y^2=x^{2g+1}+\pi^b,\>\>ab\equiv 1\bmod 4g\!+\!2\>(=m)$, e.g.
$$
  \Psi/\Q_p: y^2=x^7+p^3 \qquad \qquad
  \cbox{
\begin{tikzpicture}[xscale=1,yscale=0.9] \draw[thick,shorten >=2] (4/15,0)--(121/30,0) node[inner sep=2pt,blue,above left,scale=0.8] {14} node[inner sep=2pt,below left,scale=0.8] {$\Gamma_1$}; \draw[shorten <=-3pt, shorten >=-3pt] (3/5,0)--node[inner sep=2pt,scale=0.8,blue,right] {5} (3/5,3/5); \draw[shorten <=-3pt] (3/5,3/5)--node[inner sep=2pt,scale=0.8,blue,above] {1} (0,3/5); \draw[shorten <=-3pt] (7/5,0)--node[inner sep=2pt,scale=0.8,blue,right] {2} (7/5,3/5); \draw[shorten <=-3pt] (11/5,0)--node[inner sep=2pt,scale=0.8,blue,right] {7} (11/5,3/5); \end{tikzpicture}
  }
$$
\begin{center}
\figuretitle Special fibre with maximal $m=14$ in genus 3.
\end{center}
These curves are $\Delta_v$-regular in every residue characteristic, and the special fibre
can be computed with \cite[Thm 1.1]{newton}. When $g=1$, these are elliptic curves
$y^2=x^3+p$ and $y^2=x^3+p^5$ that realise Kodaira types \II{} and \IIS{} with $m=6$.
\end{remark}

\forcecomment
\begin{notation}
\label{notprin2}
A component type $[d]\mgOL$, where $\mgOL$ is primitive with core $\Psi$ is denoted
$$
  [d]C_{{\rm g}g}\tfrac{l_1}{}\tfrac{l_2}{}\cdots \tfrac{l_{t}}{}==\cdots=
$$
$l_1,...,l_t$ are the non-zero elements of $d\cL$, and each `$=$' stands for a zero in $d\cL$.
The genus subscript is omitted if $g=0$, and $[d]$ is omitted when $d=1$.
\end{notation}

\begin{example}
See Example \ref{exellmrnc} and Tables.
\end{example}

\begin{theorem}
\label{compmain}
Let $\G=\mgOL$ be a principal component type. Then
\begin{enumerate}
\item
$\chi_\G \ge 0$.
\item
$\chi_\G=0$ if and only if $\G=[d]\G_0$ for some $d\ge 1$ with $\G_0$ one of
the 9 primitive component types (see also Table \ref{comps0} on page \pageref{comps0})
\begin{center}
  \rm \redtype D:;\qquad\redtype 1_{g1};\qquad\IZS\qquad\IV\qquad\IVS\qquad\III\qquad\IIIS\qquad\II{}\qquad\IIS
\end{center}
\item
If $\G$ is primitive,
then $m\le 2\chi_\G+6$, with equality if and only if $\G$ is of type
$$
  m^{a,\frac m2-a,\frac m2}, \qquad m\equiv 2\bmod 4,\>\gcd(a,m)=1.
$$
\item
For any $B>0$, there are finitely many component types $\G$ with $0<\chi_\G<B$.
\end{enumerate}
\end{theorem}

\begin{proof}

\emph{1. Primitivity.}
Because $\chi_{[d]\G}=d\chi_\G$, it is enough to prove (1)--(4) for primitive types.
So assume from now on that $\G$ is primitive.

\emph{2. Reduction to $|\cO|+|\cL|\ge 3$. }
If $|\cO|+|\cL|\le 2$, then $g>0$ (stability), and $\sum_{o\in\cO}o+\sum_{l\in\cL}l\equiv 0\mod m$ leaves
the following four possibilities:
$$
\begin{array}{|@{\quad}c@{\quad}|@{\quad}c@{\quad}|@{\quad}c@{\quad}|@{\quad}c@{\quad}|@{\quad}c@{\quad}|}
\hline
\G      & 1_{{\rm g}g} & 1_{{\rm g}g}{-} & 1_{{\rm g}g}{-}{-} & m_{{\rm g}g}^{a,-a}\>\>(m>1, \gcd(a,m)=1) \TBcr
\hline
\chi_\G & 2g\!-\!2\ge 0 & 2g\!-\!1>0 & 2g>0 & 2gm\!-\!2\ge 2m-2>0 \TBcr
\hline
\end{array}
$$
In every case, (1)--(4) hold for these, including the strict inequality $m<2\chi_\G+6$.
So assume from now on that $|\cO|+|\cL|\ge 3$ (in particular, stability is automatic).

\emph{3. Reduction of (1)--(3) to $|\cO|+|\cL|=3$. } Suppose $|\cO|+|\cL|\ge 4$. Then
$$
\begin{array}{llllll}
   \chi_\G &\ge& \displaystyle (|\cO|+|\cL|\!-\!2)m - \sum_{o\in\cO} \gcd(o,m) \ge (|\cO|+|\cL|\!-\!2)m - \frac{|\cO|}2 m \cr
    &=& \displaystyle \frac{m}{2}(|\cO|+2|\cL|-4) \ge \frac m2|\cL| \ge 0,\cr
\end{array}
$$
so (1) holds in this case. Moreover, $\chi_\G=0$ means that all inequalities are equalities,
so $g=0$, $\cL=\emptyset$, $|\cO|=4$, and all $\gcd(o,m)=\frac m2$. Then $\G$ is of type $2^{1,1,1,1}$,
proving (2).

Now we check $m<2\chi_\G+6$. When $g\!>\!0$ or $|\cL|\!>\!0$ or $|\cO|\!>\!4$, we have $\chi_\G\ge\frac m2$
by the same chain of inequalities as above.
Otherwise, write $\cO=\{o_1,o_2,o_3,o_4\}$, with $\cL=\emptyset$. Then
$$
  \chi_G = m (2 - \frac{1}{m_1} - ... - \frac{1}{m_4}), \qquad m_i=\frac{m}{\gcd(o_i,m)}\ge 2,
$$
and $m=\lcm(m_1,...,m_4)$ by primitivity. Observe that if some prime power $p^a$ divides $m_1$ (say),
equivalently the $p$-adic valuation $v_p(o_1)\le v_p(m)-a$,
then one of the $o_2,o_3,o_4$ must also satisfy $v_p(o_i)\le v_p(m)-a$ by the strong triangle inequality,
and so $p^a|m_i$. So

\smallskip
\hbox to\textwidth{($\dagger$) \hfill every prime power that divides one of the $m_i$ divides at least two of them.\hfill}
\smallskip

Now if the largest $p^a$ dividing two of the $m_i$ is $\ge 4$, then
$\chi_\G\ge m(2-1/2-1/2-1/4-1/4)\ge m/2$ again. Otherwise, all $m_i\in\{2,3\}$, and $m\le 6<2\chi_\G+6$.

\emph{4. Reduction of (1)--(3) to $|\cO|=3$, $\cL=\emptyset$, $g=0$.}
Suppose either $|\cO|\le 2$ or $g>0$. Then, similarly to Step 3,
$$
  \chi_\G \ge (2g-2+3)m - \frac{|\cO|}2 m = \frac m2(4g+2-|\cO|) \ge 0,
$$
with equalities if and only if $g=0$, $m=2$, $|\cL|=1$, $\cO=\{1,1\}$ (Type $2^{1,1}$=).
Also, if $g>0$, then $\chi_\G\ge 2m$.

Now suppose $\G$ has $g=0$ and $|\cO|\le 2$, and let $\G'$ be the type obtained from $\G$
by moving all inner multiplicities into outer ones; thus,
$m_{\G'}\!=\!m$, $g_{\G'}\!=\!g\!=\!0$, $\cO_{\G'}=\cO\cup\cL$ and $\cL_{\G'}=\emptyset$. Then
$$
  \chi_{\G} = \chi_{\G'} + \sum_{o\in \cL} \gcd(o,m) > \chi_G'.
$$
If $m<2\chi_{\G'}+6$, then the same holds for $\chi_\G$. So, we are reduced to the case $\cL=\emptyset$.

\emph{5. Proof of (1)--(3).}
From now on, $g=0$, $|\cO|=3$, $\cL=\emptyset$. Similarly to Step 3, write
$$
  \chi_G = m (1 \!-\! \frac{1}{m_1} \!-\! \frac{1}{m_2} \!-\! \frac{1}{m_3}), \qquad m_i=\frac{m}{\gcd(o_i,m)}\ge 2,\>m=\lcm(m_1,m_2,m_3),
$$
and the $m_i$ satisfy ($\dagger$).
It is well-known (see e.g. \cite[\S 4.2]{BMS}) that triples $m_1,m_2,m_3\ge 2$ with $\sum \frac{1}{m_i}\ge 1$ are,
up to reordering,
$$
\begin{array}{ll}
  (2,2,n),\>\>(2,3,3),\>\>(2,3,4),\>\>(2,3,5) & (\text{`Spherical case' }\sum \frac{1}{m_i}>1) \cr
  (2,3,6),\>\>(2,4,4),\>\>(3,3,3)             & (\text{`Parabolic case' }\sum \frac{1}{m_i}=1) \cr
\end{array}
$$
The only spherical one satisfying ($\dagger$) is (2,2,2), but in this case $m=2$, all $o_i=1$ and
$o_1+o_2+o_3\not\equiv 0\mod m$. So $\chi_\G\ge 0$. The parabolic ones give two types each with $\chi_\G=0$,
and these are the last six in (2).

In all other cases, $\chi_\G>0$, and it remains to show $m\le 2\chi_{\G}+6$ and establish when equality occurs.
Using ($\dagger$) and grouping
primes powers dividing each of the pairs of the $m_i$ and all three, we can write $m=n_1n_2n_3c$ with
$n_i$ pairwise coprime, and
$$
  m_1=cn_2n_3,\>\>m_2=cn_1n_3,\>\>m_3=cn_1n_2, \qquad n_1,n_2,n_3,c\ge 1.
$$
After permuting the $m_i$ if necessary, we may assume that $n_1\le n_2\le n_3$. Now, if
$n_1\ge 2$ then $3\le n_2<n_3$ as the $n_i$ are coprime, which gives $m_1,m_2,m_3\ge 6$
and $\chi_G\ge\frac m2$. So assume $n_1=1$, and note that either $n_2\ge 2$ or $c\ge 2$ (else $m_3=cn_1n_2=1$).
Now,
$$
  \chi_\G = cn_2n_3 (1-\frac{1}{cn_2n_3}-\frac{1}{cn_2}-\frac{1}{cn_3}) = cn_2n_3-1-n_2-n_3,
$$
and
$$
  2\chi_\G + 6 - m  = cn_2n_3 + 4 - 2 n_3 - 2 n_2 = (c-1)n_2n_3 + (n_2-2)(n_3-2).
$$
If $n_2\ge 2$, then this expression is $\ge 0$, with equality when $c=1$ and $n_2=2$. In this case,
$$
  2\nmid n_3, \qquad m=2n_3, \qquad (m_1,m_2,m_3) = (m,m/2,2),
$$
which gives the asserted type $m^{a,\frac m2-a,\frac m2}$ in (3). Finally, if $n_2=1$ and $c\ge 2$, then
$$
  2\chi_\G + 4 - m  = cn_3 + 2 - 2 n_3 - 2 = (c-2)n_3 \ge 0,
$$
proving that $m\le 2\chi_\G+4$ in this case.

\emph{6. Proof of (4).}
It remains to show that for every $B>0$ there are only finitely primitive types $\G$ with $0<\chi_\G<B$.
By (3), $m\le 2B+4$ is bounded with $B$. Also, as in Step 3, we find
$$
  B \ge \chi_\G \ge \Bigl(2g-2+|\cL|+\frac{|\cO|}2\Bigr)m,
$$
so that $g$, $|\cL|$ and $|\cO|$ are also bounded with $B$. As all elements of $\cO$ and $\cL$ are
in $\Z/m\Z$, we get the asserted finiteness.
\end{proof}
\endcomment

\section{Invariants of principal components}
\label{sprininv}

We now consider principal components $\G$, and their primary invariants:
multiplicity $m$, genus $g$, initial multiplicities $\cO$ of outer chains,
and $\cL$ of inner chains from $\G$ (see \ref{defprineuler}, \ref{princore}).
For this section, it is convenient to define a \textbf{multiple} $[c]\G$ ($c\ge 1$) as a component with
all multiplicties multiplied by $c$, that is with invariants $(cm, g, c\cO, c\cL)$.
Note that $\chi([c]\Gamma)=c\chi(\Gamma)$.

From the finiteness theorem for cores (\ref{coremain}) we deduce the corresponding statement for
tuples $\mgOL$ with fixed $\chi\le 0$ (\ref{cormgol}) (which together with the classfication of chains
shows that there are finitely many reduction families in a fixed genus, see \ref{thmfamilies}). We go a bit further
and classify tuples with $m$ `large' relative to $-\chi$. Possible values $\frac{m}{-\chi}>2$ turn out to be
discrete, with 2 as the largest accumulation point, and we classify all such tuples $\mgOL$ in \ref{thmlargem}.

We also deduce that the largest possible component multiplicity in a reduction type of genus $g\ge 2$ is $12(g-1)$,
and find the corresponding families (\ref{corlargem}, \ref{remlargem}).

\newpage

\begin{theorem}
\label{cormgol}
Let $\G$ be a principal component with invariants $\mgOL$. Then
\begin{enumerate}
\item
$\chi_\G\le 0$.
\item
$\chi_\G=0$ if and only if either

\begin{tabular}{cl}
{\rm (D-tail)} & $g\!=\!0$, $\core(\G)={\rm D}$ and $\cL\!=\!\{m\}$, \cr
{\rm ($\rm I_{\rm g1}$ case)} & $g\!=\!1$ and $\cO\!=\!\cL\!=\!\{\}$, \cr
{\rm (all other)} & $g\!=\!0$, $\core(\G)\in\{\IZS, \IV, \IVS, \III, \IIIS, \II, \IIS\}$ and $\cL\!=\!\{\}$.\cr
\end{tabular}

\noindent
In other words, $\G$ is a multiple of one of the following (outer multiplicities $o\in\cO$ are shown pointing
down, inner multiplicities $l\in\cL$ pointing up):

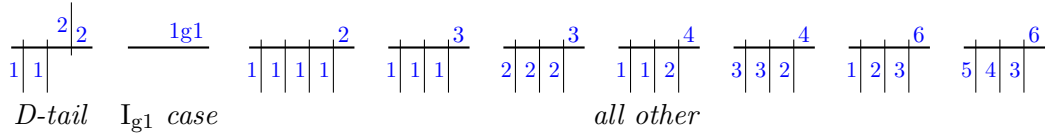
\begin{center}
\begin{tabular}{ccccccccc}
\begin{tikzpicture}[xscale=0.4,yscale=1,baseline=0] \draw[thick] (0,0)--(2.7,0) node[blue,scale=0.8,inner sep=2pt,above left] {$2$} node[inner sep=2pt,below left=0pt and -1pt,scale=0.8] {$$}; \draw[shorten <=-3pt] (0.4,0)--node[blue,scale=0.77,inner sep=1.5pt,left] {$1$} (0.4,-0.6); \draw[shorten <=-3pt] (1.2,0)--node[blue,scale=0.77,inner sep=1.5pt,left] {$1$} (1.2,-0.6); \draw[shorten <=-3pt] (2,0)--node[blue,scale=0.77,inner sep=1.5pt,left] {$2$} (2,0.6); \end{tikzpicture}&\begin{tikzpicture}[xscale=0.4,yscale=1,baseline=0] \draw[thick] (0,0)--(2.7,0) node[blue,scale=0.8,inner sep=2pt,above left] {$1{\rm{} g}1$} node[inner sep=2pt,below left=0pt and -1pt,scale=0.8] {$$}; \end{tikzpicture}&\begin{tikzpicture}[xscale=0.4,yscale=1,baseline=0] \draw[thick] (0,0)--(3.5,0) node[blue,scale=0.8,inner sep=2pt,above left] {$2$} node[inner sep=2pt,below left=0pt and -1pt,scale=0.8] {$$}; \draw[shorten <=-3pt] (0.4,0)--node[blue,scale=0.77,inner sep=1.5pt,left] {$1$} (0.4,-0.6); \draw[shorten <=-3pt] (1.2,0)--node[blue,scale=0.77,inner sep=1.5pt,left] {$1$} (1.2,-0.6); \draw[shorten <=-3pt] (2,0)--node[blue,scale=0.77,inner sep=1.5pt,left] {$1$} (2,-0.6); \draw[shorten <=-3pt] (2.8,0)--node[blue,scale=0.77,inner sep=1.5pt,left] {$1$} (2.8,-0.6); \end{tikzpicture}&\begin{tikzpicture}[xscale=0.4,yscale=1,baseline=0] \draw[thick] (0,0)--(2.7,0) node[blue,scale=0.8,inner sep=2pt,above left] {$3$} node[inner sep=2pt,below left=0pt and -1pt,scale=0.8] {$$}; \draw[shorten <=-3pt] (0.4,0)--node[blue,scale=0.77,inner sep=1.5pt,left] {$1$} (0.4,-0.6); \draw[shorten <=-3pt] (1.2,0)--node[blue,scale=0.77,inner sep=1.5pt,left] {$1$} (1.2,-0.6); \draw[shorten <=-3pt] (2,0)--node[blue,scale=0.77,inner sep=1.5pt,left] {$1$} (2,-0.6); \end{tikzpicture}&\begin{tikzpicture}[xscale=0.4,yscale=1,baseline=0] \draw[thick] (0,0)--(2.7,0) node[blue,scale=0.8,inner sep=2pt,above left] {$3$} node[inner sep=2pt,below left=0pt and -1pt,scale=0.8] {$$}; \draw[shorten <=-3pt] (0.4,0)--node[blue,scale=0.77,inner sep=1.5pt,left] {$2$} (0.4,-0.6); \draw[shorten <=-3pt] (1.2,0)--node[blue,scale=0.77,inner sep=1.5pt,left] {$2$} (1.2,-0.6); \draw[shorten <=-3pt] (2,0)--node[blue,scale=0.77,inner sep=1.5pt,left] {$2$} (2,-0.6); \end{tikzpicture}&\begin{tikzpicture}[xscale=0.4,yscale=1,baseline=0] \draw[thick] (0,0)--(2.7,0) node[blue,scale=0.8,inner sep=2pt,above left] {$4$} node[inner sep=2pt,below left=0pt and -1pt,scale=0.8] {$$}; \draw[shorten <=-3pt] (0.4,0)--node[blue,scale=0.77,inner sep=1.5pt,left] {$1$} (0.4,-0.6); \draw[shorten <=-3pt] (1.2,0)--node[blue,scale=0.77,inner sep=1.5pt,left] {$1$} (1.2,-0.6); \draw[shorten <=-3pt] (2,0)--node[blue,scale=0.77,inner sep=1.5pt,left] {$2$} (2,-0.6); \end{tikzpicture}&\begin{tikzpicture}[xscale=0.4,yscale=1,baseline=0] \draw[thick] (0,0)--(2.7,0) node[blue,scale=0.8,inner sep=2pt,above left] {$4$} node[inner sep=2pt,below left=0pt and -1pt,scale=0.8] {$$}; \draw[shorten <=-3pt] (0.4,0)--node[blue,scale=0.77,inner sep=1.5pt,left] {$3$} (0.4,-0.6); \draw[shorten <=-3pt] (1.2,0)--node[blue,scale=0.77,inner sep=1.5pt,left] {$3$} (1.2,-0.6); \draw[shorten <=-3pt] (2,0)--node[blue,scale=0.77,inner sep=1.5pt,left] {$2$} (2,-0.6); \end{tikzpicture}&\begin{tikzpicture}[xscale=0.4,yscale=1,baseline=0] \draw[thick] (0,0)--(2.7,0) node[blue,scale=0.8,inner sep=2pt,above left] {$6$} node[inner sep=2pt,below left=0pt and -1pt,scale=0.8] {$$}; \draw[shorten <=-3pt] (0.4,0)--node[blue,scale=0.77,inner sep=1.5pt,left] {$1$} (0.4,-0.6); \draw[shorten <=-3pt] (1.2,0)--node[blue,scale=0.77,inner sep=1.5pt,left] {$2$} (1.2,-0.6); \draw[shorten <=-3pt] (2,0)--node[blue,scale=0.77,inner sep=1.5pt,left] {$3$} (2,-0.6); \end{tikzpicture}&\begin{tikzpicture}[xscale=0.4,yscale=1,baseline=0] \draw[thick] (0,0)--(2.7,0) node[blue,scale=0.8,inner sep=2pt,above left] {$6$} node[inner sep=2pt,below left=0pt and -1pt,scale=0.8] {$$}; \draw[shorten <=-3pt] (0.4,0)--node[blue,scale=0.77,inner sep=1.5pt,left] {$5$} (0.4,-0.6); \draw[shorten <=-3pt] (1.2,0)--node[blue,scale=0.77,inner sep=1.5pt,left] {$4$} (1.2,-0.6); \draw[shorten <=-3pt] (2,0)--node[blue,scale=0.77,inner sep=1.5pt,left] {$3$} (2,-0.6); \end{tikzpicture}
\\[16pt]
D-tail & $\rm I_{\rm g1}$ case & \multicolumn{7}{c}{all other} \cr
\end{tabular}
\end{center}

\item
For every $B<0$, there are finitely possible tuples of invariants $\mgOL$
of principal components $\G$ with $\chi_\G=B$.
\end{enumerate}
\end{theorem}

\begin{proof}
It is enough to prove the statements when $\gcd(m,\cO,\cL)=1$, since $\chi([c]\Gamma)=c\chi(\Gamma)$ for $c>1$.
Let $o_1,...,o_k$ be the nonzero elements of $\cO\cup\cL$, with multiplicities, and $n$ the number of zeroes in $\cL$.
(There are no zeroes in $\cO$ by \ref{outermain}.)
Then $\Psi=m^{o_1,...,o_k}$ is the core of $\G$, and
$$
  \chi_\G \le \chi(\Psi) - m(n+2g),
$$
with equality if and only if $\cL=\emptyset$.

Because $\G$ is principal (minimality condition in \ref{defredtype}), by Theorem \ref{coremain} either
\begin{itemize}
\item[(a)]
$\chi(\Psi)\le 0$, or
\item[(b)]
$\Psi=\redtype{I}$ and either $g\ge 1$ or $n\ge 3$, or
\item[(c)]
$\Psi=m^{a,m-a}$ with $m\ge 2$, $\gcd(a,m)=1$ and either $g\ge 1$ or $n\ge 1$.
\end{itemize}
Then (a) $\chi_\G\le \chi(\Psi)\le 0$ or (b) $\chi_\G\le 2-2g-n\le 0$ or
(c) $\chi_\G\le 2-m(n+2g)\le 0$.
So $\chi_\G\le 0$, and also there are finitely many $\G$ with a given $\chi=B<0$
in all three cases.

Finally, the nine possibilities with $\chi_\G=0$ come from
the core I with $g=1$, core D with $n=1$, and seven cores $\IZS, \IV, \IVS, \III, \IIIS, \II{}, \IIS$
with $\chi(\Psi)=0$, see \ref{coremain} (2, 3).
\end{proof}

\begin{theorem}
\label{thmlargem}
Let $\G$ be a principal component with invariants $\mgOL$ and $\chi=\chi_\G$. If $$\frac{m}{-\chi}>2$$ then
$\frac{m}{-\chi}$ is of the form $2+\frac 4k$ or $2+\frac 3k$ for some $k\ge 1$, and every such fraction\footnote
{Explicitly, these are fractions $
  6, 5, 4, 3\frac{1}{2}, 3\frac{1}{3}, 3, 2\frac{4}{5}, 2\frac{3}{4}, 2\frac{2}{3}, 2\frac{3}{5}, 2\frac{4}{7},
  2\frac{1}{2}, 2\frac{4}{9}, 2\frac{3}{7}, 2\frac{2}{5}, 2\frac{3}{8}, 2\frac{4}{11}, 2\frac{1}{3}, 2\frac{4}{13},
  2\frac{3}{10}, \ldots
$}
occurs in this way.
See Table \ref{tablargem} (p. \pageref{tablargem}) for a complete list of such tuples $\mgOL$.
\end{theorem}

\begin{proof}
From the proof of Theorem \ref{cormgol} it follows that principal components
with $\frac{m}{-\chi}>2$ have cores that satisfy the same inequality.
Such cores are a subset of those listed in Table \ref{tabMChi}, and analysing
them we find all the possible tuples $\mgOL$.
\end{proof}

\begin{corollary}[Large multiplicities]
\label{corlargem}
Let $R$ be a reduction type of total genus $g>1$. If a principal component $\G$ of $R$ has
multiplicity $m>2(2g-2)$, it is listed in Table \ref{tablargem}.
\end{corollary}

\begin{proof}
This follows directly from the theorem, since $0\le -\chi_\G\le 2g-2$ by \ref{cormgol} and \ref{eulerformula}.
\end{proof}

\begin{remark}
\label{remlargem}
For example, there are exactly two families in every genus $g\ge 2$ that have a component of the
multiplicity $6(2g-2)$, which is largest possible.
They are obtained by taking $\G_0$ from the top row of Table \ref{tablargem}
with $k=1$ ($m=6$, $\frac{m}{-\chi}=6$) and $u\in(\Z/6\Z)^\times$, letting $\G=[c]\G_0$ with $c=2g-2$,
and attaching a D-tail of some depth $n$:

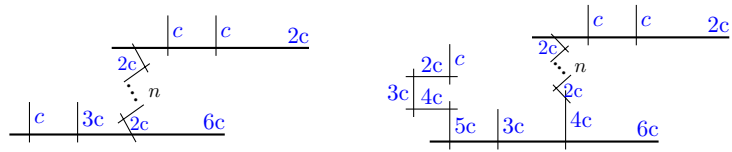
\begin{center}
\begin{tikzpicture}[xscale=0.8,yscale=0.72] \draw[thick,shorten >=2] (-1/3,0)--(10/3,0) node[inner sep=2pt,blue,above left,scale=0.8] {6c}; \draw[shorten <=-3pt] (0,0)--node[inner sep=2pt,scale=0.8,blue,right] {$c$} (0,3/5); \draw[shorten <=-3pt] (4/5,0)--node[inner sep=2pt,scale=0.8,blue,right] {3c} (4/5,3/5); \draw (17/10,0) ++(0.07,-0.14)--++(-0.26,0.53)++(0.12,-0.2) ++(0,0.02) node[blue,scale=0.7,inner sep=1,anchor=west]{2c} ++(0,-0.02) ++(-0.16,0.02)--++(0.43,0.34)++(-0.13,0.07) node{$\cdot$} ++(-0.06,0.13) node{$\cdot$} ++(-0.06,0.13) node{$\cdot$}++(0.17,-0.1)    node[auto,inner sep=5,scale=0.7,anchor=west]{$n$} ++(-0.25,0.19)--++(0.43,0.34)++(-0.18,0.04) ++(0,-0.04) node[blue,scale=0.7,inner sep=1,anchor=east]{2c} ++(0,0.04) ++(0.18,-0.04)++(-0.03,-0.16)--++(-0.26,0.53); \draw[thick,shorten >=2] (41/30,8/5)--(71/15,8/5) node[inner sep=2pt,blue,above left,scale=0.8] {2c}; \draw[shorten <=-3pt] (23/10,8/5)--node[inner sep=2pt,scale=0.8,blue,right] {$c$} (23/10,11/5); \draw[shorten <=-3pt] (31/10,8/5)--node[inner sep=2pt,scale=0.8,blue,right] {$c$} (31/10,11/5); \end{tikzpicture}
\qquad
\begin{tikzpicture}[xscale=0.8,yscale=0.72] \draw[thick,shorten >=2] (4/15,0)--(499/120,0) node[inner sep=2pt,blue,above left,scale=0.8] {6c}; \draw[shorten <=-3pt, shorten >=-3pt] (3/5,0)--node[inner sep=2pt,scale=0.8,blue,right] {5c} (3/5,3/5); \draw[shorten <=-3pt, shorten >=-3pt] (3/5,3/5)--node[inner sep=2pt,scale=0.8,blue,above] {4c} (0,3/5); \draw[shorten <=-3pt, shorten >=-3pt] (0,3/5)--node[inner sep=2pt,scale=0.8,blue,left] {3c} (0,6/5); \draw[shorten <=-3pt, shorten >=-3pt] (0,6/5)--node[inner sep=2pt,scale=0.8,blue,above] {2c} (3/5,6/5); \draw[shorten <=-3pt] (3/5,6/5)--node[inner sep=2pt,scale=0.8,blue,right] {$c$} (3/5,9/5); \draw[shorten <=-3pt] (7/5,0)--node[inner sep=2pt,scale=0.8,blue,right] {3c} (7/5,3/5); \draw[shorten <=-3pt, shorten >=-3pt] (101/40,0)--node[inner sep=2pt,scale=0.8,blue,right] {4c} (101/40,4/5); \draw (101/40,4/5) ++(0.08,-0.09)--++(-0.29,0.34)++(0.14,-0.12) ++(0,0.01) node[blue,scale=0.7,inner sep=1,anchor=west]{2c} ++(0,-0.01) ++(-0.14,0)--++(0.3,0.3)++(-0.12,0.03) node{$\cdot$} ++(-0.07,0.08) node{$\cdot$} ++(-0.07,0.08) node{$\cdot$}++(0.16,-0.05)    node[auto,inner sep=5,scale=0.7,anchor=west]{$n$} ++(-0.24,0.1)--++(0.3,0.3)++(-0.16,0) ++(0.01,-0.03) node[blue,scale=0.7,inner sep=1,anchor=east]{2c} ++(-0.01,0.03) ++(0.16,0)++(0,-0.12)--++(-0.29,0.34); \draw[thick,shorten >=2] (59/30,77/40)--(16/3,77/40) node[inner sep=2pt,blue,above left,scale=0.8] {2c}; \draw[shorten <=-3pt] (29/10,77/40)--node[inner sep=2pt,scale=0.8,blue,right] {$c$} (29/10,101/40); \draw[shorten <=-3pt] (37/10,77/40)--node[inner sep=2pt,scale=0.8,blue,right] {$c$} (37/10,101/40); \end{tikzpicture}
\end{center}

\noindent
These families are denoted \redtype{[c]II_D} and \redtype{[c]II^*_D}, respectively, see \S\ref{slabel}.

Similarly, we find (again in the notation of \S\ref{slabel}):

\noindent
In genus $g=2$, reduction families with a component of multiplicity $m>3(2g-2)=6$ are
\begin{center}
\redtype{[2]II^*_D}, \redtype{[2]II_D}, \redtype{10^{1,4,5}}, \redtype{10^{3,2,5}}, \redtype{10^{7,8,5}}, \redtype{10^{9,6,5}}, \redtype{8^{1,3,4}}, \redtype{8^{5,7,4}}, \redtype{[2]III^*_D}, \redtype{[2]III_D}.
\end{center}

\noindent
In genus $g=3$, reduction families with a component of multiplicity $m>3(2g-2)=12$ are
\begin{center}
\begin{minipage}{13.4cm}
\redtype{[4]II^*_D}, \redtype{[4]II_D}, \redtype{[2]10^{1,4,5}}, \redtype{[2]10^{3,2,5}}, \redtype{[2]10^{7,8,5}}, \redtype{[2]10^{9,6,5}}, \redtype{[3]II^*\e T}, \redtype{[3]II\e T}, \redtype{[2]8^{1,3,4}}, \redtype{[2]8^{5,7,4}}, \redtype{[4]III^*_D}, \redtype{[4]III_D}, \redtype{14^{1,6,7}}, \redtype{14^{11,10,7}}, \redtype{14^{13,8,7}}, \redtype{14^{3,4,7}}, \redtype{14^{5,2,7}}, \redtype{14^{9,12,7}}.
\end{minipage}
\end{center}

\noindent
In genus $g=4$, reduction families with a component of multiplicity $m>3(2g-2)=18$ are
\begin{center}
\begin{minipage}{15.3cm}
\redtype{[6]II^*_D}, \redtype{[6]II_D}, \redtype{[3]10^{1,4,5}}, \redtype{[3]10^{3,2,5}}, \redtype{[3]10^{7,8,5}}, \redtype{[3]10^{9,6,5}}, \redtype{[3]8^{1,3,4}}, \redtype{[3]8^{5,7,4}}, \redtype{[4]II^*\e 4^{1,3}}, \redtype{[4]II^*\e [2]D_D}, \redtype{[4]II\e 4^{1,3}}, \redtype{[4]II\e [2]D_D}, \redtype{[6]III^*_D}, \redtype{[6]III_D}, \redtype{[2]10^{1,4,5}_D}, \redtype{[2]10^{3,2,5}_D}, \redtype{[2]10^{7,8,5}_D}, \redtype{[2]10^{9,6,5}_D}.
\end{minipage}
\end{center}
\end{remark}

\begin{corollary}[Bounds]
\label{prinbounds}
For a principal component $\G$ with invariants $\mgOL$ and $\chi_\G<0$,
the following bounds hold:
\begin{enumerate}
\item
$m\le -6\chi_\G$,
\item
$g\le \lfloor 1-\chi_\G/2\rfloor$,
\item
$m(|\cO|+2|\cL|-4)\le -2\chi_\G$,
\item
$|\cO|\le 4-\chi_\G$.
\item
$|\cL|\le 2-\frac{\chi_\G}{2}$, unless $m=1$, $\cO=\emptyset$ and $|\cL|=2\!-\!\chi_\G\!-\!2g\le 2\!-\!\chi_\G$.
\end{enumerate}
\end{corollary}

\begin{proof}
(1) Immediate from \ref{thmlargem}, since the largest fraction $m/-\chi_\G$ is 6.

\noindent
(3) From the definition of $\chi$ (\ref{defprineuler}) and noting that each $\gcd(o,m)\le m/2$, we get
$$
  \chi =(2\!-\!2g\!-\!|\cO|-\!|\cL|)m + \sum_{o\in\cO} \gcd(o,m) \le
  (2-2g-\frac{|\cO|}{2}-|\cL|)m,
$$
and a similar computation proves (2) as well (or apply \ref{corabtor} below).

\noindent
(4) This follows from (3) when $m\ge 2$. Otherwise, $m=1$, $\cO=\emptyset$ and the claim holds as well.

\noindent
(5) When $|\cL|>\frac{-\chi_\G}{2}+2$, from (3) we get $m=1$, and in this case $\cO=\emptyset$.
\end{proof}

\begin{remark}
The bounds on $m$, $g$, $|\cO|$, $|\cL|$ in \ref{prinbounds} are sharp.
For example, reduction families of the following form achieve them (in the notation of \S\ref{slabel}):
\begin{center}
$[c]$\redtype{II_D} (maximal $m$), \quad
\redtype{I_{g\text{$g$}}} (maximal $g$), \quad
\redtype{2^{1,...,1}} (maximal $|\cO|$), \quad
I$\frac{1}{}\cdots\frac{1}{}$ (maximal $|\cL|$).
\end{center}
\end{remark}

\section{Finiteness and classification in genus 1}
\label{sfamilies}

We recover the classical classification of families in genus 1 and finiteness in genus $>1$:

\begin{corollary}
\noindent\par\noindent
\begin{enumerate}
\item (Kodaira--N\'eron classification \cite{Ko})
Reduction families of elliptic curves are the ten Kodaira-N\'eron families
$\redtype{I_{g1}}(=\IZ)$, $\In{{\sl n}}$, $\IZS$, $\InS{{\sl n}}$, \IV, \IVS, \III, \IIIS, \II{} and \IIS{} (p.\pageref{kodairatypes}).
\item (Saito \cite{Saito})
Reduction types of genus 1 curves are multiples \redtype{[\hbox{$c$}]I_{g1}}, \redtype{[\hbox{$c$}]I_{{\hbox{$\scriptstyle n$}}}}, \ldots, \redtype{[\hbox{$c$}]II^*}
of Kodaira types, obtained by multiplying all multiplicities by some $c\ge 1$.
\item (Artin--Winters \cite{AW})
There are finitely many reduction families in fixed genus $g\ge 2$.
\end{enumerate}
\label{thmfamilies}
\end{corollary}

\begin{proof}
Let $R$ be a reduction type of genus $g$. Recall that principal components $\G$ have $\chi_\G\le 0$ (Corollary \ref{cormgol}(1)),
and $\sum_\G\chi_\G=2-2g$ (Theorem \ref{eulerformula}).

(2)
As $2-2g=0$, every principal component $\G$ must have $\chi_\G=0$, and these are classified in \ref{cormgol}(2).
The eight with $\cL=0$ give eight types $\In0=\redtype{I_{g1}}$, \IZS, \IV, \IVS, \III, \IIIS, \II, \IIS{}
and their multiples. Apart from the exceptional family $\In{n}$ and its multiples, the only other possibility
is to glue two core D types to a multiple of $\InS{n}$.

(1)
An elliptic curve has a rational point at infinity, which reduces to a component of multiplicity 1.
By (2), its reduction type is in one of Kodaira--N\'eron families.

(3)
Take a reduction family in genus $g\ge 2$, and consider its principal components $\G$ with $\chi_\G<0$.
There are at most $2g-2$ of them by the total genus formula, and finitely many possibilities for its type
by Corollary \ref{cormgol} (3). Because a chain between components is determined (up to depth) by its outgoing multiplicites,
there are finitely many possible ways to construct inner chains between these components. Finally, all unused outgoing
multiplicities have to become loops (from a component to itself) or D-tails (from a component to one with $\chi=0$ as in
\ref{cormgol} (2)), and again there are finitely many options for these, up to depths.
\end{proof}

\begin{corollary}
\label{corabtor}
Let $R$ be a reduction type. Write
\begin{center}
\begin{tabular}{llllll}
$V$ &=& set of principal components of $R$, \cr
$E$ &=& set of inner chains of $R$, \cr
$a$ &=& $\sum_{\G\in V} g_\G$ & (`abelian part'),\cr
$t$ &=& $|E|-|V|+1$ & (`toric part').\cr
\end{tabular}
\end{center}
Then
$$
  g(R) \ge a + t,
$$ 
with equality if and only if either
\begin{itemize}
\item
Every component of $R$ has multiplicity 1 (i.e. $R$ is semistable, cf. \ref{defsemistable}); or
\item
$R=\redtype{[\hbox{$c$}]I_{{\hbox{$\scriptstyle n$}}}}$ or $\redtype{[\hbox{$c$}]I^*_{{\hbox{$\scriptstyle n$}}}}$,
of total genus 1.
\end{itemize}
\end{corollary}

\begin{proof}
When $g(R)=1$, the claim is true by inspection, from the classification in \ref{cormgol} (2), (3).
Suppose $g(R)>1$, and
consider a (multi)graph $G$ with vertex set $V$ and edge set $E$. Then $t$ is the dimension
of its first homology group, the `number of loops in $G$'. We have
$$
  2g(R)-2
    \>\>{\buildrel\ref{eulerformula}\over=}\>\>
    \sum_{\G\in V} m_\G\Bigl(2g_\G-2+|\cL_\G|+\sum_{o\in \cO(\G)} (1-\frac{\gcd(o,m_\G)}{m_\G}) \Bigr)
    \>\>{\buildrel\ref{cormgol}(1)\over\ge}\>\>
    \sum_{\G\in V} (2g_\G-2+|\cL_\G|)
$$
$$
  = (2\sum_{\G\in V} g_\G) -2|V|+2|E| = 2a+2t-2,
$$
proving the inequality. If equality holds, then $\cO_\G=\emptyset$ for every $\G\in V$.
In particular there are no D-tails, and none of the cases in \ref{cormgol} (2) are possible
since $g(R)>1$. So $\chi_G<0$ for every $\G\in V$ and the equality can only hold
when $m_\G=1$ for every $\G\in V$. By the structure of chains this implies that
all components of $R$ have multiplicity 1. Conversely, if all multiplicities are 1
(which automatically implies $\cO_\G=\emptyset$ for all $\G\in V$), the equality does hold.
\end{proof}

\section{Principal types and shapes}
\label{sshapes}

Frow now on, suppose $R$ is a reduction type of genus $>1$. Thus, $\chi(R)<0$
and $R$ has principal components $\G$ with $\chi_\G<0$.
We break $R$ into pieces centered at such components, incorporating chains of $\P^1$s from $\G$ to itself (`loops')
and from $\G$ to components with $\chi=0$ (`D-tails') into the data for $\G$, and call these `principal types'.
These are connected by inner chains~(`edges').

\begin{definition}[Principal type]
\label{defcomfam}
Let $\G$ be a principal component with $\chi_\G<0$.
\begin{itemize}
\item
A \textbf{loop} is an inner chain from $\G$ to itself.
\item
A \textbf{D-tail} is an inner chain from $\G$ to a principal component with $\chi_{\G'}=0$.\\
(By Corollary \ref{cormgol} (2) these look like in the picture below, with $2n_i=\gcd(m,l''_i)$.)
\item
An \textbf{edge} is an inner chain from $\G$ to a principal component $\G'\ne \G$ with $\chi_{\G'}<0$.
\end{itemize}
Say $\G$ has multiplicity~$m$, genus~$g$, initial outer chain multiplicities~$\cO$, and
initial inner chain multiplicities~$\cL$ (\ref{defprineuler}). This trichotomy splits $\cL$ into three disjoint parts:
$$
  \cL \>=\> \LL \>\>\bigcup\>\> \LD \>\>\bigcup\>\> \LM,
$$
with $\LL$ grouped naturally into pairs (initial multiplicities on the two sides of the loop).
The \textbf{weight} of $l\in\LM$ is $\gcd(l,m)$, the weight of the corresponding inner chain.
The data
\begin{center}
$T_\G=\mgOLfull$
\end{center}
is called a \textbf{principal type}, of \textbf{Euler characteristic} $\chi_{T_\G}=\chi_\G$,
\textbf{weights} $(\gcd(l,m))_{l\in\LM}$ and \textbf{shape}
  ($\chi_\G$; weights).
$$
\pbox[c]{18cm}{\begin{tikzpicture}[TPMain,scale=1.4,edgemults/.style={above right=-0.9 and -0.14,scale=0.8,blue},
    edgemult/.style={scale=0.9},inner/.style={purple!80!black},otherprin/.style={very thick,inner}]
  \draw[-] (0,0)
  edge[very thick] ++(25,0)
  ++(0.5,1) edge[ethinedge] ++(0.3,-1) node[edgemults] {$\!o_1$}
    node[scale=0.9,below right=28pt and 3pt] {\rm outer chains $\cO$}
  ++(0,0) edge[thinedge] ++(0.3,1) ++(0.3,1) edge[thinedge] ++(-0.3,1) ++(0.2,-2)
  ++(0.5,1) edge[ethinedge] ++(0.3,-1) node[edgemults] {$\!o_2$} ++(0,0) edge[thinedge] ++(0.3,1) ++(0.3,1) ++(0.2,-2)
  node[edgemult,above right=0.3 and 0.05] {$...$} ++(1,0)
  ++(0.5,1) edge[ethinedge] ++(0.3,-1) 
  ++(0,0) edge[thinedge] ++(0.3,1) ++(0.3,1) edge[thinedge] ++(-0.3,1) ++(0.2,-2)
  ++(13.5,0)
  ++(0.5,1) node[edgemults] (CX) {$l_1$}
    node[scale=0.9,below right=28pt and 3pt,inner] (LMlabel) {\rm edges $\LM$}
  ++(0,0) edge[inner,ethinedge] ++(0.3,-1) edge[inner,thinedge] ++(0.3,1) ++(0.3,1) edge[inner,thinedge] ++(-0.3,1) ++(-1.5,1) edge[otherprin] ++(1.5,0) ++(1.5,-1) ++(0.2,-2)
  ++(0.5,1) node[edgemults] {$l_2$} ++(0,0) edge[inner,ethinedge] ++(0.3,-1) edge[inner,thinedge] ++(0.3,1) ++(0.3,1) edge[inner,thinedge] ++(-0.3,1) ++(-0.6,1) edge[otherprin] ++(1.6,0) ++(0.6,-1) ++(0.2,-2)
  node[edgemult,above right=0.3 and 0.05] {$...$} ++(1,0)
  ++(0.5,1)
  ++(0,0) edge[inner,ethinedge] ++(0.3,-1) edge[inner,thinedge] ++(0.3,1) ++(-0.1,1) edge[otherprin] ++(1.6,0) ++(0.1,-1) ++(0.3,1) ++(0.2,-2)
  ++(2.5,0)
  node[mainmult,above=-1.5pt,scale=1.1,blue] {$m\ {\rm g}g$}    
  ++(-25.0,0.2) node[scale=1] {$\G$}
  ++(4,-2)
  ++(1.5,-1)
  ;
  \node [below=1pt of LMlabel, inner sep=0pt, scale=0.9, inner] {weights $(\gcd(l_1,m), \gcd(l_2,m),\ldots)$};
  \path[draw,xscale=1.6,yscale=3] (4,0) ++(0,-0.1)--++(115:0.5)++(295:0.1)++(245:0.10)--++(65:0.5)++(245:0.1)++(180:0.1)
            --++(0:0.18)++(0:0.07)--++(0:0.05)++(0:0.07)--++(0:0.18)++(180:0.1)++(115:0.1)--++(295:0.5)
            ++(115:0.1)++(65:0.1)--++(245:0.5)++(-0.13,0.15) node[below=3pt,scale=0.8,blue]{$l'_1,\>\>l'_2$}
                node[scale=0.9,below right=17pt and 3pt] {\rm loops $\LL$};
  \path[draw,xscale=1.6,yscale=3] (5,0) ++(0,-0.1)--++(115:0.5)++(295:0.1)++(245:0.10)--++(65:0.5)++(245:0.1)++(180:0.1)
            --++(0:0.18)++(0:0.07)--++(0:0.05)++(0:0.07)--++(0:0.18)++(180:0.1)++(115:0.1)--++(295:0.5)
            ++(115:0.1)++(65:0.1)--++(245:0.5)++(-0.13,0.15) node[below=3pt,scale=0.8,blue]{$l'_3,l'_4$};
  \draw[-] (9,0) node[edgemult,above right=0.3 and 0.05] {$...$};
  \path[draw] (12,0) node[below=2pt,scale=0.8,blue] {$l_1''$}
                node[scale=0.9,below right=15pt and 3pt] {\rm D-tails $\LD$}
     ++(285:0.3)--++(105:1.6)++(285:0.3)
     ++(255:0.3)--++(75:1.6)++(255:0.3)
     ++(285:0.3)--++(105:1.6)++(285:0.3)
     ++(255:0.3)--++(75:1.6)++(255:0.3)
     ++(180:0.7) edge[very thick] ++(0:1.4) node[right=23pt, scale=0.6] {$2n_1$} ++(0:0.3) ++(-90:0.3)-- node[near end, left, scale=0.6] {$n_1$} ++(90:1.2)
     ++(-90:1.2) ++(0:0.8) --node[near end, right, scale=0.6] {$n_1$} ++(90:1.2);
  \path[draw] (11.8,0)
     node[edgemult,above right=0.3 and 1.7] {$...$} ++(2,0)
     node[below=2pt,scale=0.8,blue] {$l_2''$}
     ++(285:0.3)--++(105:1.2)++(285:0.3)
     ++(255:0.3)--++(75:1.2)++(255:0.3)
     ++(285:0.3)--++(105:1.2)++(285:0.3)
     ++(180:0.7) edge[very thick] ++(0:1.6) node[right=25pt, scale=0.6] {$2n_2$} ++(0:0.3) ++(-90:0.3)-- node[near end, left, scale=0.6] {$n_2$} ++(90:1.2)
     ++(-90:1.2) ++(0:1.0) --node[near end, right, scale=0.6] {$n_2$} ++(90:1.2);
  \node[scale=0.9,otherprin] at (20.7,4.3) {Other principal components $\G_i$};
  \node[scale=0.9,otherprin] at (23.1,3.3) {with $\chi_{\G_i}<0$};
\end{tikzpicture}}
$$
\begin{center}
\figuretitle Principal type $T_\G=\mgOLfull$
\end{center}
\end{definition}

Note immediate corollaries of the total genus formula \ref{eulerformula} and of \ref{cormgol} (3):

\begin{corollary}
\label{totalgenusprintype}
Let $R$ be a reduction type of total genus $g\ge 2$. Then
$$
  2-2g = \sum_T \chi_T,
$$
the sum taken over the principal types of all components $\G$ in $R$ with $\chi_\G<0$.
\end{corollary}

\begin{corollary}
There are finitely many possible principal types $T=\mgOLfull$ with fixed $\chi_T$.
\end{corollary}

\begin{notation}
\label{notprintype}
We draw a principal type $T=(m,g,\cO,\LL,\LD,\LM)$ as:
\begin{itemize}
\item
An arc marked `$m$' if $g=0$ and `$m\>{\rm g}g$' if $g>0$, for the principal component.
\item
Outer multiplicity $o\in\cO$ as a segment pointing down, marked $o$.
\item
Pair $l,l'\in\LL$ as a loop pointing up, marked $l-l'$.
\item
Element $l\in\LD$ as a segment pointing up, marked $l{\rm D}$.
\item
Element $l\in\LM$ as a longer segment pointing up, marked $l$.
\end{itemize}
We draw the shape of $T$ as a shaded disk marked with $\chi$, with a segment for each edge marked with
its weight.
\end{notation}

\begin{example}
The following reduction type has principal types
$T_{\G_1}\!=\!(6,1,\{5\},\{6,6\},\{4\},\{3\})$ and $T_{\G_2}\!=\!(3,2,\{\},\{\},\{\},\{3\})$,
of Euler characteristic
$\chi_{\G_1}\!=\!-29$, $\chi_{\G_2}\!=\!-9$ and total genus 16:
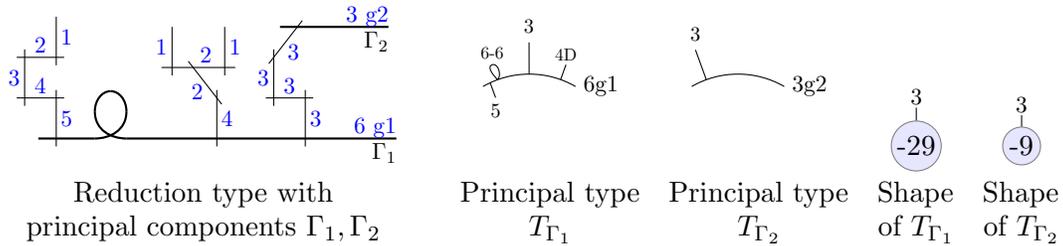
\begin{center}
\begin{tabular}{c@{\qquad}c@{\quad}c@{\quad}c@{\quad}c}
\begin{tikzpicture}[xscale=0.7,yscale=0.9] \draw[thick,shorten >=2] (0.26,0)--(7.16,0) node[inner sep=2pt,blue,above left,scale=0.8] {6 \smash{g}1} node[inner sep=2pt,below left,scale=0.8] {$\Gamma_1$}; \draw[shorten <=-3pt, shorten >=-3pt] (3/5,0)--node[inner sep=2pt,scale=0.8,blue,right] {5} (3/5,3/5); \draw[shorten <=-3pt, shorten >=-3pt] (3/5,3/5)--node[inner sep=2pt,scale=0.8,blue,above] {4} (0,3/5); \draw[shorten <=-3pt, shorten >=-3pt] (0,3/5)--node[inner sep=2pt,scale=0.8,blue,left] {3} (0,6/5); \draw[shorten <=-3pt, shorten >=-3pt] (0,6/5)--node[inner sep=2pt,scale=0.8,blue,above] {2} (3/5,6/5); \draw[shorten <=-3pt] (3/5,6/5)--node[inner sep=2pt,scale=0.8,blue,right] {1} (3/5,9/5); \path[draw,thick] (4/3,0) edge[white,line width=2] ++(0.6,0) edge[out=0,in=-90] ++(0.6,0.4) ++(0.6,0.4) edge[out=90,in=0] ++(-0.3,0.3) ++(-0.3,0.3) edge[out=180,in=90] ++(-0.3,-0.3) ++(-0.3,-0.3) edge[out=-90,in=180] ++(0.6,-0.4); \draw[shorten <=-3pt, shorten >=-3pt] (3.65,0)--node[inner sep=2pt,scale=0.8,blue,right] {4} (3.65,3/5); \draw[shorten <=-3pt, shorten >=-3pt] (3.65,3/5)--node[inner sep=2pt,scale=0.8,blue,below left=-1pt and -1pt] {2} (3.2,21/20); \draw (2.6,21/20)--node[inner sep=2pt,scale=0.8,blue,above,pos=0.6] {2} (4,21/20); \draw[shorten <=-3pt] (2.8,21/20)--node[inner sep=2pt,scale=0.8,blue,left] {1} (2.8,33/20); \draw[shorten <=-3pt] (3.8,21/20)--node[inner sep=2pt,scale=0.8,blue,right] {1} (3.8,33/20); \draw[shorten <=-3pt, shorten >=-3pt] (5.33,0)--node[inner sep=2pt,scale=0.8,blue,right] {3} (5.33,3/5); \draw[shorten <=-3pt, shorten >=-3pt] (5.33,3/5)--node[inner sep=2pt,scale=0.8,blue,above] {3} (4.73,3/5); \draw[shorten <=-3pt, shorten >=-3pt] (4.73,3/5)--node[inner sep=2pt,scale=0.8,blue,left] {3} (4.73,6/5); \draw[shorten <=-3pt, shorten >=-3pt] (4.73,6/5)--node[inner sep=2pt,scale=0.8,blue,below right=-1pt and -1pt] {3} (5.18,33/20); \draw[thick,shorten >=2] (4.85,33/20)--(7.01,33/20) node[inner sep=2pt,blue,above left,scale=0.8] {3 \smash{g}2} node[inner sep=2pt,below left,scale=0.8] {$\Gamma_2$}; \end{tikzpicture}
&
\begin{tikzpicture}[scale=1]
\draw(0,0)
\shapearc{35pt}{6-6,3,.4}{5}{6g1}
;
\end{tikzpicture}
&
\begin{tikzpicture}[scale=1]
\draw(0,0)
\shapearc{35pt}{3}{}{3g2}
;
\end{tikzpicture}
&
\begin{tikzpicture}[baseline=-8] \node[pPShapeChi] (1) at (0,0) {-29}; \draw[shorten >= 4pt] (1) -- ++(90:0.65) node[pPShapeWeight] {3}; \end{tikzpicture}
&
\begin{tikzpicture}[baseline=-8] \node[pPShapeChi] (2) at (0,0) {-9}; \draw[shorten >= 4pt] (2) -- ++(90:0.55) node[pPShapeWeight] {3}; \end{tikzpicture}
\cr
Reduction type with                & Principal type  & Principal type  &   Shape          &  Shape \cr
principal components $\G_1, \G_2$  &   $T_{\G_1}$    &   $T_{\G_2}$    &   of $T_{\G_1}$  &  of $T_{\G_2}$ \cr
\end{tabular}
\end{center}
\end{example}

Principal components of the same shape behave like interchangeable pieces: one can be swapped by another
in a reduction type without changing the total genus (see \S\ref{sintro}).
In our example, there are many
reduction types of genus 16 with two principal types of shape $(-29; 3)$ and $(-9; 3)$
and an edge between them of weight 3.
A configuration like this is called the shape of a reduction type, and is drawn as follows:

\begin{center}
\scalebox{1.1}{\begin{tikzpicture}[xscale=0.8,yscale=0.8,auto=left, v/.style={pShapeVertex}, l/.style={pShapeEdge}] \node[v] (1) at (1,1) {-29}; \node[v] (2) at (2,1) {-9}; \draw[] (1)--node[l] {3} (2); \end{tikzpicture}}
\end{center}

\begin{definition}[Shape]
\label{defshape}
The \textbf{shape} of a reduction type of genus $g>1$ is a graph whose
\begin{itemize}
\item
vertices are principal components $\G$ with $\chi_\G<0$, marked with $\chi_\G$, and
\item
edges are inner chains between them, marked with their weights (`1' is usually omitted).
\end{itemize}
\end{definition}

\begin{example}[Genus 2 classification]
\label{exg2}
Let $R$ be a reduction family of total genus 2. Then $\chi(R)=-2=\sum_T \chi_T$ as in \ref{totalgenusprintype},
and each principal type $T$ in $R$ has $\chi_T<0$. So $R$ consists of
\begin{itemize}
\item
one principal type $T$ with $\chi_T=-2$, with no edges; or
\item
two principal types $T_1, T_2$ with $\chi_{T_1}=\chi_{T_2}=-1$ and identical weights, linked together.
\end{itemize}

As in the proof of \ref{cormgol}, we can find all principal types $T$ with $\chi_T=-1$ with edges, and with $\chi_T=-2$
and no edges. There are 13 and 46 respectively (p.~\pageref{tabPS1}), and the 13 are the 10
`Kodaira types with an edge of weight 1' (what Ogg \cite{Og} denotes Kod($\G$)) plus three singletons:

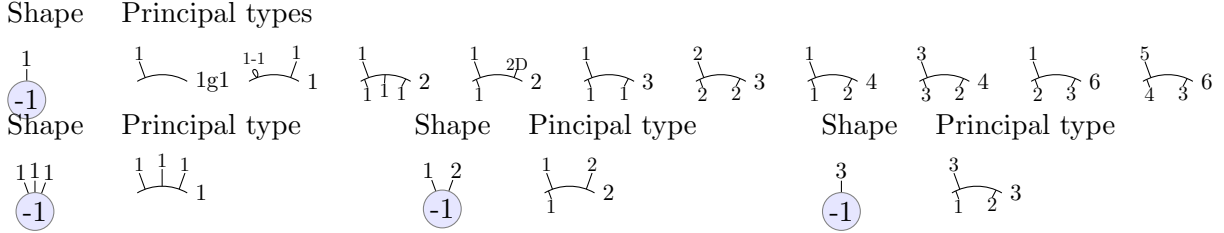
\begin{center}
\def\shapearclabel#1#2#3#4#5{\shapearc{#1}{#2}{#3}{#4}}
\begin{tabular}{p{0.75cm}@{\qquad}p{14.4cm}} Shape & Principal types\vphantom{$\int^X$}\\[-6pt] \vtop{\null \begin{tikzpicture}[baseline=0] \node[pPShapeChi] (1) at (0,0) {-1}; \draw[shorten >= 4pt] (1) -- ++(90:0.55) node[pPShapeWeight] {1}; \end{tikzpicture}} & \vtop{\null \begin{tikzpicture}[scale=0.54] \path[use as bounding box] (-1,-0.09) rectangle (26.5,2); \draw(0,0) \shapearclabel{35pt}{1}{}{1g1}{\redtype{I_{g1}\e{1}}} ++(2.75,0) \shapearclabel{35pt}{1-1,1}{}{1}{\redtype{I_1\e{1}}} ++(2.75,0) \shapearclabel{35pt}{1}{1,1,1}{2}{\redtype{I^*_0\e{1}}} ++(2.75,0) \shapearclabel{35pt}{1,.2}{1}{2}{\redtype{I^*_1\e{1}}} ++(2.75,0) \shapearclabel{35pt}{1}{1,1}{3}{\redtype{IV\e{1}}} ++(2.75,0) \shapearclabel{35pt}{2}{2,2}{3}{\redtype{IV^*\e{2}}} ++(2.75,0) \shapearclabel{35pt}{1}{1,2}{4}{\redtype{III\e{1}}} ++(2.75,0) \shapearclabel{35pt}{3}{3,2}{4}{\redtype{III^*\e{3}}} ++(2.75,0) \shapearclabel{35pt}{1}{2,3}{6}{\redtype{II\e{1}}} ++(2.75,0) \shapearclabel{35pt}{5}{4,3}{6}{\redtype{II^*\e{5}}}; \end{tikzpicture} }\\[-10pt]\end{tabular}\par\vskip-1pt\noindent\begin{tabular}{p{0.75cm}@{\qquad}p{3.544cm}p{0.75cm}@{\qquad}p{3.544cm}p{0.75cm}@{\qquad}p{3.544cm}}Shape & Principal type & Shape & Pincipal type & Shape & Principal type\vphantom{$\int^X$}\\[-6pt] \vtop{\null \begin{tikzpicture}[baseline=0] \node[pPShapeChi] (1) at (0,0) {-1}; \draw[shorten >= 4pt] (1) -- ++(110:0.55) node[pPShapeWeight] {1}; \draw[shorten >= 4pt] (1) -- ++(90:0.55) node[pPShapeWeight] {1}; \draw[shorten >= 4pt] (1) -- ++(70:0.55) node[pPShapeWeight] {1}; \end{tikzpicture}} & \vtop{\null \begin{tikzpicture}[scale=0.54] \path[use as bounding box] (-1,-0.09) rectangle (1.75,2); \draw(0,0) \shapearclabel{35pt}{1,1,1}{}{1}{\redtype{I\e{1}\e{1}\e{1}}}; \end{tikzpicture} }& \vtop{\null \begin{tikzpicture}[baseline=0] \node[pPShapeChi] (1) at (0,0) {-1}; \draw[shorten >= 4pt] (1) -- ++(110:0.55) node[pPShapeWeight] {1}; \draw[shorten >= 4pt] (1) -- ++(70:0.55) node[pPShapeWeight] {2}; \end{tikzpicture}} & \vtop{\null \begin{tikzpicture}[scale=0.54] \path[use as bounding box] (-1,-0.09) rectangle (1.75,2); \draw(0,0) \shapearclabel{35pt}{1,2}{1}{2}{\redtype{D\e{1}\e{2}}}; \end{tikzpicture} }& \vtop{\null \begin{tikzpicture}[baseline=0] \node[pPShapeChi] (1) at (0,0) {-1}; \draw[shorten >= 4pt] (1) -- ++(90:0.55) node[pPShapeWeight] {3}; \end{tikzpicture}} & \vtop{\null \begin{tikzpicture}[scale=0.54] \path[use as bounding box] (-1,-0.09) rectangle (1.75,2); \draw(0,0) \shapearclabel{35pt}{3}{1,2}{3}{\redtype{T\e{3}}}; \end{tikzpicture} }\\ \end{tabular}
\end{center}

\noindent
Therefore, genus 2 reduction types come in 5 shapes
\smallskip

\begin{center}
\cbox{\begin{tikzpicture}[xscale=0.7,yscale=0.7,auto=left, v/.style={pShapeVertex}] \node[v] (1) at (1,1) {-2}; \end{tikzpicture}}\qquad\cbox{\begin{tikzpicture}[xscale=0.8,yscale=0.8,auto=left, v/.style={pShapeVertex}, l/.style={pShapeEdge}] \node[v] (1) at (1,1) {-1}; \node[v] (2) at (2,1) {-1}; \draw[] (1)--node[l] {} (2); \end{tikzpicture}}\qquad\cbox{\begin{tikzpicture}[xscale=1,yscale=1,auto=left, v/.style={pShapeVertex}, l/.style={pShapeEdge}, e/.style={pe3}] \node[v] (1) at (1,1) {-1}; \node[v] (2) at (2,1) {-1}; \draw[e] (1)--node[l] {1,1,1} (2); \end{tikzpicture}}\qquad\cbox{\begin{tikzpicture}[xscale=0.9,yscale=0.9,auto=left, v/.style={pShapeVertex}, l/.style={pShapeEdge}, e/.style={pe2}] \node[v] (1) at (1,1) {-1}; \node[v] (2) at (2,1) {-1}; \draw[e] (1)--node[l] {1,2} (2); \end{tikzpicture}}\qquad\cbox{\begin{tikzpicture}[xscale=0.8,yscale=0.8,auto=left, v/.style={pShapeVertex}, l/.style={pShapeEdge}] \node[v] (1) at (1,1) {-1}; \node[v] (2) at (2,1) {-1}; \draw[] (1)--node[l] {3} (2); \end{tikzpicture}}
\end{center}

\smallskip
\noindent
with 46, $\binom{10+1}2=55$, 1, 1, 1 families of each shape, respectively (104 in total).

\begin{center}
\begin{tabular}{c@{\qquad}c@{\qquad}c@{\qquad}c@{\qquad}c}
\cbox{\begin{tikzpicture}[xscale=0.55,yscale=0.7] \draw[thick,shorten >=2] (-1/3,0)--(11/6,0) node[inner sep=2pt,blue,above left,scale=0.8] {1 \smash{g}2}; \end{tikzpicture} }&\cbox{\begin{tikzpicture}[xscale=0.55,yscale=0.7] \draw[thick,shorten >=2] (-1/3,0)--(11/6,0) node[inner sep=2pt,blue,above left,scale=0.8] {1 \smash{g}1}; \draw[blue!80!white,densely dashed] (0,0)--(0,4/5);\draw[thick,shorten >=2] (-1/3,4/5)--(11/6,4/5) node[inner sep=2pt,blue,above left,scale=0.8] {1 \smash{g}1}; \end{tikzpicture} }&\cbox{\begin{tikzpicture}[xscale=0.55,yscale=0.7] \draw[thick,shorten >=2] (-1/3,0)--(89/30,0) node[inner sep=2pt,blue,above left,scale=0.8] {1}; \draw[blue!80!white,densely dashed] (0,0)--(0,4/5);\draw[blue!80!white,densely dashed] (2/3,0)--(2/3,4/5);\draw[blue!80!white,densely dashed] (4/3,0)--(4/3,4/5);\draw[thick,shorten >=2] (-1/3,4/5)--(89/30,4/5) node[inner sep=2pt,blue,above left,scale=0.8] {1}; \end{tikzpicture} }&\cbox{\begin{tikzpicture}[xscale=0.55,yscale=0.7] \draw[thick,shorten >=2] (-1/3,0)--(7/2,0) node[inner sep=2pt,blue,above left,scale=0.8] {2}; \draw[shorten <=-3pt] (0,0)--node[inner sep=2pt,scale=0.8,blue,right] {1} (0,3/5); \draw[shorten <=-3pt, shorten >=-3pt] (17/15,0)--node[inner sep=2pt,scale=0.8,blue,right] {1} (17/15,4/5); \draw[blue!80!white,densely dashed] (28/15,0)--(28/15,4/5);\draw[thick,shorten >=2] (4/5,4/5)--(41/10,4/5) node[inner sep=2pt,blue,above left,scale=0.8] {2}; \draw[shorten <=-3pt] (37/15,4/5)--node[inner sep=2pt,scale=0.8,blue,right] {1} (37/15,7/5); \end{tikzpicture} }&\cbox{\begin{tikzpicture}[xscale=0.55,yscale=0.7] \draw[thick,shorten >=2] (-1/3,0)--(25/6,0) node[inner sep=2pt,blue,above left,scale=0.8] {3}; \draw[shorten <=-3pt] (0,0)--node[inner sep=2pt,scale=0.8,blue,right] {1} (0,3/5); \draw[shorten <=-3pt, shorten >=-3pt] (7/5,0)--node[inner sep=2pt,scale=0.8,blue,right] {2} (7/5,3/5); \draw[shorten <=-3pt] (7/5,3/5)--node[inner sep=2pt,scale=0.8,blue,above] {1} (4/5,3/5); \draw[blue!80!white,densely dashed] (38/15,0)--(38/15,4/5);\draw[thick,shorten >=2] (11/5,4/5)--(37/6,4/5) node[inner sep=2pt,blue,above left,scale=0.8] {3}; \draw[shorten <=-3pt] (47/15,4/5)--node[inner sep=2pt,scale=0.8,blue,right] {1} (47/15,7/5); \draw[shorten <=-3pt, shorten >=-3pt] (68/15,4/5)--node[inner sep=2pt,scale=0.8,blue,right] {2} (68/15,7/5); \draw[shorten <=-3pt] (68/15,7/5)--node[inner sep=2pt,scale=0.8,blue,above] {1} (59/15,7/5); \end{tikzpicture} }
\Bcr
one of & one of & unique & unique & unique\Tcr
46 families & 55 families & family & family & family\cr
\end{tabular}\\[4pt]
\tabletitle\label{tabSE2} Simplest reduction family of each shape in genus 2.
\end{center}

\noindent
Table \ref{tabNU} lists all 104 families with their labels (see \S\ref{slabel})
and the Namikawa--Ueno names \cite{NU}.
\end{example}

The same approach works in every genus:

\begin{algorithm}
\label{algmain}
(Computing all reduction families in a given genus $g\ge 2$)
\noindent
\begin{enumerate}
\item Compute all cores $\Psi=m^{o_1,...,o_k}$ with $2-2g\le \chi(\Psi)\le 0$, using the bounds $2\le m\le 6-2\chi(\Psi)$ on $m$,
and $3\le k\le 4-\chi(\Psi)$ on $k$ (Theorem \ref{coremain}).
\item For each $1\le B\le 2g\!-\!2$ and each core $\Psi$ from (1) and cores with $\chi(\Psi)=2$,
determine possible tuples of invariants $\mgOL$
of principal components $\G$ with $\chi_\G=-B$ and core $\Psi$, following the proof of \ref{cormgol}.
\item Compute all principal types from the tuples in (2), by splitting $\cL$ into $\LD$, $\LL$, $\LM$
in possible ways (Definition \ref{defcomfam}). When $\chi_\G=2\!-\!2g$ (maximal possible),
we require $\LM=\emptyset$.
\item Group principal types by their shape $(\chi; \text{weights})$.
\item For each number of vertices $N\in[1,...,2g-2]$:\\
For each integer partition $p_1,...,p_n$ of $2g-2$ into $N$ positive parts:\\
For each connected simple graph $G$ with $N$ vertices $v_1,...,v_N$:
\item Construct all possible shapes $S$ with underlying graph $G$ and principal shapes from (3) with $\chi=-p_i$
  at a vertex $v_i$, up to isomorphism.
\item Loop through possible allocations of a principal type at every vertex $v$ with assigned shape, matching its
edges to the endpoints of $v$ in possible ways. Each such assignment gives a reduction family.
\item End of loop (5).
\end{enumerate}
\end{algorithm}

\begin{example}[Genus 3 classification]
\label{exg3}
This is how the algorithm works in genus 3:

\smallskip\noindent\textbf{Step 1.}
Using the bounds $m\le 6-2\cdot(-4)=14$ and the bound on $k$ we find

\begin{tabular}{llllll}
7
&cores\ \
${\rm I}_0^*$, ${\rm IV}$, ${\rm IV}^*$, ${\rm III}$, ${\rm III}^*$, ${\rm II}$, ${\rm II}^*$
& with $\chi=0$ & (Table \ref{tabC0}),\cr
16
&cores\ \
$2^{1,1,1,1,1,1}$, $3^{1,1,2,2}$, $4^{1,3,2,2}$, $5^{1,1,3}$, \ldots
& with $\chi=-2$ & (Table \ref{tabC-2}), \cr
43
&cores\ \
$2^{1,1,1,1,1,1,1,1}$, $3^{1,1,1,1,2}$, $3^{1,2,2,2,2}$, \ldots
& with $\chi=-4$ & (Table \ref{tabC-4}).
\end{tabular}

\smallskip\noindent\textbf{Steps 2--4.}
Consider each core $\Psi$ from Step 1 in turn, as well as
cores $\Psi=1^\emptyset$ and $\Psi=m^{a,m-a}$ that have $\chi(\Psi)=2$.
Write $\Psi=(m_C)^{o_1,...,o_k}$, and list tuples $\mgOL$ with core $\Psi$ and $\chi\ge -4$,
using inequalities in the proof of \ref{cormgol}:
\begin{itemize}
\item[(a)]
Select a multiple $[d]\Psi$ for some $d\ge 1$. This determines $m=dm_C$ and the non-zero elements
$do_1,...,do_k$ of $\cO\coprod\cL$.
\item[(b)]
Select genus $g\ge 0$.
\item[(c)]
Select which of the $do_i$ are in $\cO$ and which ones are in $\cL$.
\item[(d)]
Add some number of $m$'s (i.e. zero elements in $\Z/m\Z$) to $\cL$.
\item[(e)]
Split $\cL$ further into $\LL$, $\LD$ and $\LM$.
\item[(f)]
Keep all the obtained principal types $T$ with $\chi_T>-4$ and $\LM\ne\emptyset$
(to be used in reduction families with multiple principal types) and those
of weight $(-4, \emptyset)$.
\end{itemize}

\noindent
The principal types $T$ with $-4\le \chi_T\le -1$ are listed in Tables \ref{tabPS1}--\ref{tabPS4}.

\smallskip\noindent\textbf{Steps 5--6.}
Take simple connected graphs with $\le 4$ vertices, allocate $\chi\le -1$ to each vertex so that they add up
to $-4$, and go through the possibilities of assigning a multiset of weights to every edge, in a way
that every vertex acquires a shape ($\chi$, weights) for which we there are principal types. This gives
the
35
possible shapes in genus 3, listed in Table \ref{tabR3} on p.\pageref{tabR3}:
\begin{center}
\begin{tikzpicture}[xscale=0.7,yscale=0.7,auto=left, v/.style={pShapeVertex}] \node[v] (1) at (1,1) {-4}; \end{tikzpicture}\quad\begin{tikzpicture}[xscale=0.8,yscale=0.8,auto=left, v/.style={pShapeVertex}, l/.style={pShapeEdge}] \node[v] (1) at (1,1) {-3}; \node[v] (2) at (2,1) {-1}; \draw[] (1)--node[l] {3} (2); \end{tikzpicture}\quad\begin{tikzpicture}[xscale=0.8,yscale=0.8,auto=left, v/.style={pShapeVertex}, l/.style={pShapeEdge}] \node[v] (1) at (1,1) {-3}; \node[v] (2) at (2,1) {-1}; \draw[] (1)--node[l] {} (2); \end{tikzpicture}\quad\begin{tikzpicture}[xscale=1,yscale=1,auto=left, v/.style={pShapeVertex}, l/.style={pShapeEdge}, e/.style={pe3}] \node[v] (1) at (1,1) {-3}; \node[v] (2) at (2,1) {-1}; \draw[e] (1)--node[l] {1,1,1} (2); \end{tikzpicture}\quad\begin{tikzpicture}[xscale=0.9,yscale=0.9,auto=left, v/.style={pShapeVertex}, l/.style={pShapeEdge}, e/.style={pe2}] \node[v] (1) at (1,1) {-3}; \node[v] (2) at (2,1) {-1}; \draw[e] (1)--node[l] {1,2} (2); \end{tikzpicture}\quad\begin{tikzpicture}[xscale=0.9,yscale=0.9,auto=left, v/.style={pShapeVertex}, l/.style={pShapeEdge}, e/.style={pe2}] \node[v] (1) at (1,1) {-2}; \node[v] (2) at (2,1) {-2}; \draw[e] (1)--node[l] {1,3} (2); \end{tikzpicture}\quad\begin{tikzpicture}[xscale=0.8,yscale=0.8,auto=left, v/.style={pShapeVertex}, l/.style={pShapeEdge}] \node[v] (1) at (1,1) {-2}; \node[v] (2) at (2,1) {-2}; \draw[] (1)--node[l] {2} (2); \end{tikzpicture}\quad\begin{tikzpicture}[xscale=0.9,yscale=0.9,auto=left, v/.style={pShapeVertex}, l/.style={pShapeEdge}, e/.style={pe2}] \node[v] (1) at (1,1) {-2}; \node[v] (2) at (2,1) {-2}; \draw[e] (1)--node[l] {1,1} (2); \end{tikzpicture}
\quad$\cdots$
\end{center}

\smallskip\noindent\textbf{Steps 7--8.}
Pick a principal type for every vertex of the prescribed shape, and assign its edges to the edges of
the vertex matching the weight. This gives one reduction family, and we find 1901 of them in total in genus 3
(again see Table \ref{tabR3}).

\end{example}

Running the algorithm, we find:

\begin{theorem}
\label{numfamilies}
The number of reduction families in genus $2\le g\le 6$ is:
\begin{itemize}
\item
\gtwonumtypes{} families in genus 2, in 5 shapes (Table \ref{tabR2} on p. \pageref{tabR2}, see Example \ref{exg2}),
\item
\gthreenumtypes{} families in genus 3, in 35 shapes (Table \ref{tabR3} on p. \pageref{tabR3}, see Example \ref{exg3}),
\item
\gfournumtypes{} families in genus 4, in 310 shapes (Table \ref{tabR4} on p. \pageref{tabR4}),
\item
\gfivenumtypes{} families in genus 5, in 3700 shapes (some in Table \ref{tabR5} on p. \pageref{tabR5}),
\item
\gsixnumtypes{} families in genus 6, in 56253 shapes (some in Table \ref{tabR6} on p. \pageref{tabR6}).
\end{itemize}
\end{theorem}

The Reduction Type library {\tt redlib} \cite{redlib} allows to do this computation in Magma \cite{Magma}, Python or Javascript,
provided that all shapes in that genus are known. The shapes were generated in Magma using the Small Graph Database
for connected graphs on up to 10 vertices (which is as far as the database goes, so $2g-2=10\iff g=6$ is the limit).

\section{Canonical graph traversing}

To give a reduction type a canonical label, we introduce a labelling convention for any graph
(finite connected undirected multigraph) with marked vertices and edges.
We assume that the possible marks have a total order of them, so that they can be compared and sorted.

\begin{definition}
\label{graphpath}
Let $G$ be a graph on $n$ vertices.
A \textbf{path} is a sequence of symbols of the form
\begin{center}
  vertex, edge, vertex, edge, \ldots, edge, vertex
\end{center}
where
\begin{enumerate}
\item
`vertex' is a symbol $i\in\{1, 2, ..., n\}$ (for reduction types we use `$c_i$' instead)
if it appeared before on the path and as an $i$th vertex
for the first time, and the mark of a vertex otherwise. (The former are considered larger than the vertex marks.)
\item
`edge' is a mark of an edge between the two vertices referred to by the two neighbouring symbols, or a symbol
`\&' for a jump between two vertices without following an edge.
(The latter is considered larger than the edge marks.)
\item
Every edge in the graph appears exactly once in the path.
\end{enumerate}
\end{definition}

\begin{definition}
\label{minpath}
The order on the marks induces a lexicographic ordering of the paths. Call a path \textbf{minimal}
if
\begin{enumerate}
\item
There are no paths of smaller length.
\item
It is lexicographically smallest among all those that satisfy (1).
\end{enumerate}

\end{definition}

\begin{example}
Here is a graph with vertices marked A, C, E and edges marked b, d.
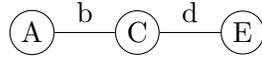
$$
\begin{tikzpicture}[scale=0.7,n/.style={circle,draw,inner sep=2pt}]
  \path(1,1) node[n] (1) {A};
  \path(3,1) node[n] (2) {C};
  \path(5,1) node[n] (3) {E};
  \draw[-,color=black]
    (1) -- node[above]{b} (2) -- node[above]{d} (3);
\end{tikzpicture}
$$
It has two paths that satisfy (1),
\begin{center}
A, b, C, d, E \qquad\text{and}\qquad E, d, C, b, A.
\end{center}
The first one is minimal for the usual lexicographic ordering of the marks.
Similarly,
$$
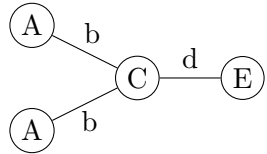
\begin{tikzpicture}[scale=0.7,n/.style={circle,draw,inner sep=2pt}]
  \path(1,2) node[n] (1) {A};
  \path(3,1) node[n] (2) {C};
  \path(5,1) node[n] (3) {E};
  \path(1,0) node[n] (4) {A};
  \draw[-,color=black]
    (1) -- node[pos=0.4,above right, inner sep=2pt]{b} (2) -- node[above]{d} (3);
  \draw[-,color=black]
    (2) -- node[pos=0.6,below right, inner sep=1pt]{b} (4);
\end{tikzpicture}
$$
has the minimal path (that can be traversed in two different ways in the graph)
\begin{center}
A, b, C, b, A, \&, E, d, 2.
\end{center}
The non-shortest paths have more jumps, for example adding `\&, A' at the end.

\end{example}

\begin{lemma}
Suppose $G$ has $2o$ vertices of odd degree, and $m>0$ edges.
Then the minimal paths in $G$ have length $2m+1$ if $o=0$ and $2m+o$ otherwise.
\end{lemma}

\begin{proof}
As in the examples above, a minimal path minimises the number of jumps between vertices.
A path with no jumps
has length $2m+1$, and it exists if only if $G$ has an Eulerian path, equivalently $o\le 1$.
When $o>1$, the minimal number of jumps is $o-1$ (this is known as the Chinese Postman Problem),
and the total length is $2m+o$.
\end{proof}

\begin{remark}
There is an unmarked analogue as well, replacing the vertex marks by `v' and/or edge marks by `e'.
Also, one can extend the lemma to unconnected graphs.
\end{remark}

\begin{remark}
A path describes how to construct a graph, by adding vertices (for every `v') and
edges (for every `e'). In particular,
\begin{enumerate}
\item
Any path determines the graph up to isomorphism.
\item
The minimal path $P(G)$ is unique, and can be used as a fingerprint
for isomorphism testing of marked multigraphs:
$$
  G_1 \iso G_2 \qquad\liff\qquad P(G_1)=P(G_2).
$$
\item
The ways to traverse a minimal path through vertices and edges of $G$ form an $\Aut G$-torsor,
so there are many of them if $\Aut G$ is large.
\end{enumerate}
\end{remark}

\begin{remark}
An old algorithm of Fleury \cite{Fleury} can be adapted to construct
a minimal path in a a marked multigraph. For reduction types, it is polynomial in the genus,
and quite efficient in practice. The idea is to use odd degree vertices as starting
points (or all vertices if all degrees are even), and extend possible paths inductively,
keeping only the ones that could eventually extend to a shortest path, discarding all but
the lexicographically smallest ones along the way.
\end{remark}

\section{Canonical label}
\label{slabel}

\def\loops{\langle\text{loops}\rangle}
\def\dtails{\langle\text{D-tails}\rangle}
\def\edges{\langle\text{edges}\rangle}
\def\outedges{\langle\text{outgoing edges}\rangle}
\def\innertypes{\langle\text{inner types}\rangle}
\def\sprintype{\langle\text{standalone principal type}\rangle}
\def\printype{\langle\text{principal type}\rangle}
\def\rtype{\langle\text{reduction type}\rangle}

Let $R$ be a reduction type. We can view it as a simple graph $G$ whose vertices
are marked with principal types of components $\G$ with $\chi_\G<0$ and edges with sets of inner chains between them.
We will attach labels to those, so that every path through $G$ (\ref{graphpath}) determines a valid
label for $R$. Then we order everything, and use minimal path (\ref{minpath}) for the canonical label.

\begin{notation}[Label]
\label{redlabel}
If $\G_1-\G_2-\ldots$ is a path through $G$, the corresponding label of $R$ is

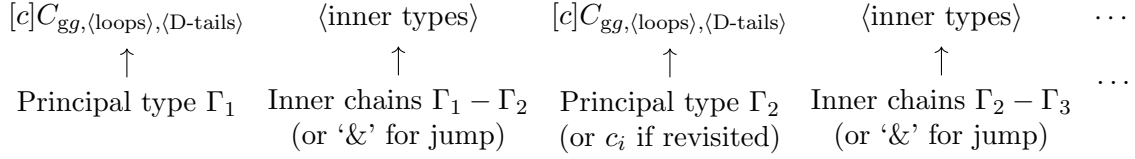
\begin{tikzpicture}[node distance=1em and 1.5em, every node/.style={anchor=center, inner sep=5pt}]
  \node (n1) at (0,0) {$[c] C_{{\rm g} g, \loops, \dtails}$};
  \node (n2) [right=of n1] {$\innertypes$};
  \node (n3) [right=of n2] {$[c] C_{{\rm g} g, \loops, \dtails}$};
  \node (n4) [right=of n3] {$\innertypes$};
  \node (n5) [right=of n4] {$\cdots$};
  \node (l1) [below=of n1, align=center] {Principal type $\G_1$};
  \node (l2) [below=of n2, align=center] {Inner chains $\G_1-\G_2$\\(or `$\&$' for jump)};
  \node (l3) [below=of n3, align=center] {Principal type $\G_2$\\(or $c_i$ if revisited)};
  \node (l4) [below=of n4, align=center] {Inner chains $\G_2-\G_3$\\(or `$\&$' for jump)};
  \node (l5) [below=of n5, align=center] {$\ldots$};
  \draw[->] (l1) -- (n1);
  \draw[->] (l2) -- (n2);
  \draw[->] (l3) -- (n3);
  \draw[->] (l4) -- (n4);
\end{tikzpicture}

\end{notation}

\subsection{Ordering multiplicities}\label{ssOrder}

We represent elements of $\Z/m\Z$, as $1,2,...,m\subset\N$,
and sort them first by gcd with $m$, and then by the element itself, in increasing order. For example,
$$
  m=4, \quad \cL=\{1,1,3,3,4\}, \quad \cO=\{1,3,2,2\}.
$$
We call this `sorting by $\gcd(l,m), l$' and sort loops and D-tails similarly below.

\subsection{Labels for principal types}

Let $\mgOLfull$ be the principal type of some component $\G_i$ with $\chi<0$.
(In genus 1, the labelling scheme needs one tweak.
There are no components with $\chi_\G<0$, so we pick a principal one $\G_1$
(or, in the exceptional case $[c]\In{n}$, any $\G_1$), and label its principal type.
This gives Kodaira types or their multiples, as expected.)

\subsection{Weight, core and genus}

In \ref{redlabel},

\begin{itemize}[leftmargin=\parindent]
\item
$[c]$ shows the weight
$c\!=\!\gcd(m,\cO,\LL,\LD,\LM)$ of $\G$ (cf. \ref{princore}), omitted if $c=1$.
\item
$C$ is the \textbf{core} of $\G$ as in \ref{princore}.
\item
${\rm g}g$ shows the geometric genus of $\G$, omitted if $g=0$.
\end{itemize}

\noindent
Here are some reduction types with no loops, D-tails or inner chains, and their labels:
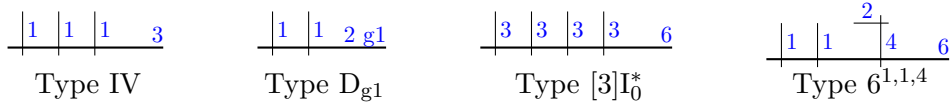
\begin{center}
\begin{tabular}{c@{\qquad\quad}c@{\qquad\quad}c@{\qquad\quad}c}
\cbox{\begin{tikzpicture}[xscale=0.6,yscale=0.8] \draw[thick,shorten >=2] (-1/3,0)--(97/30,0) node[inner sep=2pt,blue,above left,scale=0.8] {3}; \draw[shorten <=-3pt] (0,0)--node[inner sep=2pt,scale=0.8,blue,right] {1} (0,3/5); \draw[shorten <=-3pt] (4/5,0)--node[inner sep=2pt,scale=0.8,blue,right] {1} (4/5,3/5); \draw[shorten <=-3pt] (8/5,0)--node[inner sep=2pt,scale=0.8,blue,right] {1} (8/5,3/5); \end{tikzpicture} }&\cbox{\begin{tikzpicture}[xscale=0.6,yscale=0.8] \draw[thick,shorten >=2] (-1/3,0)--(79/30,0) node[inner sep=2pt,blue,above left,scale=0.8] {2 \smash{g}1}; \draw[shorten <=-3pt] (0,0)--node[inner sep=2pt,scale=0.8,blue,right] {1} (0,3/5); \draw[shorten <=-3pt] (4/5,0)--node[inner sep=2pt,scale=0.8,blue,right] {1} (4/5,3/5); \end{tikzpicture} }&\cbox{\begin{tikzpicture}[xscale=0.6,yscale=0.8] \draw[thick,shorten >=2] (-1/3,0)--(121/30,0) node[inner sep=2pt,blue,above left,scale=0.8] {6}; \draw[shorten <=-3pt] (0,0)--node[inner sep=2pt,scale=0.8,blue,right] {3} (0,3/5); \draw[shorten <=-3pt] (4/5,0)--node[inner sep=2pt,scale=0.8,blue,right] {3} (4/5,3/5); \draw[shorten <=-3pt] (8/5,0)--node[inner sep=2pt,scale=0.8,blue,right] {3} (8/5,3/5); \draw[shorten <=-3pt] (12/5,0)--node[inner sep=2pt,scale=0.8,blue,right] {3} (12/5,3/5); \end{tikzpicture} }&\cbox{\begin{tikzpicture}[xscale=0.6,yscale=0.8] \draw[thick,shorten >=2] (-1/3,0)--(23/6,0) node[inner sep=2pt,blue,above left,scale=0.8] {6}; \draw[shorten <=-3pt] (0,0)--node[inner sep=2pt,scale=0.8,blue,right] {1} (0,3/5); \draw[shorten <=-3pt] (4/5,0)--node[inner sep=2pt,scale=0.8,blue,right] {1} (4/5,3/5); \draw[shorten <=-3pt, shorten >=-3pt] (11/5,0)--node[inner sep=2pt,scale=0.8,blue,right] {4} (11/5,3/5); \draw[shorten <=-3pt] (11/5,3/5)--node[inner sep=2pt,scale=0.8,blue,above] {2} (8/5,3/5); \end{tikzpicture} }\\[8pt]Type \redtype{IV} & Type \redtype{D_{g1}} & Type \redtype{[3]I^*_0} & Type \redtype{6^{1,1,4}}
\end{tabular}
\end{center}

\subsection{Loops and D-tails}

In \ref{redlabel},

\begin{itemize}[leftmargin=\parindent]
\item
$\loops$ is a sequence $\{l-l'\}\,n,...$, one for every loop $(l,l')\in\LL$, of depth $n$ in $\G$,\\
sorted by $\gcd(l,m), l, l'$ (and ordering $l\le l'$),\\
$\{l-l'\}$ is omitted if $l, l'$ is the smallest pair in $\cL$ with $\gcd(l,m)=\gcd(l',m)$
unused by the previous loops;
\item
$\dtails$ is a sequence $\{l\}\,n\,{\rm D},...$, one for every D-tail $l\in\LD$, of depth $n$ in $\G$,\\
sorted by $\gcd(l,m)$, $l$, $n$,\\
$\{l\}$ is omitted if it is the smallest even index in $\cL$ unused by the previous loops and D-tails.
\item
depth $n$ in an inner type $m\edge{l}{l'}{n}m'$ is omitted in $\loops$ and $\dtails$
if it is minimal possible among $m\edge{l}{l'}{*}m'$, except when $l,l'$ has been omitted
(so that the loop label is non-empty),
\item
For compatibility with Kodaira-N\'eron classification,
$$
  \rm D_{\cdots,\emph{n\,}\{2\}D,\cdots} \qquad\text{is written}\qquad
  \rm I^*_{\emph{n\,},\cdots\>\cdots},
$$
and similarly for its multiples.
\end{itemize}

\noindent
Here are some reduction types with loops and D-tails, and their labels:

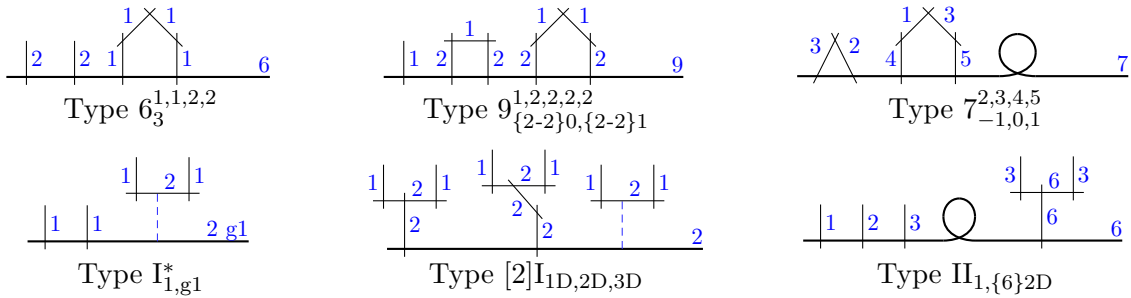
\begin{center}
\begin{tabular}{c@{\qquad\quad}c@{\qquad\quad}c@{\qquad\quad}c}
\cbox{\begin{tikzpicture}[xscale=0.8,yscale=0.8] \draw[thick,shorten >=2] (-1/3,0)--(62/15,0) node[inner sep=2pt,blue,above left,scale=0.8] {6}; \draw[shorten <=-3pt] (0,0)--node[inner sep=2pt,scale=0.8,blue,right] {2} (0,3/5); \draw[shorten <=-3pt] (4/5,0)--node[inner sep=2pt,scale=0.8,blue,right] {2} (4/5,3/5); \draw[shorten <=-3pt, shorten >=-3pt] (8/5,0)--node[inner sep=2pt,scale=0.8,blue,left] {1} (8/5,3/5); \draw[shorten <=-3pt, shorten >=-3pt] (8/5,3/5)--node[inner sep=2pt,scale=0.8,blue,above left=0pt and -1pt] {1} (41/20,21/20); \draw[shorten <=-3pt, shorten >=-3pt] (5/2,3/5)--node[inner sep=2pt,scale=0.8,blue,above right=0pt and -1pt] {1} (41/20,21/20); \draw[shorten <=-3pt, shorten >=-3pt] (5/2,0)--node[inner sep=2pt,scale=0.8,blue,right] {1} (5/2,3/5); \end{tikzpicture} }&\cbox{\begin{tikzpicture}[xscale=0.8,yscale=0.8] \draw[thick,shorten >=2] (-1/3,0)--(71/15,0) node[inner sep=2pt,blue,above left,scale=0.8] {9}; \draw[shorten <=-3pt] (0,0)--node[inner sep=2pt,scale=0.8,blue,right] {1} (0,3/5); \draw[shorten <=-3pt, shorten >=-3pt] (4/5,0)--node[inner sep=2pt,scale=0.8,blue,left] {2} (4/5,3/5); \draw[shorten <=-3pt, shorten >=-3pt] (4/5,3/5)--node[inner sep=2pt,scale=0.8,blue,above] {1} (7/5,3/5); \draw[shorten <=-3pt, shorten >=-3pt] (7/5,0)--node[inner sep=2pt,scale=0.8,blue,right] {2} (7/5,3/5); \draw[shorten <=-3pt, shorten >=-3pt] (11/5,0)--node[inner sep=2pt,scale=0.8,blue,left] {2} (11/5,3/5); \draw[shorten <=-3pt, shorten >=-3pt] (11/5,3/5)--node[inner sep=2pt,scale=0.8,blue,above left=0pt and -1pt] {1} (53/20,21/20); \draw[shorten <=-3pt, shorten >=-3pt] (31/10,3/5)--node[inner sep=2pt,scale=0.8,blue,above right=0pt and -1pt] {1} (53/20,21/20); \draw[shorten <=-3pt, shorten >=-3pt] (31/10,0)--node[inner sep=2pt,scale=0.8,blue,right] {2} (31/10,3/5); \end{tikzpicture} }&\cbox{\begin{tikzpicture}[xscale=0.8,yscale=0.8] \draw[thick,shorten >=2] (-0.33,0)--(5.26,0) node[inner sep=2pt,blue,above left,scale=0.8] {7}; \draw[shorten <=-3pt, shorten >=-3pt] (3/5,0)--node[inner sep=2pt,scale=0.8,blue,above right] {2} (3/10,3/5); \draw[shorten <=-3pt, shorten >=-3pt] (0,0)--node[inner sep=2pt,scale=0.8,blue,above left] {3} (3/10,3/5); \draw[shorten <=-3pt, shorten >=-3pt] (7/5,0)--node[inner sep=2pt,scale=0.8,blue,left] {4} (7/5,3/5); \draw[shorten <=-3pt, shorten >=-3pt] (7/5,3/5)--node[inner sep=2pt,scale=0.8,blue,above left=0pt and -1pt] {1} (37/20,21/20); \draw[shorten <=-3pt, shorten >=-3pt] (23/10,3/5)--node[inner sep=2pt,scale=0.8,blue,above right=0pt and -1pt] {3} (37/20,21/20); \draw[shorten <=-3pt, shorten >=-3pt] (23/10,0)--node[inner sep=2pt,scale=0.8,blue,right] {5} (23/10,3/5); \path[draw,thick] (91/30,0) edge[white,line width=2] ++(0.6,0) edge[out=0,in=-90] ++(0.6,0.4) ++(0.6,0.4) edge[out=90,in=0] ++(-0.3,0.3) ++(-0.3,0.3) edge[out=180,in=90] ++(-0.3,-0.3) ++(-0.3,-0.3) edge[out=-90,in=180] ++(0.6,-0.4); \end{tikzpicture} }\\[12pt]Type \redtype{6^{1,1,2,2}_3} & Type \redtype{9^{1,2,2,2,2}_{\{2\m2\}0,\{2\m2\}1}} & Type \redtype{7^{2,3,4,5}_{-1,0,1}}
\\[8pt]
\cbox{\begin{tikzpicture}[xscale=0.7,yscale=0.8] \draw[thick,shorten >=2] (-1/3,0)--(119/30,0) node[inner sep=2pt,blue,above left,scale=0.8] {2 \smash{g}1}; \draw[shorten <=-3pt] (0,0)--node[inner sep=2pt,scale=0.8,blue,right] {1} (0,3/5); \draw[shorten <=-3pt] (4/5,0)--node[inner sep=2pt,scale=0.8,blue,right] {1} (4/5,3/5); \draw[blue!80!white,densely dashed] (32/15,0)--(32/15,4/5);\draw (23/15,4/5)--node[inner sep=2pt,scale=0.8,blue,above,pos=0.6] {2} (44/15,4/5); \draw[shorten <=-3pt] (26/15,4/5)--node[inner sep=2pt,scale=0.8,blue,left] {1} (26/15,7/5); \draw[shorten <=-3pt] (41/15,4/5)--node[inner sep=2pt,scale=0.8,blue,right] {1} (41/15,7/5); \end{tikzpicture} }&\cbox{\begin{tikzpicture}[xscale=0.7,yscale=0.8] \draw[thick,shorten >=2] (4/15,0)--(191/30,0) node[inner sep=2pt,blue,above left,scale=0.8] {2}; \draw[shorten <=-3pt, shorten >=-3pt] (3/5,0)--node[inner sep=2pt,scale=0.8,blue,right] {2} (3/5,4/5); \draw (0,4/5)--node[inner sep=2pt,scale=0.8,blue,above,pos=0.6] {2} (7/5,4/5); \draw[shorten <=-3pt] (1/5,4/5)--node[inner sep=2pt,scale=0.8,blue,left] {1} (1/5,7/5); \draw[shorten <=-3pt] (6/5,4/5)--node[inner sep=2pt,scale=0.8,blue,right] {1} (6/5,7/5); \draw[shorten <=-3pt, shorten >=-3pt] (187/60,0)--node[inner sep=2pt,scale=0.8,blue,right] {2} (187/60,3/5); \draw[shorten <=-3pt, shorten >=-3pt] (187/60,3/5)--node[inner sep=2pt,scale=0.8,blue,below left=-1pt and -1pt] {2} (8/3,21/20); \draw (31/15,21/20)--node[inner sep=2pt,scale=0.8,blue,above,pos=0.6] {2} (52/15,21/20); \draw[shorten <=-3pt] (34/15,21/20)--node[inner sep=2pt,scale=0.8,blue,left] {1} (34/15,33/20); \draw[shorten <=-3pt] (49/15,21/20)--node[inner sep=2pt,scale=0.8,blue,right] {1} (49/15,33/20); \draw[blue!80!white,densely dashed] (71/15,0)--(71/15,4/5);\draw (62/15,4/5)--node[inner sep=2pt,scale=0.8,blue,above,pos=0.6] {2} (83/15,4/5); \draw[shorten <=-3pt] (13/3,4/5)--node[inner sep=2pt,scale=0.8,blue,left] {1} (13/3,7/5); \draw[shorten <=-3pt] (16/3,4/5)--node[inner sep=2pt,scale=0.8,blue,right] {1} (16/3,7/5); \end{tikzpicture} }&\cbox{\begin{tikzpicture}[xscale=0.7,yscale=0.8] \draw[thick,shorten >=2] (-0.33,0)--(5.83,0) node[inner sep=2pt,blue,above left,scale=0.8] {6}; \draw[shorten <=-3pt] (0,0)--node[inner sep=2pt,scale=0.8,blue,right] {1} (0,3/5); \draw[shorten <=-3pt] (4/5,0)--node[inner sep=2pt,scale=0.8,blue,right] {2} (4/5,3/5); \draw[shorten <=-3pt] (8/5,0)--node[inner sep=2pt,scale=0.8,blue,right] {3} (8/5,3/5); \path[draw,thick] (7/3,0) edge[white,line width=2] ++(0.6,0) edge[out=0,in=-90] ++(0.6,0.4) ++(0.6,0.4) edge[out=90,in=0] ++(-0.3,0.3) ++(-0.3,0.3) edge[out=180,in=90] ++(-0.3,-0.3) ++(-0.3,-0.3) edge[out=-90,in=180] ++(0.6,-0.4); \draw[shorten <=-3pt, shorten >=-3pt] (4.2,0)--node[inner sep=2pt,scale=0.8,blue,right] {6} (4.2,4/5); \draw (3.6,4/5)--node[inner sep=2pt,scale=0.8,blue,above,pos=0.6] {6} (5,4/5); \draw[shorten <=-3pt] (3.8,4/5)--node[inner sep=2pt,scale=0.8,blue,left] {3} (3.8,7/5); \draw[shorten <=-3pt] (4.8,4/5)--node[inner sep=2pt,scale=0.8,blue,right] {3} (4.8,7/5); \end{tikzpicture} }\\[16pt]Type \redtype{I^*_{1,g1}} & Type \redtype{[2]I_{1D,2D,3D}} & Type \redtype{II_{1,\{6\}2D}}
\end{tabular}
\end{center}

\subsection{Inner chains}

In \ref{redlabel},

\begin{itemize}[leftmargin=\parindent]
\item
Inner chains between components $\G$ of multiplicity $m$ and $\G'$ of multiplicity $m'$
are labelled with their inner types $\edge{d}{d'}{n}$, and sorted by $\gcd(d,m)\>(\!=\!\gcd(d',m')), d, d'$,\\
$d-d'$ is omitted if $l, l'$ is the smallest pair in $\cL\times\cL'$ with $\gcd(d,m)=\gcd(d',m)$
unused by the loops and D-tails of $\G$ and $\G'$;\\
$n$ is omitted it is minimal possible among $m\edge{d}{d'}{*}m'$.
\end{itemize}

\noindent
Here are some reduction types with inner chains and their labels (see also Example \ref{ex23chain}),

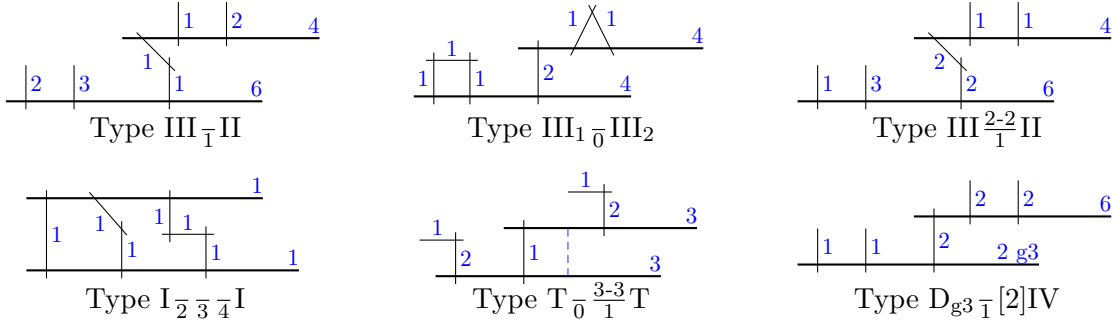
\begin{center}
\begin{tabular}{c@{\qquad\quad}c@{\qquad\quad}c@{\qquad\quad}c}
\cbox{\begin{tikzpicture}[xscale=0.8,yscale=0.8] \draw[thick,shorten >=2] (-1/3,0)--(241/60,0) node[inner sep=2pt,blue,above left,scale=0.8] {6}; \draw[shorten <=-3pt] (0,0)--node[inner sep=2pt,scale=0.8,blue,right] {2} (0,3/5); \draw[shorten <=-3pt] (4/5,0)--node[inner sep=2pt,scale=0.8,blue,right] {3} (4/5,3/5); \draw[shorten <=-3pt, shorten >=-3pt] (143/60,0)--node[inner sep=2pt,scale=0.8,blue,right] {1} (143/60,3/5); \draw[shorten <=-3pt, shorten >=-3pt] (143/60,3/5)--node[inner sep=2pt,scale=0.8,blue,below left=-1pt and -1pt] {1} (29/15,21/20); \draw[thick,shorten >=2] (8/5,21/20)--(149/30,21/20) node[inner sep=2pt,blue,above left,scale=0.8] {4}; \draw[shorten <=-3pt] (38/15,21/20)--node[inner sep=2pt,scale=0.8,blue,right] {1} (38/15,33/20); \draw[shorten <=-3pt] (10/3,21/20)--node[inner sep=2pt,scale=0.8,blue,right] {2} (10/3,33/20); \end{tikzpicture} }&\cbox{\begin{tikzpicture}[xscale=0.8,yscale=0.8] \draw[thick,shorten >=2] (-1/3,0)--(101/30,0) node[inner sep=2pt,blue,above left,scale=0.8] {4}; \draw[shorten <=-3pt, shorten >=-3pt] (0,0)--node[inner sep=2pt,scale=0.8,blue,left] {1} (0,3/5); \draw[shorten <=-3pt, shorten >=-3pt] (0,3/5)--node[inner sep=2pt,scale=0.8,blue,above] {1} (3/5,3/5); \draw[shorten <=-3pt, shorten >=-3pt] (3/5,0)--node[inner sep=2pt,scale=0.8,blue,right] {1} (3/5,3/5); \draw[shorten <=-3pt, shorten >=-3pt] (26/15,0)--node[inner sep=2pt,scale=0.8,blue,right] {2} (26/15,4/5); \draw[thick,shorten >=2] (7/5,4/5)--(137/30,4/5) node[inner sep=2pt,blue,above left,scale=0.8] {4}; \draw[shorten <=-3pt, shorten >=-3pt] (44/15,4/5)--node[inner sep=2pt,scale=0.8,blue,above right] {1} (79/30,7/5); \draw[shorten <=-3pt, shorten >=-3pt] (7/3,4/5)--node[inner sep=2pt,scale=0.8,blue,above left] {1} (79/30,7/5); \end{tikzpicture} }&\cbox{\begin{tikzpicture}[xscale=0.8,yscale=0.8] \draw[thick,shorten >=2] (-1/3,0)--(241/60,0) node[inner sep=2pt,blue,above left,scale=0.8] {6}; \draw[shorten <=-3pt] (0,0)--node[inner sep=2pt,scale=0.8,blue,right] {1} (0,3/5); \draw[shorten <=-3pt] (4/5,0)--node[inner sep=2pt,scale=0.8,blue,right] {3} (4/5,3/5); \draw[shorten <=-3pt, shorten >=-3pt] (143/60,0)--node[inner sep=2pt,scale=0.8,blue,right] {2} (143/60,3/5); \draw[shorten <=-3pt, shorten >=-3pt] (143/60,3/5)--node[inner sep=2pt,scale=0.8,blue,below left=-1pt and -1pt] {2} (29/15,21/20); \draw[thick,shorten >=2] (8/5,21/20)--(149/30,21/20) node[inner sep=2pt,blue,above left,scale=0.8] {4}; \draw[shorten <=-3pt] (38/15,21/20)--node[inner sep=2pt,scale=0.8,blue,right] {1} (38/15,33/20); \draw[shorten <=-3pt] (10/3,21/20)--node[inner sep=2pt,scale=0.8,blue,right] {1} (10/3,33/20); \end{tikzpicture} }\\[16pt]Type \redtype{III\e(1)II} & Type \redtype{III_1\e(0)III_2} & Type \redtype{III\e{2-2}(1)II}
\\[8pt]
\cbox{\begin{tikzpicture}[xscale=0.8,yscale=0.8] \draw[thick,shorten >=2] (-1/3,0)--(257/60,0) node[inner sep=2pt,blue,above left,scale=0.8] {1}; \draw[shorten <=-3pt, shorten >=-3pt] (0,0)--node[inner sep=2pt,scale=0.8,blue,right] {1} (0,6/5); \draw[shorten <=-3pt, shorten >=-3pt] (5/4,0)--node[inner sep=2pt,scale=0.8,blue,right] {1} (5/4,24/35); \draw[shorten <=-3pt, shorten >=-3pt] (5/4,24/35)--node[inner sep=2pt,scale=0.8,blue,below left=-1pt and -1pt] {1} (4/5,6/5); \draw[shorten <=-3pt, shorten >=-3pt] (53/20,0)--node[inner sep=2pt,scale=0.8,blue,right] {1} (53/20,3/5); \draw[shorten <=-3pt, shorten >=-3pt] (53/20,3/5)--node[inner sep=2pt,scale=0.8,blue,above] {1} (41/20,3/5); \draw[shorten <=-3pt, shorten >=-3pt] (41/20,3/5)--node[inner sep=2pt,scale=0.8,blue,left] {1} (41/20,6/5); \draw[thick,shorten >=2] (-1/3,6/5)--(221/60,6/5) node[inner sep=2pt,blue,above left,scale=0.8] {1}; \end{tikzpicture} }&\cbox{\begin{tikzpicture}[xscale=0.8,yscale=0.8] \draw[thick,shorten >=2] (4/15,0)--(41/10,0) node[inner sep=2pt,blue,above left,scale=0.8] {3}; \draw[shorten <=-3pt, shorten >=-3pt] (3/5,0)--node[inner sep=2pt,scale=0.8,blue,right] {2} (3/5,3/5); \draw[shorten <=-3pt] (3/5,3/5)--node[inner sep=2pt,scale=0.8,blue,above] {1} (0,3/5); \draw[shorten <=-3pt, shorten >=-3pt] (26/15,0)--node[inner sep=2pt,scale=0.8,blue,right] {1} (26/15,4/5); \draw[blue!80!white,densely dashed] (37/15,0)--(37/15,4/5);\draw[thick,shorten >=2] (7/5,4/5)--(47/10,4/5) node[inner sep=2pt,blue,above left,scale=0.8] {3}; \draw[shorten <=-3pt, shorten >=-3pt] (46/15,4/5)--node[inner sep=2pt,scale=0.8,blue,right] {2} (46/15,7/5); \draw[shorten <=-3pt] (46/15,7/5)--node[inner sep=2pt,scale=0.8,blue,above] {1} (37/15,7/5); \end{tikzpicture} }&\cbox{\begin{tikzpicture}[xscale=0.8,yscale=0.8] \draw[thick,shorten >=2] (-1/3,0)--(113/30,0) node[inner sep=2pt,blue,above left,scale=0.8] {2 \smash{g}3}; \draw[shorten <=-3pt] (0,0)--node[inner sep=2pt,scale=0.8,blue,right] {1} (0,3/5); \draw[shorten <=-3pt] (4/5,0)--node[inner sep=2pt,scale=0.8,blue,right] {1} (4/5,3/5); \draw[shorten <=-3pt, shorten >=-3pt] (29/15,0)--node[inner sep=2pt,scale=0.8,blue,right] {2} (29/15,4/5); \draw[thick,shorten >=2] (8/5,4/5)--(149/30,4/5) node[inner sep=2pt,blue,above left,scale=0.8] {6}; \draw[shorten <=-3pt] (38/15,4/5)--node[inner sep=2pt,scale=0.8,blue,right] {2} (38/15,7/5); \draw[shorten <=-3pt] (10/3,4/5)--node[inner sep=2pt,scale=0.8,blue,right] {2} (10/3,7/5); \end{tikzpicture} }\\[16pt]Type \redtype{I\e(2)\e(3)\e(4)I} & Type \redtype{T\e(0)\e{3-3}(1)T} & Type \redtype{D_{g3}\e(1)[2]IV}
\end{tabular}
\end{center}

\subsection{Families}

Reduction family is denoted by its reduction type with minimal depths (\S\ref{intronotation}).

\subsection{Jumping and revisting components}

Here are a few examples of Notation \ref{redlabel} when the graph $G$ is not Eulerian (so that jumps `$\&$'
between vertices are needed), and/or components revisited (marked by $c_i$ if it is the $i$th principal type
printed in the label).

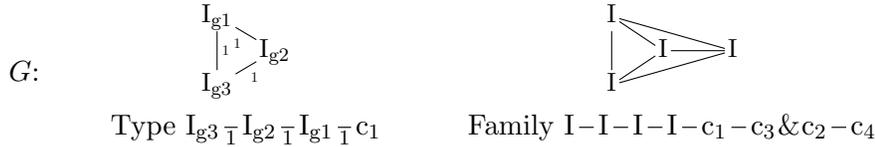
\begin{center}
$G$:\qquad
\begin{tabular}{c@{\qquad\quad}c@{\qquad\quad}c@{\qquad\quad}c}
\cbox{\begin{tikzpicture}[xscale=1.2,yscale=1.2, sup/.style={midway,auto,scale=0.5}, sub/.style={sup,swap}, lrg/.style={scale=0.9,inner sep=0.1em,text depth=0.5ex,text height=1.6ex}] \node[lrg] at (0.64,0.64) (1) {$\redtype{I_{g3}}$}; \node[lrg] at (1.28,1.02) (2) {$\redtype{I_{g2}}$}; \node[lrg] at (0.64,1.4) (3) {$\redtype{I_{g1}}$}; \draw (1) edge node [sub] {1} (3); \draw (3) edge node [sub] {1} (2); \draw (1) edge node [sub] {1} (2); \end{tikzpicture} }&\cbox{\begin{tikzpicture}[xscale=1.2,yscale=1.2, sup/.style={midway,auto,scale=0.5}, sub/.style={sup,swap}, lrg/.style={scale=0.9,inner sep=0.1em,text depth=0.5ex,text height=1.6ex}] \node[lrg] at (0.64,0.64) (1) {$\redtype{I}$}; \node[lrg] at (1.21,1.02) (2) {$\redtype{I}$}; \node[lrg] at (0.64,1.4) (3) {$\redtype{I}$}; \node[lrg] at (1.98,1.02) (4) {$\redtype{I}$}; \draw (1) edge (3); \draw (1) edge (4); \draw (2) edge (4); \draw (3) edge (2); \draw (1) edge (2); \draw (3) edge (4); \end{tikzpicture} }\\[16pt]Type \redtype{I_{g3}\e(1)I_{g2}\e(1)I_{g1}\e(1)c_1} & Family \redtype{I\e I\e I\e I\e c_1\e c_3\&c_2\e c_4}
\end{tabular}
\end{center}

\subsection{Standalone principal types}

When the component type is not labelled as part of the graph but we wish to refer to it in a standalone
way, we represent its edges
\begin{itemize}
\item
as a sequence $\frac{\>l\>}{}$ following the component type
for every edge $l$.
\end{itemize}
See Tables \ref{tabPS1}--\ref{tabPS6} on pp.\pageref{tabPS1}--\pageref{tabPS6} for examples.

\subsection{Canonical label}
\label{canlabel}

To find the canonical label, we attach scores to principal types (vertices of the graph $G$ associated to $R$),
inner chains (edges of $G$) and also vertices and edges of the shape $S$ associated to $R$.
In all cases, they are sequences of integers ordered lexicographically:

\begin{itemize}
\item[(a)]
Score of a principal type \mgOLfull{} in $G$:
$$
  (\chi_\G,m,-g,|\LM|,|\LD|,|\LL|,|\cO|,\cO,\LL,\LD,\LM, \text{depths of loops}, \text{depths of D-tails}).
$$
\item[(b)]
Score of a (directed) edge of inner type $m\edge{d}{d'}nm'$ in $G$:
$$
  (\gcd(m,d)=\gcd(m',d'), d, d', n)
$$
\item[(c)]
Score of a vertex of $S$ is the Euler characteristic $\chi$ of the corresponding principal type.
\item[(d)]
Score of an edge of $S$ is the sequence of weights of all inner chains between the principal types
at the two vertices.
\end{itemize}

Then we find the minimal path in $S$ (see \ref{minpath}) and all ways to traverse it. Going back to $G$, we
find the lexicographically smallest path in $G$ among them. We then follow it to get the label.

\noindent\par\medskip\noindent

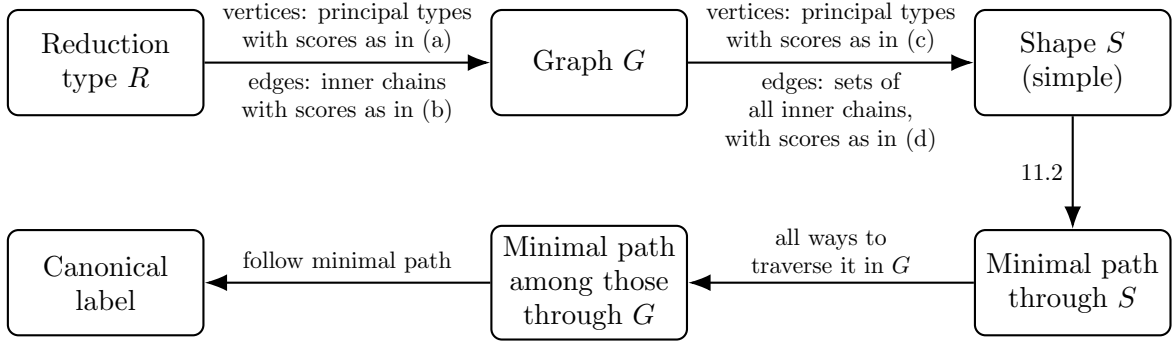
\begin{tikzpicture}[
  node distance=1.5cm and 3.8cm,
  obj/.style={draw, rounded corners, minimum height=1.4cm, minimum width=2.6cm, align=center},
  every path/.style={draw, -{Latex[length=3mm]}, thick},
  alabel/.style={midway, auto, swap, font=\small, scale=0.9, align=center},
  blabel/.style={alabel, swap}
]

\node[obj] (R) {Reduction\\type $R$};
\node[obj] (G) [right=of R] {Graph $G$};
\node[obj] (S) [right=of G] {Shape $S$\\(simple)};
\node[obj] (MS) [below=of S] {Minimal path\\through $S$};
\node[obj] (MG) [left=of MS] {Minimal path\\among those\\through $G$};
\node[obj] (CL) [left=of MG] {Canonical\\label};

\draw (R)  -- node[blabel] {vertices: principal types\\with scores as in (a)}
              node[alabel] {edges: inner chains\\with scores as in (b)} (G);
\draw (G)  -- node[blabel] {vertices: principal types\\with scores as in (c)}
              node[alabel] {edges: sets of\\all inner chains,\\with scores as in (d)} (S);
\draw (S)  -- node[alabel] {\ref{minpath}}  (MS);
\draw (MS) -- node[alabel] {all ways to\\traverse it in $G$}       (MG);
\draw (MG) -- node[alabel] {follow minimal path}            (CL);

\end{tikzpicture}

\begin{example}[Canonical label]
There are ${10+2\choose 3}=220$ families of genus 4 families of shape
$$
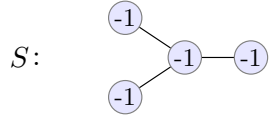
\cbox{$S$}\colon\qquad\cbox{%
\scalebox{1.1}{\begin{tikzpicture}[xscale=0.8,yscale=0.8,auto=left, v/.style={pShapeVertex}, l/.style={pShapeEdge}] \node[v] (1) at (1,1) {-1}; \node[v] (2) at (1.9,1.6) {-1}; \node[v] (3) at (1,2.2) {-1}; \node[v] (4) at (2.9,1.6) {-1}; \draw[] (1)--node[l] {} (2); \draw[] (2)--node[l] {} (3); \draw[] (2)--node[l] {} (4); \end{tikzpicture}}
}
$$

\noindent
A reduction type from these families has a component of multiplicity 1 as a `root' and three Kodaira types as `leaves',
linked with chains of weight 1. Let $R$ be one of those, say

\begin{center}
\cbox{
\begin{tikzpicture}[xscale=1.2,yscale=1.2, sup/.style={midway,auto,scale=0.5}, sub/.style={sup,swap}, lrg/.style={scale=0.9,inner sep=0.1em,text depth=0.5ex,text height=1.6ex}] \node[lrg] at (0.8,0.8) (1) {$\redtype{I^*_0}$}; \node[lrg] at (1.52,1.28) (2) {$\redtype{I}$}; \node[lrg] at (0.8,1.76) (3) {$\redtype{II}$}; \node[lrg] at (2.32,1.28) (4) {$\redtype{III^*}$}; \draw (2) edge node [sub] {0} (4); \draw (3) edge node [sub] {1} (2); \draw (1) edge node [sub] {1} (2); \end{tikzpicture}
}\qquad\qquad\cbox{
\begin{tikzpicture}[xscale=0.7,yscale=0.8] \draw[thick,shorten >=2] (-1/3,0)--(329/30,0) node[inner sep=2pt,blue,above left,scale=0.8] {1} node[inner sep=2pt,below left,scale=0.8] {$\Gamma_2$}; \draw[shorten <=-3pt, shorten >=-3pt] (0,0)--node[inner sep=2pt,scale=0.8,blue,right] {1} (0,4/5); \draw[thick,shorten >=2] (-1/3,4/5)--(97/30,4/5) node[inner sep=2pt,blue,above left,scale=0.8] {6} node[inner sep=2pt,below left,scale=0.8] {$\Gamma_3$}; \draw[shorten <=-3pt] (3/5,4/5)--node[inner sep=2pt,scale=0.8,blue,right] {2} (3/5,7/5); \draw[shorten <=-3pt] (7/5,4/5)--node[inner sep=2pt,scale=0.8,blue,right] {3} (7/5,7/5); \draw[shorten <=-3pt, shorten >=-3pt] (301/60,0)--node[inner sep=2pt,scale=0.8,blue,right] {2} (301/60,3/5); \draw[shorten <=-3pt, shorten >=-3pt] (301/60,3/5)--node[inner sep=2pt,scale=0.8,blue,below left=-1pt and -1pt] {3} (137/30,21/20); \draw[thick,shorten >=2] (127/30,21/20)--(39/5,21/20) node[inner sep=2pt,blue,above left,scale=0.8] {4} node[inner sep=2pt,below left,scale=0.8] {$\Gamma_4$}; \draw[shorten <=-3pt, shorten >=-3pt] (31/6,21/20)--node[inner sep=2pt,scale=0.8,blue,right] {3} (31/6,33/20); \draw[shorten <=-3pt, shorten >=-3pt] (31/6,33/20)--node[inner sep=2pt,scale=0.8,blue,above] {2} (137/30,33/20); \draw[shorten <=-3pt] (137/30,33/20)--node[inner sep=2pt,scale=0.8,blue,left] {1} (137/30,9/4); \draw[shorten <=-3pt] (179/30,21/20)--node[inner sep=2pt,scale=0.8,blue,right] {2} (179/30,33/20); \draw[shorten <=-3pt, shorten >=-3pt] (137/15,0)--node[inner sep=2pt,scale=0.8,blue,right] {1} (137/15,4/5); \draw[thick,shorten >=2] (44/5,4/5)--(79/6,4/5) node[inner sep=2pt,blue,above left,scale=0.8] {2} node[inner sep=2pt,below left,scale=0.8] {$\Gamma_1$}; \draw[shorten <=-3pt] (146/15,4/5)--node[inner sep=2pt,scale=0.8,blue,right] {1} (146/15,7/5); \draw[shorten <=-3pt] (158/15,4/5)--node[inner sep=2pt,scale=0.8,blue,right] {1} (158/15,7/5); \draw[shorten <=-3pt] (34/3,4/5)--node[inner sep=2pt,scale=0.8,blue,right] {1} (34/3,7/5); \end{tikzpicture}
}\\
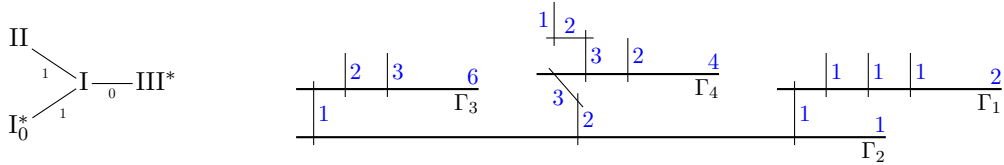
Reduction type $R$ as a graph $G$ and as a special fibre
\end{center}

\noindent
Here is how we determine its canonical label.

First, all vertices of the shape $S$ have score $-1$ ($\chi$'s of principal types) and all edges $(1)$ (weight sequence).
If we write `v' for a new vertex with $\chi=-1$, `-' for an edge of weight (1), and \& for a jump, we find
that the shortest path $P$ is v-v-v\&v-$c_2$. Recall that actual edges are preferred to jumps,
and going to a new vertex preferred to revisiting an old one. (In general, also vertices with smaller $\chi$ come first,
if possible, as they have smaller scores.) For example,

\smallskip
\noindent
\begin{center}
\begin{tabular}{llllll}
v-v-v\&v-$c_2$ &$<$& v-v\&v-$c_2$-v && (jumps are larger than edge marks)\cr
v-v-v\&v-$c_2$ &$<$& v-v-v\&$c_2$-v && (repeated vertex indices are larger than vertex marks)\cr
\end{tabular}
\end{center}

This path can be used to construct the shape, and determines it up to isomorphism. There are $|\Aut S|=6$ ways to trail $S$ in accordance with this path, and as far the shape is concerned, they are completely identical.
This gives six possible labels for our reduction family that all traverse the shape according to path $P$:

\begin{center}
\redtype{I^*_0}-\redtype{I}-\redtype{II}\&\redtype{III^*}-$c_2$\quad\redtype{I^*_0}-\redtype{I}-\redtype{III^*}\&\redtype{II}-$c_2$\quad\redtype{II}-\redtype{I}-\redtype{I^*_0}\&\redtype{III^*}-$c_2$\quad\redtype{II}-\redtype{I}-\redtype{III^*}\&\redtype{I^*_0}-$c_2$\quad\redtype{III^*}-\redtype{I}-\redtype{II}\&\redtype{I^*_0}-$c_2$\quad\redtype{III^*}-\redtype{I}-\redtype{I^*_0}\&\redtype{II}-$c_2$
\end{center}

Now we use (a),(b) above to assign scores to vertices and edges that characterise the actual principal types
(rather than just their $\chi$) and inner chains (rather than just their weight):

\begin{center}
\cbox{
$\begin{array}{ll}
\text{Component} & \text{score}\cr
\redtype{I^*_0\e{1}} & (-1, 2, 0, 1, 0, 0, 3, 1, 1, 1, 1)\cr
\redtype{II\e{1}} & (-1, 6, 0, 1, 0, 0, 2, 2, 3, 1)\cr
\redtype{III^*\e{3}} & (-1, 4, 0, 1, 0, 0, 2, 3, 2, 3)\cr
\redtype{I\e{1}\e{1}\e{1}} & (-1, 1, 0, 3, 0, 0, 0, 1, 1, 1)\cr
\end{array}$
}
\qquad
\qquad
\cbox{
$\begin{array}{ll}
\text{Edge} & \text{score}\cr
\redtype{I}-\redtype{I^*_0} & (1, 1, 1) \cr
\redtype{I}-\redtype{II} & (1, 1, 1) \cr
\redtype{I}-\redtype{III^*} & (1, 3, 0) \cr
\end{array}$
}
\end{center}
The component score is $(\chi,m,-g,...)$, so when all components have the same $\chi$ like in this example, the ones with smaller multiplicity $m$ have smaller score. Because $m(\II)=6$, $m(\IIIS)=4$, $m(\IZS)=2$, the
first two labels are preferred to the last four. They both start with a component \IZS, then an edge \IZS--1 and a component I. After that they differ in that the first one traverses an edge I--II and the second one I--\IIIS. The edge score is smaller for I--II because of the incoming multiplicity, so the first one is the minimal path, and it determines the label for $R$:
\begin{center}
\redtype{I^*_0\e(1)I\e(1)II\&III^*\e(0)c_2}
\end{center}
\end{example}

\section{Counts}
\label{ssprincounts}

As the computation for $2\le g\le 6$ suggests (\ref{numfamilies}), the number of reduction families
in genus $g$ grows fast with $g$.
We show that it is at least $10^g$ for every $g$ (\ref{thm10g}), grows at least like $g^{g-o(g)}$
as $g\to\infty$ (\ref{thmgg}), and analyse the shapes with most families (\ref{fixedfamcount}).
The $10^g$ bound exploits a small number of shapes that have `almost'
$10^g$ families for each one of them.
The superexponential growth $g^{g-o(g)}$ is caused by the number of shapes, which becomes
dominant when $g$ is large.

\begin{theorem}
\label{corecount}
Fix $B\ge 0$. The number $N_B$ of cores $\Psi$ with $\chi(\Psi)\ge -B$ satisfies
$$
  N_B \le 16^{\sqrt B}.
$$
\end{theorem}

\begin{proof}
Recall that
a core requires $k$ non-zero elements $o_1,\dots,o_k\in \Z/m\Z$ with
$\sum_{i=1}^k o_i\equiv 0\pmod m$. Choose a multiset (unordered, repetitions allowed) $o_1,...,o_{k-1}$
from $\{1,\dots,m-1\}$ arbitrarily, and define the final element $o_k$ by
$$
  o_k \equiv -\sum_{i=1}^{k-1} o_i \pmod m.
$$
It yields a valid choice when $o_k\not\equiv 0\pmod m$.
The number of multisets of size $k\!-\!1$ from an $m\!-\!1$-element set is
$\binom{(m-1)+(k-1)-1}{k-1}=\binom{m+k-3}{k-1}$, and this is a bound
on the number $N_{m,k}$ of cores with given $m$ and $k$.

Denote $C=2B+8$ and apply the bound \ref{coremain}(4) for $m$ in terms of $B$. We get
$$
  N_B = \sum_{k=3}^\infty \sum_{m=2}^{\lfloor \frac{C}{k-2}\rfloor} N_{m,k} \le
    \sum_{k=3}^\infty \sum_{m=2}^{\lfloor \frac{C}{k-2}\rfloor} \binom{m+k-3}{k-1}.
$$
By a standard identity $\sum_{n=a}^{b} \binom{n}{a} = \binom{b+1}{a+1}$, the inner sums collapses to
one binomial coefficient:
\begin{equation}\label{nbbound1}
  N_B \>\>\le\>\> \sum_{k=3}^\infty \binom{\lfloor \frac{C}{k-2}\rfloor+k-2}{k} =
    \sum_{l=1}^{\lfloor\frac C2\rfloor} \binom{\lfloor \frac{C}{l}\rfloor+l}{l+2}.
\end{equation}
Compute the terms with $l=1,2,3$ and estimate the others
using the bound $\binom{a+b}{b} \le 4^{\sqrt{a b}}$ with $a = \lfloor C/l \rfloor$ and $b = l + 2$,
using that $(l+2)/l \le 1.5$ for $l \ge 4$:
$$
  \binom{\lfloor C/l \rfloor + l}{\,l+2\,}
     \le 4^{\sqrt{\lfloor C/l \rfloor\,(l+2)}}
     \le 4^{\sqrt{1.5\,C}}.
$$
There are $\le C/2$ terms with $l\ge 4$, and we get
$$
  N_B \le \text{\small$\displaystyle
   \binom{C+1}{3} +
   \binom{\lfloor C/2 \rfloor + 2}{4} +
   \binom{\lfloor C/3 \rfloor + 3}{5} $} + \frac{C}2 4^{\sqrt{1.5\,C}}.
$$
For large $C$, the polynomial terms are negligible and
$N_B\le 4^{2\sqrt{C/2-4}}=4^{2\sqrt{B}}$ because $\sqrt{1.5}<2/\sqrt{2}$. In
fact, one checks that the right hand-side $\le 4^{2\sqrt{B}}$ for $B>246$,
and for $0<B\le 246$ one verifies directly that the right-hand side of \eqref{nbbound1}
is $\le 4^{2\sqrt{B}}=16^{\sqrt{B}}$.
\end{proof}

\begin{remark}
(a)
The number of cores does grow exponentially in $\sqrt{B}$, at least modulo
the number-theoretic fluctuations (visible in Tables \ref{tabC2}--\ref{tabC-8} and \ref{tabP}). For example,
if $m=p$ is prime, the number $N_{p,k}$ of cores $p^{o_1,...,o_k}$ with given $k\ge 3$ can be shown to be
exactly
$$
  N_{p,k} = \frac 1p\Bigg(\binom{p+k-2}{k}+(p-1)\,\delta(k)\Bigg), \qquad
  \delta(k) = \begin{cases}
    1, & k\equiv 0\pmod p,\\[-3pt]
    -1,& k\equiv 1\pmod p,\\[-3pt]
    0,& \text{otherwise.}
  \end{cases}
$$
When $k=p$, we get $N_{p,p}\approx \binom{2p}{p}\approx 4^p/\sqrt{\pi p}$ (by Stirling's formula) with
$B=-\chi=p(p-2)+p$, so $p\approx\sqrt{B}$
and we get a lower bound $N_B\ge 4^{\sqrt{B}+o(\sqrt{B})}$.

(b)
The count for the number of principal types with $\chi\ge -B$ is similar, but with a twist: it is still
sub-exponential but the growth is $c^{B^{5/6}}$ rather than $c^{B^{1/2}}$.
As in \ref{corecount}, we can see that the number $P_{m,k}$ of principal types with given $m$ and $k=|\cO\cup\cL|$
satisfies
$$
  P_{m,k} \le \binom{4m+k-2}{k}\binom{\frac{m(m+1)}{2}+k-1}{k}. 
$$
The first factor counts the number of ways to select $k$ elements from $\Z/m\Z$ and
decide which ones become outer multiplicities or those of D-tails, loops or
edges. This is the number of multisets of size $k$ in the set $(\Z/m\Z) \times \{O,\text{loops},D,\text{edges}\}$ of size $4m$.
It is again $O(c^{\sqrt{B}})$ for some $c>0$.

The second factor bounds the number of ways to allocates loops by selecting up to $k$ pairs of elements
of $\Z/m\Z$, i.e. counts multisets of size $k$ chosen from a set with $\frac{m(m+1)}{2}$ elements.
If we consider principal types with $k\approx m^2/2$ even and only loops ($\cO=\LD=\LM=\emptyset$) that maximise
this binomial coefficient, the count becomes $c^{B^{5/6}+o(B^{5/6})}$.

(c)
When we count families of a specific shape by allocating principal types with given $\chi$ to its vertices,
the ordering of the edges starts to matter, and this increases the count once more. For example, just
considering principal types $S=\mgOLfull$ with
$$
  m=3, \quad g=0, \quad \cO=\LL=\LD=\emptyset, \quad \LM=\{\text{$k$ \textbf{ordered} elements of $(\Z/3\Z)^\times$}\}.
$$
Then $\chi(S)=6-3k$ and the number of such types is $\approx 2^k/3$ (the factor 3 coming from the sum 0 condition),
so the number of such types is \emph{exponential} in $B=-\chi$.

We now give a precise exponential upper bound
that we will need in order to understand the shapes with largest number of families.
\end{remark}

\begin{proposition}
\label{princountmain}
Fix $\chi\le -1$ and a weight vector $w=(w_1,w_2,...)$ with $w_i\in\N$.
Let $N^\chi_w$ be the number of principal types of Euler charactristic $\chi$ given together
with an ordering of $\LM=(l_1,l_2,...)$ such that the weight vector of $(l_i)$ is $w$.
Then
$$
  N^\chi_w \le \sqrt{10}^{2-\chi-|w|}.
$$
Moreover, $N^\chi_w \le \frac12\sqrt{10}^{2-\chi-|w|}$ unless $\chi$ and $w$ are one of the following:
$$
\begin{array}{c@{\quad}c@{\quad}c@{\quad}c@{\quad}l}
\chi & |w| & w & N^\chi_w & \sqrt{10}^{2-\chi-|w|} \cr
\text{\rm any} & 2\!-\!\chi & (1,...,1) & 1 & 1 \cr
-1 & 1 & (1) & 10 & 10 \cr
-2 & 1 & (2) & 18 & 10^{3/2}\approx 31.6 \cr
-2 & 2 & (1,1) & 8 & 10 \cr
\end{array}
$$
\end{proposition}

\begin{proof}
Let $S=\mgOLfull{}$ be a principal type with $|\cO|\!=\!n$, $|\cL|\!=\!k$ and $B\!=\!-\chi(S)$.
Recall that $m\le 6B$ and $m(n+2k-4)\le 2B$ by \ref{prinbounds}(1,3). If $k>\frac B2+2$,
then $m=1$, $n=0$ and $\LD=\emptyset$ by \ref{prinbounds}(5);
varying $g$ and the number of loops gives all such types, and their
count satisfies the asserted bound (cf. `right half' of the table below). Thus, assume $k\le\frac B2+2$.

For $B\le 31$, we simply compute all $N_w^k$ and check the claims. Table \ref{tabP} lists
$n_b^\chi=\max_{|w|=b} N_w^\chi$ for varying $\chi$ and $b$, with the exceptions to
the stronger bound $n^\chi_b \le \frac12\sqrt{10}^{2-\chi-b}$ in boldface.

\medskip

\noindent
\begin{center}
{\footnotesize
\begin{tabular}{r|c@{ }c@{ }c@{ }c@{ }c@{ }c@{\>}c@{\>}c@{\>}c@{\>}c@{\>}c@{\>}c@{\>}c@{}c@{}c@{}c@{}c@{}c@{}c@{}c@{}c@{}c@{}c@{}c@{}c}
$\chi$ & $n^\chi_{0}$ & $n^\chi_{1}$ & $n^\chi_{2}$ & $n^\chi_{3}$ & $n^\chi_{4}$ & $n^\chi_{5}$ & $n^\chi_{6}$ & $n^\chi_{7}$ & $n^\chi_{8}$ & $n^\chi_{9}$ & $n^\chi_{10}$ & $n^\chi_{11}$ & $n^\chi_{12}$ & $n^\chi_{13}$ & $n^\chi_{14}$ & $n^\chi_{15}$ & $n^\chi_{16}$ & $n^\chi_{17}$ & $n^\chi_{18}$ & $n^\chi_{19}$ & $n^\chi_{20}$ & $n^\chi_{21}$ & $n^\chi_{22}$ & $n^\chi_{23}$ & $n^\chi_{24}$ \\\hline
-1 & 0&\bf 10&1&\bf 1 \\
-2 & 46&\bf 18&\bf 8&1&\bf 1 \\
-3 & 0&39&6&5&0&\bf 1 \\
-4 & 150&29&37&4&3&0&\bf 1 \\
-5 & 0&108&14&23&1&2&0&\bf 1 \\
-6 & 282&85&100&10&11&1&2&0&\bf 1 \\
-7 & 0&193&32&56&2&4&0&2&0&\bf 1 \\
-8 & 466&124&173&24&28&2&4&0&2&0&\bf 1 \\
-9 & 0&294&64&117&6&16&1&3&0&2&0&\bf 1 \\
-10 & 702&259&298&50&91&6&6&1&3&0&2&0&\bf 1 \\
-11 & 0&635&94&188&18&33&2&5&0&3&0&2&0&\bf 1 \\
-12 & 1323&331&705&82&106&18&33&2&5&0&3&0&2&0&\bf 1 \\
-13 & 0&880&178&585&28&48&6&7&1&4&0&3&0&2&0&\bf 1 \\
-14 & 1419&617&856&152&553&28&58&6&7&1&4&0&3&0&2&0&\bf 1 \\
-15 & 0&1246&330&586&60&357&10&54&2&6&0&4&0&3&0&2&0&\bf 1 \\
-16 & 2373&745&1310&256&364&54&123&10&12&2&6&0&4&0&3&0&2&0&\bf 1 \\
-17 & 0&1988&407&1083&100&232&20&103&6&8&1&5&0&4&0&3&0&2&0&\bf 1 \\
-18 & 3447&1263&2400&401&1171&108&264&20&103&6&8&1&5&0&4&0&3&0&2&0&\bf 1 \\
-19 & 0&3215&690&2367&192&1199&62&156&10&13&2&7&0&5&0&4&0&3&0&2&0&\bf 1 \\
-20 & 4754&1642&3789&632&2651&190&1293&62&198&10&13&2&7&0&5&0&4&0&3&0&2&0&\bf 1 \\
-21 & 0&4032&1098&3421&293&1804&82&421&20&188&6&9&1&6&0&5&0&4&0&3&0&2&0&\bf 1 \\
-22 & 5930&2529&4528&980&3271&295&720&82&295&20&18&6&9&1&6&0&5&0&4&0&3&0&2&0&\bf 1 \\
$\cdots$ & $\cdots$\\
-37 & 0&60610&20010&76255&11307&70282&5414&83851&2640&17164&708&14491&105\,&8276&50&43&20&22&6&13&1&10&0&9&$\cdots$
\end{tabular}}\\[5pt]
\textsc{Table} \customlabel{tabP}{P} of $n_b^\chi=\max_{|w|=b} N_w^\chi$ for varying $\chi$ and $b$.
\end{center}

Suppose $B>31$.
For simplicity, we will overcount by ordering everything rather than just edges.
The number of principal types with given $B$, $m$, $n$, $k$ and ordered $\cO$, $\LL$, $\LD$, $\LM$ is
$
  \le \binom{k+2}{2} m^{n+k} 
$
where the binomial coefficient determines splitting of $\cL$ into $\LL$, $\LD$ and $\LM$, and $m^{n+k}$
selects $n+k$ ordered numbers from $\Z/m\Z$. 
The genus $g$ is determined uniquely from the other invariants.
So we need to show that for $B>31$ and $0\le k\le\frac B2+2$, we have
$$
 \sum_{m=1}^{6B}\>\>\>
 \sum_{n=0}^{\lfloor\frac{2B}m\rfloor-2k+4} m^{n+k} \le \frac12\sqrt{10}^{2-\chi-|w|}/{\textstyle \binom{k+2}{2}}.
$$
We can check this by hand for $31<B<120$, and for $B\ge 120$ we will prove that
$$
  X_{B,k}=\frac{\text{number of terms in the double sum$\>\cdot\>$maximum of the terms}}{\text{right-hand side}}
$$
satisfies $X_{B,k}<1$. Indeed,

\begin{itemize}
\item
the number of terms $(m,n)$ in the double sum is less than $6B\cdot 3B=18B^2$.
\item
$
  m^{n+k} \le m^{(\lfloor\frac{2B}m\rfloor-2k+4) + k} = m^{\frac{2B}m-k+4} = (m^{1/m})^{2B} m^{4-k} \le (3^{1/3})^{2B} m^{4-k} = 3^{\frac{2B}3}m^{4-k},
$ so
\begin{itemize}
\item[] Case (a) $m^{n+k}\le 3^{\frac{2B}3}2^{4-k}$ when $m>1$, $k>4$.
\item[] Case (b) $m^{n+k}\le 3^{\frac{2B}3}(6B)^4$ when $m>1$, $k\le 4$.
\item[] Case (c) $m^{n+k}=1$ when $m=1$.
\end{itemize}
\item
$
  \frac 12\sqrt{10}^{2-\chi-|w|}/{\textstyle \binom{k+2}{2}}
      \ge \sqrt{10}^{2+B-k} / {\textstyle \binom{k+2}{2}}
      \ge \sqrt{10}^{2+B-k} / \frac{(\frac B2+4)(\frac B2+3)}2 \ge \sqrt{10}^{2+B-k} / \frac{B^2}2.
$
\end{itemize}
In case (a),
$$
\begin{aligned}
  X_{B,k}
  &\le 18B^2\cdot 3^{\frac{2B}3}2^{4-k} \cdot
  \frac{B^2/2}{\sqrt{10}^{2+B-k}} =
  \frac{18\cdot 2^4}{2\cdot 10} B^4 \Bigl(\frac{\sqrt{10}}2\Bigr)^k \Bigl(\frac{3^{2/3}}{\sqrt{10}}\Bigr)^B\\
  &\le \frac{18\cdot 2^4}{2\cdot 10}\cdot\frac{10}{4} B^4 \Bigl(\frac{\sqrt{10}}2\Bigr)^{B/2} \Bigl(\frac{3^{2/3}}{\sqrt{10}}\Bigr)^B
  = 36\cdot B^4 \cdot \Bigl(\frac{3^{2/3}}{2^{1/2} 10^{1/4}}\Bigr)^B.
\end{aligned}
$$
The right-hand side is $\approx 36B^4\cdot 0.827^B$ and it is $<1$ for $B\ge 120$. Similarly,
$$
\begin{aligned}
  X_{B,k} \le 18B^2 \cdot 3^{\frac{2B}3}\cdot (6B)^4 {\textstyle\binom{4+2}{2}} / {\sqrt{10}^{2+B-4}} &\quad<1\qquad\text{for $B>102$ in case (b),}\\
  X_{B,k} \le 18B^2 \cdot \frac{B^2/2}{\sqrt{10}^{2+B-k}} \le \frac{9 B^4}{\sqrt{10}^{B/2}}\qquad\>\>\, &\quad<1\qquad\text{for $B>26$ in case (c).}\\
\end{aligned}
$$
\end{proof}

\begin{theorem}
\label{fixedfamcount}
The number $N(S)$ of families of a fixed shape $S$ of genus $g$ satisfies
$$
N(S)
\le 10^g \cdot
\begin{cases}
0.512, & \text{if } g\ge 6,\\[-3pt]
0.44,  & \text{if } g=3,4,5,\\[-3pt]
0.55,  & \text{if } g=2.
\end{cases}
$$
In genus $2\le g\le 7$ there is a unique shape that attains these bounds. In every genus $g\ge 6$, there are
$\lceil\frac{(g-5)^2}{6}\rceil$ such shapes (in particular more than one when $g\ge 8$), and they are as follows:
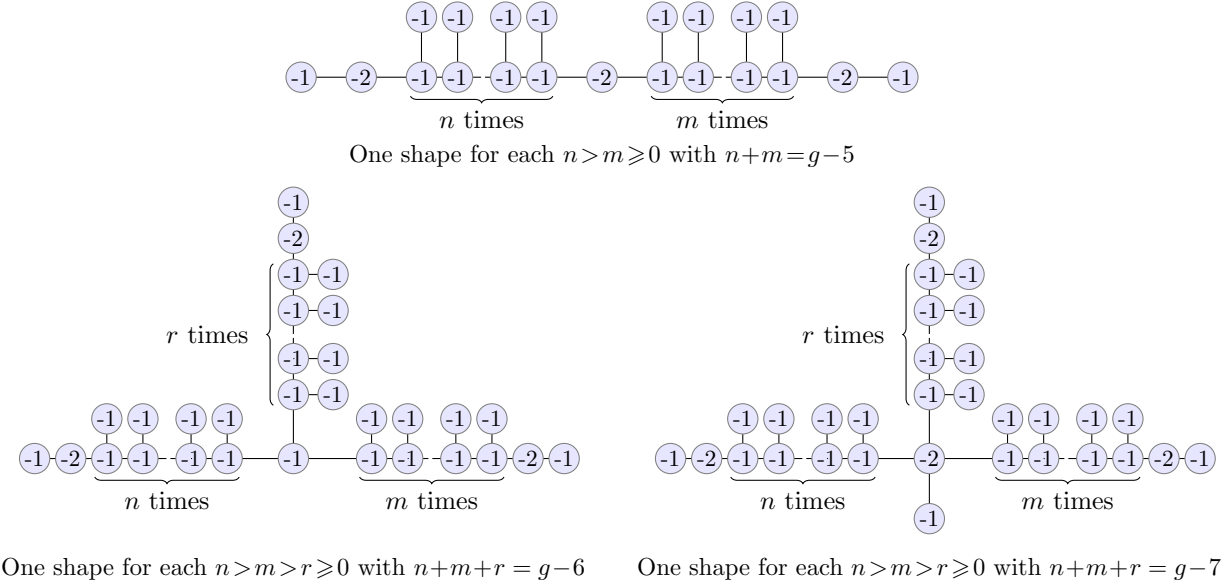
\tikzset{
  comb/.pic={
    \node[v] (1) at (-1.2,0) {-1};
    \node[v] (2) at (-0.6,0) {-2};
    \node[v] (r1) at (0,0) {-1};
    \node[v] (r2) at (0,0.66) {-1};
    \node[v] (r1a) at (0.6,0) {-1};
    \node[v] (r2b) at (0.6,0.66) {-1};
    \node[v] (r3a) at (1.4,0) {-1};
    \node[v] (r4b) at (1.4,0.66) {-1};
    \node[v] (r3) at (2,0) {-1};
    \node[v] (r4) at (2,0.66) {-1};
    \draw (1)--(2)--(r1)--(r1a)--(r2b);
    \draw (r1)--(r2);
    \draw[dash pattern=on 3pt off 1.5pt] (r1a)--(r3a);
    \draw (r4b)--(r3a)--(r3a)--(r3)--(r4);
  }
}
\rm

\begin{center}
\begin{tikzpicture}[
    xscale=0.8,yscale=0.8,v/.style={pShapeVertex},
    dashedline/.style={dash pattern=on 3pt off 1.5pt}
  ]
  \node[v] (a1) at (1,1) {-1};   \node[v] (a2) at (2,1) {-2};   \node[v] (l1) at (3,1) {-1};   \node[v] (l2) at (3,2) {-1};   \node[v] (l1a) at (3.6,1) {-1};   \node[v] (l2b) at (3.6,2) {-1};   \node[v] (l3a) at (4.4,1) {-1};   \node[v] (l4b) at (4.4,2) {-1};   \node[v] (l3) at (5,1) {-1};   \node[v] (l4) at (5,2) {-1};   \node[v] (2) at (6,1) {-2};   \node[v] (r1) at (7,1) {-1};   \node[v] (r2) at (7,2) {-1};   \node[v] (r1a) at (7.6,1) {-1};   \node[v] (r2b) at (7.6,2) {-1};   \node[v] (r3a) at (8.4,1) {-1};   \node[v] (r4b) at (8.4,2) {-1};   \node[v] (r3) at (9,1) {-1};   \node[v] (r4) at (9,2) {-1};   \node[v] (8)  at (10,1) {-2};   \node[v] (9)  at (11,1) {-1};
  \draw (r1)--(2);   \draw (2)--(l3);   \draw (l3)--(l4);   \draw (l1)--(a2);   \draw (l1)--(l2);   \draw (l1)--(l1a);   \draw (l1a)--(l2b);   \draw[dashedline] (l1a)--(l3a);   \draw (l3a)--(l4b);   \draw (l3a)--(l3);   \draw (r3)--(r4);   \draw (r1)--(r2);   \draw (r1)--(r1a);   \draw (r1a)--(r2b);   \draw[dashedline] (r1a)--(r3a);   \draw (r3a)--(r4b);   \draw (r3a)--(r3);   \draw (a2)--(a1);   \draw (r3)--(8);   \draw (8)--(9);   \draw[decorate,decoration={brace,mirror},transform canvas={yshift=-3pt}]
    (l1.south west) -- (l3.south east)
    node[midway,scale=0.9,below=2pt] {$n$ times};   \draw[decorate,decoration={brace,mirror},transform canvas={yshift=-3pt}]
    (r1.south west) -- (r3.south east)
    node[midway,scale=0.9,below=2pt] {$m$ times};
  \node[scale=0.85] at (6,-0.3) {One shape for each $n\!>\!m\!\ge\! 0$ with $n\!+\!m\!=\!g\!-\!5$};
\end{tikzpicture}
\end{center}

\noindent
\begin{tikzpicture}[xscale=0.8,yscale=0.7,v/.style={pShapeVertex},every label/.append style={scale=0.8},baseline=0]

  \pic (C1) at (0.4,0) [scale=0.8] {comb};
  \pic (C2) at (3.5,3.5) [rotate=-90, scale=0.8] {comb};
  \pic (C3) at (6.8,0) [xscale=-1, scale=0.8] {comb};

  \node[v] (M) at (C1r3-|C2r3) {-1};   

  \draw (M)--(C1r3);
  \draw (M)--(C2r3);
  \draw (M)--(C3r3);

  \node[below=1cm of M,scale=0.85] {One shape for each $n\!>\!m\!>\!r\!\ge\!0$ with $n\!+\!m\!+\!r=g\!-\!6$};

  \draw[decorate,decoration={brace,mirror},transform canvas={yshift=-3pt}]
    ($(C1r1)+(-0.2,-0.2)$)--($(C1r3)+(0.2,-0.2)$)
    node[midway,below=2pt,scale=0.9] {$n$ times};

  \draw[decorate,decoration={brace,mirror},transform canvas={yshift=-3pt}]
    ($(C3r3)+(-0.2,-0.2)$)--($(C3r1)+(0.2,-0.2)$)
    node[midway,below=2pt,scale=0.9] {$m$ times};

  \draw[decorate,decoration={brace,mirror},transform canvas={xshift=-3pt}]
    ($(C2r1)+(-0.2,0.2)$)--($(C2r3)+(-0.2,-0.2)$)
    node[midway,left=6pt,scale=0.9] {$r$ times};

\end{tikzpicture}
\quad
\begin{tikzpicture}[xscale=0.8,yscale=0.7,v/.style={pShapeVertex},every label/.append style={scale=0.8},baseline=0]

  \pic (C1) at (0.4,0) [scale=0.8] {comb};
  \pic (C2) at (3.5,3.5) [rotate=-90, scale=0.8] {comb};
  \pic (C3) at (6.8,0) [xscale=-1, scale=0.8] {comb};

  \node[v] (M) at (C1r3-|C2r3) {-2};       
  \node[v,below of=M] (T) {-1};

  \node[below=1cm of M,scale=0.85] {One shape for each $n\!>\!m\!>\!r\!\ge\!0$ with $n\!+\!m\!+\!r=g\!-\!7$};

  \draw (M)--(T);
  \draw (M)--(C1r3);
  \draw (M)--(C2r3);
  \draw (M)--(C3r3);

  \draw[decorate,decoration={brace,mirror},transform canvas={yshift=-3pt}]
    ($(C1r1)+(-0.2,-0.2)$)--($(C1r3)+(0.2,-0.2)$)
    node[midway,below=2pt,scale=0.9] {$n$ times};

  \draw[decorate,decoration={brace,mirror},transform canvas={yshift=-3pt}]
    ($(C3r3)+(-0.2,-0.2)$)--($(C3r1)+(0.2,-0.2)$)
    node[midway,below=2pt,scale=0.9] {$m$ times};

  \draw[decorate,decoration={brace,mirror},transform canvas={xshift=-3pt}]
    ($(C2r1)+(-0.2,0.2)$)--($(C2r3)+(-0.2,-0.2)$)
    node[midway,left=6pt,scale=0.9] {$r$ times};

\end{tikzpicture}
\end{theorem}

\begin{proof}
When $g\le 6$ the claims follow from Tables \ref{tabR2}-\ref{tabR6}. Suppose $S$ is a shape of genus $g\ge 7$.

Recall that $S$ is a graph whose vertex set $V$ is the set of principal components with $\chi_\G<0$, marked
with their $\chi_\G$. The edges correspond to multisets of inner chains, and let us split them into multiple
edges, one for each inner chain. This results in a multigraph $G$ with vertex set $V$ and edge set $E$ which
is in bijection with inner chains between principal components with $\chi<0$. Pick an ordering of the edges
from every vertex.

Suppose $N(S)\ge 0.512\cdot 10^g$. A family of shape $S$ is constructed by allocating a principal type $\G$ to every vertex $v$
and matching its edges to the edges incident to $v$. If $v$ is marked by $\chi_v$ and has $k_v$ incident edges,
the number of choices, by \ref{princountmain}, is at most $\sqrt{10}^{2-\chi_v-k_v}$. So
\begin{equation}
\label{famshapebound}
  N(S) \>\>\le\>\>  \prod_{v\in V} \sqrt{10}^{\,2-\chi_v-k_v}
    = \sqrt{10}^{\,2|V|-\sum_v{\chi_v}-\sum_v k_v}
    =\>\> 10^{|V|-\frac{2-2g}{2}-|E|} = 10^g 10^{|V|-|E|-1}.
\end{equation}
As $|V|-|E|\le 1$ in a connected multi-graph, we get $N(S)\le 10^g$. This is very close to the lower bound, so
all the inequalities must be close to equalities for $N(S)\ge 0.512\cdot 10^g$ to hold.
Let $X=\prod_{v\in V}(\cdots)$ in \eqref{famshapebound}.
The following conditions have to be met:

\begin{itemize}
\item
$G$ is a tree, for otherwise $|E|>|V|-1$ and $N(S)\le 10^{g-1}$ by \eqref{famshapebound}.
\item
In particular, $G$ is a simple graph, so $S$ and $G$ are identical.
\item
Let $(\chi_v, w_v)$ be the principal type and the weight vector of some vertex $v$ of $S$ (so $|w_v|=k_v$).
Then $(\chi_v, w_v)$ is one of the exceptional tuples from \ref{princountmain} that violate the stronger bound
$N^\chi_w \le \frac12\sqrt{10}^{2-\chi-|w|}$. Indeed, otherwise the contribution from $v$ to $X$
drops by a factor of $\ge 2$ resulting in $N(S)\le\frac 12 10^g$.
\item
Moreover, $(\chi_v, w_v)=(2,(2))$ is is impossible: the vertex $v'$ on the other side of the
edge from $v$ must have a matching weight $2\in w_{v'}$, and again $X$ drops by $(10^{3/2}/18)^2>2$.
\item
In particular, all edges have weight 1, and all leaves $v$ of $S$ have $\chi_v=-1$.
\end{itemize}
Under these conditions, there are only two reasons for $N(S)$ to drop below $10^g$.
\begin{enumerate}
\item
By \ref{princountmain},
every vertex of degree 2 in $S$ contributes a factor of 8/10
to $X/10^g$. Since $(8/10)^3=0.512$, so there are $\le 3$ such vertices.
\item
If $\Aut S\ne \{1\}$, then the count may drop.
\end{enumerate}
To quantify (2), consider the image $A$ of $\Aut S$ in
$\Aut(\text{leaves of $S$})$. If $A$ has an element of order $>2$ then
there is an orbit of length $n\ge 3$ on the leaves, and
$$
  X/10^g<\binom{10+n-1}n/10^n<1/2.
$$
Therefore $A\iso C_2^m$ for some $m\ge 0$. If $A$ has
an orbit of length $\ge 4$ or two orbits of length 2 with different stabilisers,
again we get $X/10^g<1/2$. It follows that $A$ is trivial or $A\iso C_2$.
In the latter case, $X/10^g$ drops at least by $55/100$.

Putting (1) and (2) together, we find that either
\begin{itemize}
\item[(A)]
$\Aut S=C_1$ and $S$ has at most 3 vertices of degree 2.
\item[(B)]
$\Aut S\iso C_2$ and $S$ has no vertices of degree 2.
\end{itemize}

It remains to classify connected trees $S$ that satisfy (A) or (B). Let $S_0\subset S$ be the subgraph
obtained by removing the leaves of $S$. As $g>2$, it is non-empty, and it is a connected tree.
Every leaf of $S_0$ is a either a vertex of $S$ of degree $2$ or has $\ge 2$ leaves of $S$ attached to it,
in which case swapping those leaves contributes at least $C_2$ to $\Aut S$. It follows that $S_0$ has at most
one leaf in case (B), which is impossible in genus $>3$, and at most three in case (A). So $S_0$ is one of

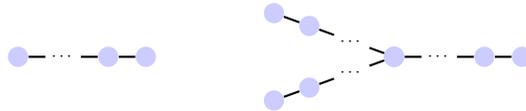
\begin{center}
  \tikzset{dottedchain/.pic = {
    \node[v] (a0) at (0,0) {};
    \node[scale=0.6] (dots) at (0.6,0) {$\cdots$};
    \node[v] (a1) at (1.2,0) {};
    \node[v] (a2) at (1.7,0) {};
    \draw[thick] (a0) -- (dots) -- (a1) -- (a2);
  }}
\begin{tikzpicture}[xscale=0.5,yscale=0.5,auto=left, v/.style={circle,scale=0.9,fill=blue!20,inner sep=3pt}]
  \pic at (0,0) {dottedchain};
  \pic at (10,0) {dottedchain};
  \pic at (10,0) [rotate=160] {dottedchain};
  \pic at (10,0) [rotate=200] {dottedchain};
\end{tikzpicture}
\end{center}

\noindent
Reattaching the leaves of $S$ gives the shapes claimed in the statement of the theorem. Finally,
it is elementary to check that

\begin{center}
\begin{tabular}{@{}l@{\quad}l@{ }l}
Number of solutions to $n\!>\!m\ge 0$ with $n\!+\!m=g\!-\!5$ & $=\lfloor\frac{g\!-\!4}{2}\rfloor$ &for $g\ge 4$, \\[1pt]
Number of solutions to $n\!>\!m\!>\!r\ge 0$ with $n\!+\!m\!+\!r=g\!-\!6$ & $=\lfloor\frac{(g\!-\!6)^2+4}{12}\rfloor$ &for $g\ge 4$, \\[2pt]
Number of solutions to $n\!>\!m\!>\!r\ge 0$ with $n\!+\!m\!+\!r=g\!-\!7$ & $=\lfloor\frac{(g\!-\!7)^2+4}{12}\rfloor$ &for $g\ge 5$, \\
\end{tabular}
\end{center}

\noindent
and that three counts add up to $\smash{\lceil\frac{(g-5)^2}{6}\rceil}$ for $g\ge 6$.
\end{proof}

\begin{theorem}
\label{thm10g}
The number of reduction families in genus $g$ is at least $10^g$.
\end{theorem}

\begin{proof}
When $g=1$, there are infinitely many reduction families when $g=1$. (In fact, elliptic curves
are enough, as there are 10 Kodaira families.)

When $2\le g\le 6$, the claim holds by \ref{numfamilies}.

When $g=7$, the following three shapes add up to $>10^7$ families:
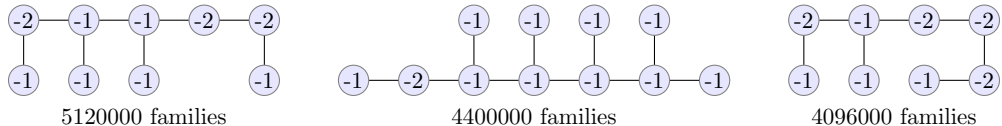
\begin{center}
\cbox{\begin{tikzpicture}[xscale=0.8,yscale=0.8,auto=left, v/.style={pShapeVertex}, l/.style={pShapeEdge}] \node[v] (1) at (5,1) {-1}; \node[v] (2) at (5,2) {-2}; \node[v] (3) at (4,2) {-2}; \node[v] (4) at (3,2) {-1}; \node[v] (5) at (2,2) {-1}; \node[v] (6) at (1,2) {-2}; \node[v] (7) at (1,1) {-1}; \node[v] (8) at (3,1) {-1}; \node[v] (9) at (2,1) {-1}; \draw[] (1)--node[l] {} (2); \draw[] (2)--node[l] {} (3); \draw[] (3)--node[l] {} (4); \draw[] (4)--node[l] {} (5); \draw[] (4)--node[l] {} (8); \draw[] (5)--node[l] {} (6); \draw[] (5)--node[l] {} (9); \draw[] (6)--node[l] {} (7); \node[scale=0.8] at (3,0.4) {5120000 families}; \end{tikzpicture}}\qquad\cbox{\begin{tikzpicture}[xscale=0.8,yscale=0.8,auto=left, v/.style={pShapeVertex}, l/.style={pShapeEdge}] \node[v] (1) at (1,1) {-1}; \node[v] (2) at (2,1) {-2}; \node[v] (3) at (3,1) {-1}; \node[v] (4) at (4,1) {-1}; \node[v] (5) at (5,1) {-1}; \node[v] (6) at (6,1) {-1}; \node[v] (7) at (7,1) {-1}; \node[v] (8) at (3,2) {-1}; \node[v] (9) at (4,2) {-1}; \node[v] (10) at (5,2) {-1}; \node[v] (11) at (6,2) {-1}; \draw[] (1)--node[l] {} (2); \draw[] (2)--node[l] {} (3); \draw[] (3)--node[l] {} (4); \draw[] (3)--node[l] {} (8); \draw[] (4)--node[l] {} (5); \draw[] (4)--node[l] {} (9); \draw[] (5)--node[l] {} (6); \draw[] (5)--node[l] {} (10); \draw[] (6)--node[l] {} (7); \draw[] (11)--node[l] {} (6); \node[scale=0.8] at (4,0.4) {4400000 families}; \end{tikzpicture}}\qquad\cbox{\begin{tikzpicture}[xscale=0.8,yscale=0.8,auto=left, v/.style={pShapeVertex}, l/.style={pShapeEdge}] \node[v] (1) at (3,1) {-1}; \node[v] (2) at (4,1) {-2}; \node[v] (3) at (4,2) {-2}; \node[v] (4) at (3,2) {-2}; \node[v] (5) at (2,2) {-1}; \node[v] (6) at (1,2) {-2}; \node[v] (7) at (1,1) {-1}; \node[v] (8) at (2,1) {-1}; \draw[] (1)--node[l] {} (2); \draw[] (2)--node[l] {} (3); \draw[] (3)--node[l] {} (4); \draw[] (4)--node[l] {} (5); \draw[] (5)--node[l] {} (6); \draw[] (5)--node[l] {} (8); \draw[] (6)--node[l] {} (7); \node[scale=0.8] at (2.5,0.4) {4096000 families}; \end{tikzpicture}}
\end{center}

When $g\ge 8$, by \ref{fixedfamcount} there are at least two shapes with $0.512\cdot 10^g$ families in each.
\end{proof}

\begin{remark}
The proof shows that shapes that are trees of Kodaira types with vertices
$$
\begin{tikzpicture}[baseline=0] \node[pPShapeChi] (1) at (0,0) {-1}; \draw[shorten >= 4pt] (1) -- ++(90:0.55) node[pPShapeWeight] {1}; \end{tikzpicture}
\quad
\text{(10 types)}
\qquad
\qquad
\begin{tikzpicture}[baseline=0] \node[pPShapeChi] (1) at (0,0) {-2}; \draw[shorten >= 4pt] (1) -- ++(110:0.55) node[pPShapeWeight] {1}; \draw[shorten >= 4pt] (1) -- ++(70:0.55) node[pPShapeWeight] {1}; \end{tikzpicture}
\quad
\text{(8 types)}
\qquad
\qquad
\begin{tikzpicture}[baseline=0] \node[pPShapeChi] (1) at (0,0) {-3}; \draw[shorten >= 4pt] (1) -- ++(110:0.55) node[pPShapeWeight] {1}; \draw[shorten >= 4pt] (1) -- ++(90:0.55) node[pPShapeWeight] {1}; \draw[shorten >= 4pt] (1) -- ++(70:0.55) node[pPShapeWeight] {1}; \end{tikzpicture}
\quad
\text{(1 type)}
$$
(cf.~Tables \ref{tabPS1}, \ref{tabPS2})
are enough to get $10^g$ families in every genus $g\ge 5$. They all have principal components
of genus 0 and multiplicity $m\le 6$, so it is a very restricted collection of families.
However, they do not exhibit the true growth of the number of families with $g$, which is superexponential.
In fact, any collection of shapes where $\chi$ of every vertex is bounded from below has a number of families
which is at most exponential in $g$, and the true number of shapes is superexponential.
To show this, we focus on shapes that correspond to semistable snc families, as it turns out that
these are enough to get the $g^{g-o(g)}$ growth.
\end{remark}

\begin{definition}
\label{defsemistable}
A reduction type $R$ is \textbf{semistable} if every component of $R$ has multiplicity~1,
and its reduction family is called a \textbf{semistable family}\footnote{
Recall that a curve has \emph{semistable reduction} (or simply is \emph{semistable}) if it has a regular
model (equivalently, mrnc model) in which every component has multiplicity 1.
Every reduction type in a semistable family is semistable (and, similarly, for snc families).}.
\end{definition}

\begin{definition}
\label{defsnc}
A reduction type $R$ has \textbf{strict normal crossings (snc)} if no component of $R$ has loops,
and its reduction family is called an \textbf{snc family}.
\end{definition}

See tables \ref{sem2table} and \ref{sem3table} for all semistable families in genus 2 and 3.

\begin{remark}
If a shape contains a semistable snc family, such a family is unique: the only principal types
are ${\rm I}_{{\rm g}g}\langle\text{edges}\rangle$, and such a type is determined by its edges and its Euler
characteristic, so a vertex in the shape determines the only principal type that can be placed at such
a vertex. In other words, there are maps

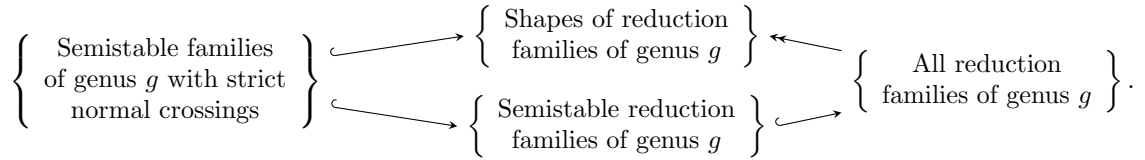
\begin{center}
\begin{tikzpicture}[>=stealth]
  \node[scale=0.9] (A) at (-1,0) {$\left\{\begin{tabular}{c}Semistable families\\of genus $g$ with strict\\normal crossings\end{tabular}\right\}$};
  \node[scale=0.9] (B) at (5,0.6) {$\left\{\begin{tabular}{c}Shapes of reduction\\families of genus $g$\end{tabular}\right\}$};
  \node[scale=0.9] (C) at (5,-0.6) {$\left\{\begin{tabular}{c}Semistable reduction\\families of genus $g$\end{tabular}\right\}$};
  \node[scale=0.9] (D) at (10,0) {$\left\{\begin{tabular}{c}All reduction\\families of genus $g$\end{tabular}\right\}.$};
  \draw[right hook->] (A) -- (B.west);
  \draw[right hook->] (A) -- (C.west);
  \draw[->>] (D) -- (B.east);
  \draw[right hook->] (C.east) -- (D);
\end{tikzpicture}
\end{center}
\noindent
The top left arrow is injective, and the top right one has fibres of size $\le 0.512\cdot 10^g$ by \ref{fixedfamcount}.
\end{remark}

\begin{theorem}
\label{thmgg}
The number of reduction families in genus $g$ is at least $g^{g-\varepsilon(g)}$
for some function $\varepsilon(g)$ with $\varepsilon(g)\to0$ as $g\to \infty$. The same holds
for the number of shapes, the number of semistable families and the number of semistable snc families.
\end{theorem}

\begin{lemma}[F. Petrov]
\label{petrov}
Fix $1/2<\lambda<1$. For large enough $n$, the number of simple graphs
on $k=\bigl\lceil 2 n^{\lambda}\bigr\rceil$ vertices with $n$ edges up to isomorphism is at least
$$
  n^{(2\lambda-1)n-k}.
$$
\end{lemma}

\begin{proof}
First, consider simple graphs on the $k$ \emph{labelled} vertices with $n$ edges.
Since $\binom{k}{2}=\frac{k(k-1)}{2}>n^{2\lambda}+n$ for large $n$,
the number of such graphs is
$$
  \binom{\binom{k}{2}}{n}
    =\frac{\prod_{j=0}^{n-1}\bigl(\binom{k}{2}-j\bigr)}{n!}
    >\frac{(n^{2\lambda})^{n}}{n!}
    >\frac{n^{2\lambda n}}{n^n}  = n^{(2\lambda-1)n}.
$$
Each unlabelled $k$-vertex graph arises from at most $k!$ labelings, so the number of distinct unlabelled graphs on
$k$ vertices with $n$ edges for $n$ large enough is at least
$$
  \frac{n^{(2\lambda-1)n}}{k!}\ge \frac{n^{(2\lambda-1)n}}{n^k}= n^{(2\lambda-1)n-k}.
$$
\end{proof}

\begin{proof}[Proof of Theorem \ref{thmgg}]
Let $g$ be large enough, and $1\le k<n<\binom{k}{2}$ be integers with $k+n+1<g$.
Write $g_0=g-n-k\ge 2$. Start with a semistable snc family that has one component of genus $g_0$ as a `root' and $k$ genus 1 components as `leaves', each connected to the root with a chain of $\P^1$s of multiplicity 1. Its shape is
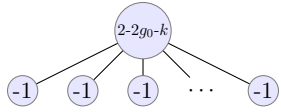
\begin{center}
\begin{tikzpicture}[xscale=0.8,yscale=0.8,auto=left, v/.style={pShapeVertex}, l/.style={pShapeEdge}]
  \node[v,scale=0.7] (1) at (3,2) {$2\text{-}2g_0\text{-}k$};
  \node[v] (2) at (1,1) {-1};
  \node[v] (3) at (2,1) {-1};
  \node[v] (4) at (3,1) {-1};
  \node[scale=0.8] (5) at (4,1) {$\cdots$};
  \node[v] (6) at (5,1) {-1};
  \draw[] (1)--node[l] {} (2);
  \draw[] (1)--node[l] {} (3);
  \draw[] (1)--node[l] {} (4);
  \draw[] (1)--node[l] {} (5);
  \draw[] (1)--node[l] {} (6);
\end{tikzpicture}
\end{center}
Next, add $n$ chains of $\P^1$s between some of the leaves,
at most one for each unordered pair.
The total genus of such a family is $k+g_0+n=g$ by \ref{corabtor}, so we have a map
\begin{center}
  \{simple graphs on $n$ vertices with $k$ edges\} $\longinjects$ \{semistable snc families of genus $g$\},
\end{center}
adding a root to a given graph on the leaves. The root is the only vertex of genus $>1$, so the map is injective.

Now we specify $n$ and $k$ in terms of $g$ and give a lower bound on the left-hand side.
Set
$$
  n = \Bigl\lfloor g - \frac{g}{\log g} \Bigr\rfloor, \qquad
  k = \bigl\lceil 2 n^\lambda \bigr\rceil \quad\text{with}\quad
  \lambda = 1-\frac{3\log\log n}{\log n}.
$$
Then, for large $g$,
$$
  n^\lambda = n^{1 - 3\log\log n / \log n} = \frac{n}{(\log n)^3}
\le \frac{g}{(\log g)^3},
$$
and
$$
  k \le 3 n^\lambda \le \frac{3g}{(\log g)^3}, \qquad
  n+k+1 \le \left(g - \frac{g}{\log g}\right) + \frac{3g}{(\log g)^3} +1 < g.
$$
With these $n$ and $k$ as functions of $g$, write
\begin{center}
  $a_g = $ number of simple connected graphs on $k$ vertices with $n$ edges.
\end{center}
By Lemma~\ref{petrov}, it is bounded from below by
$$
  a_g \ge n^{(2\lambda-1)n - k}.
$$
Hence
$$
  \frac{\log a_g}{g\log g} \ge \frac{((2\lambda-1)n - k)\log n}{g\log g}
  = \frac{(2\lambda-1)n \log n}{g\log g} - \frac{k\log n}{g\log g}.
$$
The first term satisfies
$$
  \frac{(2\lambda-1)n \log n}{g\log g} = \frac{n}{g}\frac{\log n}{\log g}\left(1 - \frac{6\log\log n}{\log n}\right) = 1 - o(1),
$$
since $n = g - g/\log g$, so $n/g = 1 + o(1)$ and $\log n/\log g = 1 + o(1)$.
The second term is $o(1)$ since $k=o(g)$ and $\log n\le \log g$. We get
$$
  \frac{\log a_g}{g\log g} \ge 1 - o(1),
$$
and $a_g \ge g^{\,g - o(g)}$, as required.
\end{proof}

\clearpage

\newgeometry{top=1.5cm, bottom=1cm, left=1cm, right=1cm}

\section{Tables}
\label{stables}

\begin{center}
$$
{\rm I}=1^\emptyset \qquad
{\rm D}=2^{1,1} \qquad
{\rm T}=3^{1,2} \qquad
4^{1,3} \qquad
5^{1,4} \qquad
5^{2,3} \qquad
6^{1,5} \qquad
7^{1,6} \qquad
7^{2,5} \qquad \cdots
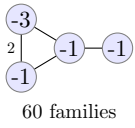
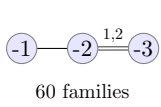
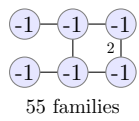
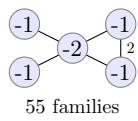
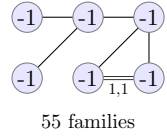
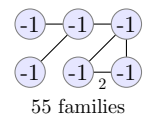
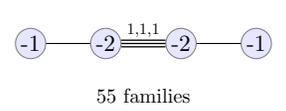
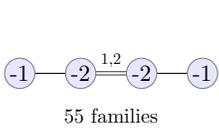
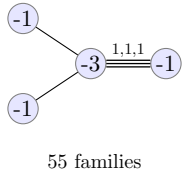
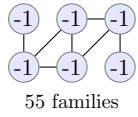
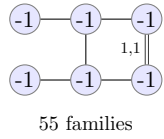
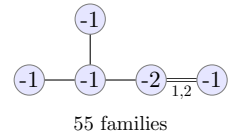
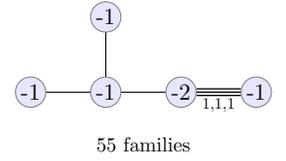
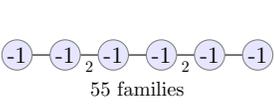
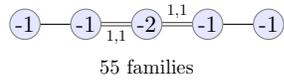
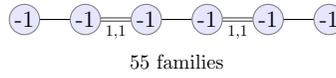
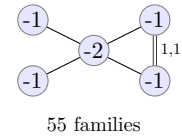
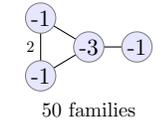
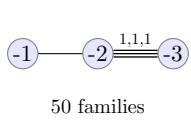
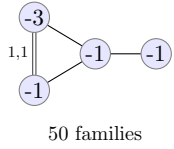
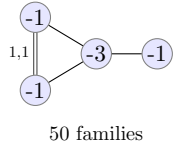
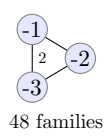
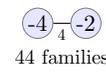
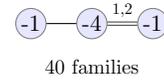
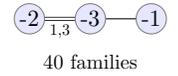
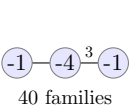
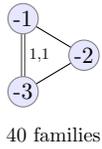
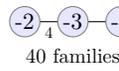
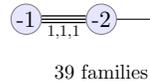
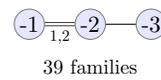
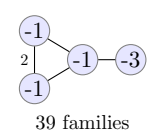
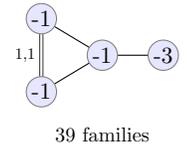
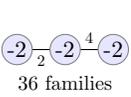
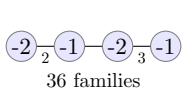
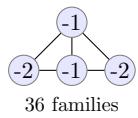
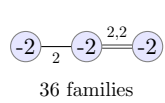
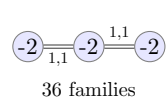
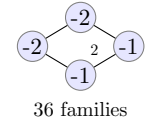
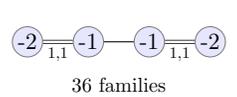
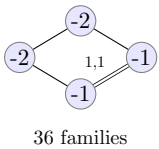
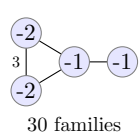
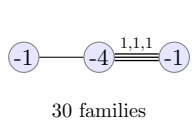
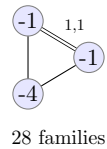
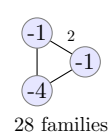
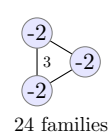
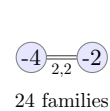
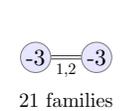
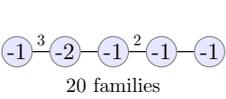
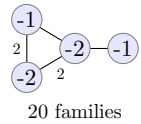
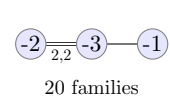
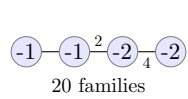
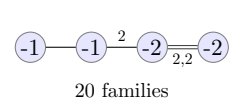
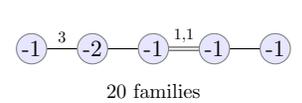
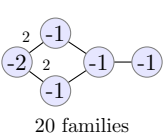
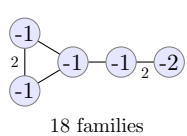
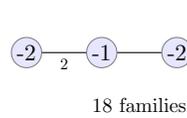
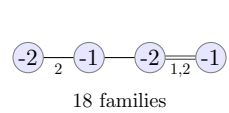
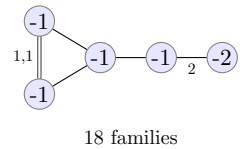
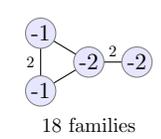
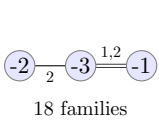
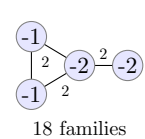
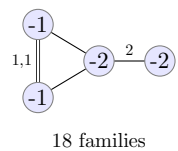
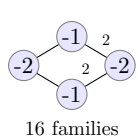
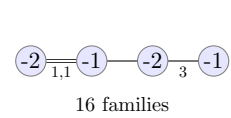
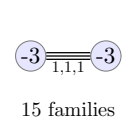
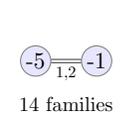
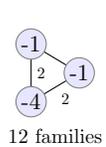
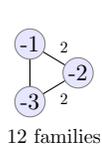
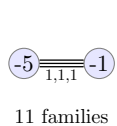
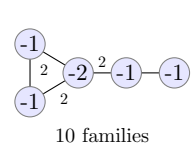
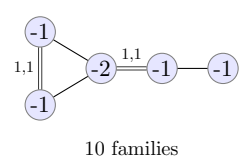
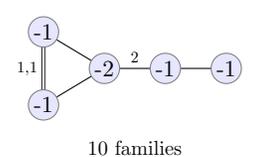
$$
\textsc{Table} \customlabel{tabC2}{C2}. (Infinitely many) cores with $\chi=2$
\end{center}

\vfillrule

\begin{center}
\\[4pt]
\textsc{Table} \customlabel{tabC0}{C0}. All cores with $\chi=0$
\end{center}

\vfillrule

\begin{center}
%
\\[4pt]
\textsc{Table} \customlabel{tabC-2}{C-2}. All cores with $\chi=-2$
\end{center}

\vfillrule

\begin{center}
%
\\[4pt]
\textsc{Table} \customlabel{tabC-4}{C-4}. All cores with $\chi=-4$
\end{center}

\vfillrule

\begin{center}
%
\\[4pt]
\textsc{Table} \customlabel{tabC-6}{C-6}. All cores with $\chi=-6$
\end{center}

\vfillrule

\begin{center}
%
\\[4pt]
\textsc{Table} \customlabel{tabC-8}{C-8}. All cores with $\chi=-8$
\end{center}

\newpage

\begin{center}
%

\vfill

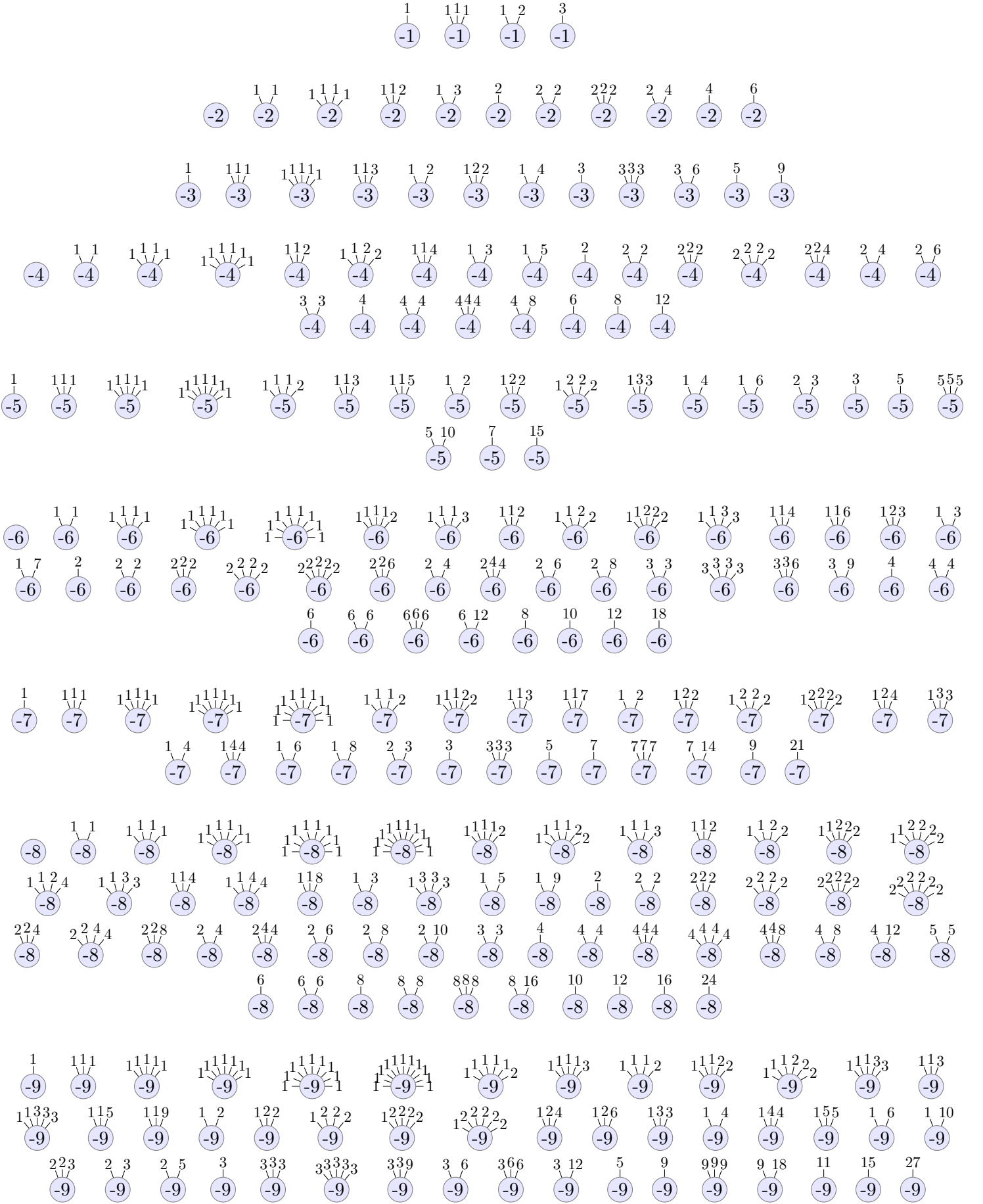
\textsc{Table} \customlabel{tabPS}{PS}. Possible shapes $(\chi_T, \text{weight})$ of principal types $T$ with $-9\le\chi_T\le -1$.
\end{center}

\newpage

\begin{center}
\begin{tabular}{|p{0.75cm}@{\qquad}p{16.5cm}|}\hline Shape & Types\vphantom{$\int^X$}\\[-6pt] {\null \begin{tikzpicture}[baseline=0] \node[pPShapeChi] (1) at (0,0) {-1}; \draw[shorten >= 4pt] (1) -- ++(90:0.55) node[pPShapeWeight] {1}; \end{tikzpicture}} & {\null \begin{tikzpicture}[scale=0.6] \path[use as bounding box] (-1,-0.09) rectangle (26.5,2); \draw(0,0) \shapearclabel{35pt}{1}{}{1g1}{\redtype{I_{g1}\e{1}}} ++(2.75,0) \shapearclabel{35pt}{1-1,1}{}{1}{\redtype{I_1\e{1}}} ++(2.75,0) \shapearclabel{35pt}{1}{1,1,1}{2}{\redtype{I^*_0\e{1}}} ++(2.75,0) \shapearclabel{35pt}{1,.2}{1}{2}{\redtype{I^*_1\e{1}}} ++(2.75,0) \shapearclabel{35pt}{1}{1,1}{3}{\redtype{IV\e{1}}} ++(2.75,0) \shapearclabel{35pt}{2}{2,2}{3}{\redtype{IV^*\e{2}}} ++(2.75,0) \shapearclabel{35pt}{1}{1,2}{4}{\redtype{III\e{1}}} ++(2.75,0) \shapearclabel{35pt}{3}{3,2}{4}{\redtype{III^*\e{3}}} ++(2.75,0) \shapearclabel{35pt}{1}{2,3}{6}{\redtype{II\e{1}}} ++(2.75,0) \shapearclabel{35pt}{5}{4,3}{6}{\redtype{II^*\e{5}}}; \end{tikzpicture} }\\ \hline\end{tabular}\par\vskip-1pt\noindent\begin{tabular}{|p{0.75cm}@{\qquad}p{4.244cm}|p{0.75cm}@{\qquad}p{4.244cm}|p{0.75cm}@{\qquad}p{4.244cm}|}Shape & Types & Shape & Types & Shape & Types\vphantom{$\int^X$}\\[-6pt] {\null \begin{tikzpicture}[baseline=0] \node[pPShapeChi] (1) at (0,0) {-1}; \draw[shorten >= 4pt] (1) -- ++(110:0.55) node[pPShapeWeight] {1}; \draw[shorten >= 4pt] (1) -- ++(90:0.55) node[pPShapeWeight] {1}; \draw[shorten >= 4pt] (1) -- ++(70:0.55) node[pPShapeWeight] {1}; \end{tikzpicture}} & {\null \begin{tikzpicture}[scale=0.6] \path[use as bounding box] (-1,-0.09) rectangle (1.75,2); \draw(0,0) \shapearclabel{35pt}{1,1,1}{}{1}{\redtype{I\e{1}\e{1}\e{1}}}; \end{tikzpicture} }& {\null \begin{tikzpicture}[baseline=0] \node[pPShapeChi] (1) at (0,0) {-1}; \draw[shorten >= 4pt] (1) -- ++(110:0.55) node[pPShapeWeight] {1}; \draw[shorten >= 4pt] (1) -- ++(70:0.55) node[pPShapeWeight] {2}; \end{tikzpicture}} & {\null \begin{tikzpicture}[scale=0.6] \path[use as bounding box] (-1,-0.09) rectangle (1.75,2); \draw(0,0) \shapearclabel{35pt}{1,2}{1}{2}{\redtype{D\e{1}\e{2}}}; \end{tikzpicture} }& {\null \begin{tikzpicture}[baseline=0] \node[pPShapeChi] (1) at (0,0) {-1}; \draw[shorten >= 4pt] (1) -- ++(90:0.55) node[pPShapeWeight] {3}; \end{tikzpicture}} & {\null \begin{tikzpicture}[scale=0.6] \path[use as bounding box] (-1,-0.09) rectangle (1.75,2); \draw(0,0) \shapearclabel{35pt}{3}{1,2}{3}{\redtype{T\e{3}}}; \end{tikzpicture} }\\ \hline\end{tabular}
\par\bigskip
\textsc{Table} \customlabel{tabPS1}{PS1} of
13
principal types with $\chi=-1$, grouped by shape.
\end{center}

\vfill

\begin{center}
\begin{tabular}{|p{0.75cm}@{\qquad}p{16.5cm}|}\hline Shape & Types\vphantom{$\int^X$}\\
\hline\end{tabular}%
\par\bigskip
\textsc{Table} \customlabel{tabPS2}{PS2} of
83
principal types with $\chi=-2$, grouped by shape.
\end{center}

\vfill
\newpage

\begin{center}
\begin{tabular}{|p{0.75cm}@{\qquad}p{16.5cm}|}\hline Shape & Types\vphantom{$\int^X$}\\
\hline\end{tabular}%
\par\bigskip
\textsc{Table} \customlabel{tabPS3}{PS3} of
75
principal types with $\chi=-3$, grouped by shape.
\end{center}

\vfill
\newpage

\begin{center}
\renewcommand{\arraystretch}{0}

\\[-1.4pt]
\begin{tabular}{|p{6cm}|p{6cm}|p{6cm}|}
\hline
Shape\>\>\>\>Types & Shape\>\>\>\>Types & Shape\>\>\>\>Types \Tcr
\begin{tikzpicture}[xscale=0.95,yscale=1,] \node at (0.15,-0.25) {}; \node[pPShapeChi,scale=0.8] (1) at (0.5,-0.25) {-4}; \draw[shorten >= 4pt] (1) -- ++(110:0.41) node[pPShapeWeight,scale=0.72] {2}; \draw[shorten >= 4pt] (1) -- ++(90:0.41) node[pPShapeWeight,scale=0.72] {2}; \draw[shorten >= 4pt] (1) -- ++(70:0.41) node[pPShapeWeight,scale=0.72] {2}; \node at (2,0) {\redtype{D\e{2}\e{2}\e{2}}}; \node at (4.1,0) {\redtype{[2]I_D\e{2}\e{2}\e{2}}}; \end{tikzpicture}&
\begin{tikzpicture}[xscale=0.95,yscale=1,] \node at (0.15,-0.25) {}; \node[pPShapeChi,scale=0.8] (1) at (0.5,-0.25) {-4}; \draw[shorten >= 4pt] (1) -- ++(110:0.41) node[pPShapeWeight,scale=0.72] {2}; \draw[shorten >= 4pt] (1) -- ++(70:0.41) node[pPShapeWeight,scale=0.72] {6}; \node at (2,0) {\redtype{[2]T\e{4}\e{6}}}; \node at (4.1,0) {\redtype{[2]T\e{2}\e{6}}}; \end{tikzpicture}&
\begin{tikzpicture}[xscale=0.95,yscale=1,] \node at (0.15,-0.25) {}; \node[pPShapeChi,scale=0.8] (1) at (0.5,-0.25) {-4}; \draw[shorten >= 4pt] (1) -- ++(110:0.41) node[pPShapeWeight,scale=0.72] {4}; \draw[shorten >= 4pt] (1) -- ++(70:0.41) node[pPShapeWeight,scale=0.72] {4}; \node at (2,0) {\redtype{[2]D\e{4}\e{4}}}; \node at (4.1,0) {\redtype{[4]I_D\e{4}\e{4}}}; \end{tikzpicture}\\
\hline
\end{tabular}
\par\bigskip
\textsc{Table} \customlabel{tabPS4}{PS4} of
277
principal types with $\chi=-4$, grouped by shape.
\end{center}

\vfill
\newpage

\begin{center}
\renewcommand{\arraystretch}{0}
\begin{tabular}{|p{18.75cm}|}
\hline
Shape\>\>\>\>Types\Tcr
\begin{tikzpicture}[xscale=0.95,yscale=1,] \node at (0.15,-0.25) {}; \node[pPShapeChi,scale=0.8] (1) at (0.5,-0.25) {-5}; \draw[shorten >= 4pt] (1) -- ++(90:0.41) node[pPShapeWeight,scale=0.72] {1}; \node at (2,0) {\redtype{I_{g3}\e{1}}}; \node at (4.79,0) {\redtype{I_{1,g2}\e{1}}}; \node at (7.58,0) {\redtype{I_{1,g1,1}\e{1}}}; \node at (10.37,0) {\redtype{I_{1,1,1}\e{1}}}; \node at (13.16,0) {\redtype{I^*_{0,g1}\e{1}}}; \node at (15.95,0) {\redtype{I^*_{1,g1}\e{1}}}; \node[scale=0.9] at (18.74,0) {\redtype{2^{1,1,1,1,1,1,1,1}\e{1}}}; \node at (2,-0.6) {\redtype{2^{1,1,1,1,1,1}_0\e{1}}}; \node at (4.79,-0.6) {\redtype{I^*_{0,\{2\m2\}}\e{1}}}; \node at (7.58,-0.6) {\redtype{2^{1,1,1,1,1,1}_D\e{1}}}; \node at (10.37,-0.6) {\redtype{I^*_{0,0,D}\e{1}}}; \node at (13.16,-0.6) {\redtype{I^*_{1,\{2\m2\}}\e{1}}}; \node at (15.95,-0.6) {\redtype{I^*_{0,D,D}\e{1}}}; \node at (18.74,-0.6) {\redtype{I^*_{1,D,D}\e{1}}}; \node at (2,-1.2) {\redtype{T_{g1}\e{2}}}; \node at (4.79,-1.2) {\redtype{T_{g1}\e{1}}}; \node at (7.58,-1.2) {\redtype{3^{1,1,1,1,2}\e{2}}}; \node at (10.37,-1.2) {\redtype{3^{1,1,1,1,2}\e{1}}}; \node at (13.16,-1.2) {\redtype{3^{1,2,2,2,2}\e{2}}}; \node at (15.95,-1.2) {\redtype{3^{1,2,2,2,2}\e{1}}}; \node at (18.74,-1.2) {\redtype{3^{1,1,2,2}_{\{1\m2\}}\e{2}}}; \node at (2,-1.8) {\redtype{3^{1,1,2,2}_{\{2\m2\}}\e{1}}}; \node at (4.79,-1.8) {\redtype{T_{\{3\m3\}}\e{2}}}; \node at (7.58,-1.8) {\redtype{3^{1,1,2,2}_0\e{2}}}; \node at (10.37,-1.8) {\redtype{3^{1,1,2,2}_{\{1\m2\}}\e{1}}}; \node at (13.16,-1.8) {\redtype{T_{\{3\m3\}}\e{1}}}; \node at (15.95,-1.8) {\redtype{4^{1,1,1,1}\e{1}}}; \node at (18.74,-1.8) {\redtype{4^{1,1,3,3}\e{3}}}; \node at (2,-2.4) {\redtype{4^{1,1,3,3}\e{1}}}; \node at (4.79,-2.4) {\redtype{4^{3,3,3,3}\e{3}}}; \node at (7.58,-2.4) {\redtype{4^{1,1,2,2,2}\e{1}}}; \node at (10.37,-2.4) {\redtype{4^{3,3,2,2,2}\e{3}}}; \node at (13.16,-2.4) {\redtype{4^{1,3,2,2}_D\e{3}}}; \node at (15.95,-2.4) {\redtype{III_{\{4\}D}\e{1}}}; \node at (18.74,-2.4) {\redtype{4^{1,3,2,2}_D\e{1}}}; \node at (2,-3) {\redtype{III^*_{\{4\}D}\e{3}}}; \node at (4.79,-3) {\redtype{5^{1,1,3}_0\e{3}}}; \node at (7.58,-3) {\redtype{5^{1,2,2}_0\e{2}}}; \node at (10.37,-3) {\redtype{5^{1,1,3}_{\{1\m3\}}\e{1}}}; \node at (13.16,-3) {\redtype{5^{1,2,2}_{\{2\m2\}}\e{1}}}; \node at (15.95,-3) {\redtype{5^{2,4,4}_{-1}\e{4}}}; \node at (18.74,-3) {\redtype{5^{3,3,4}_0\e{4}}}; \node at (2,-3.6) {\redtype{5^{3,3,4}_{\{3\m4\}}\e{3}}}; \node at (4.79,-3.6) {\redtype{5^{2,4,4}_{\{4\m4\}}\e{2}}}; \node at (7.58,-3.6) {\redtype{6^{1,5,3,3}\e{5}}}; \node at (10.37,-3.6) {\redtype{6^{5,2,2,3}\e{5}}}; \node at (13.16,-3.6) {\redtype{6^{1,4,4,3}\e{1}}}; \node at (15.95,-3.6) {\redtype{6^{1,5,3,3}\e{1}}}; \node at (18.74,-3.6) {\redtype{6^{1,1,4}_D\e{1}}}; \node at (2,-4.2) {\redtype{6^{1,5}_D\e{5}}}; \node at (4.79,-4.2) {\redtype{6^{5,5,2}_D\e{5}}}; \node at (7.58,-4.2) {\redtype{6^{1,5}_D\e{1}}}; \node at (10.37,-4.2) {\redtype{7^{1,1,5}\e{5}}}; \node at (13.16,-4.2) {\redtype{7^{1,2,4}\e{4}}}; \node at (15.95,-4.2) {\redtype{7^{1,3,3}\e{3}}}; \node at (18.74,-4.2) {\redtype{7^{1,2,4}\e{2}}}; \node at (2,-4.8) {\redtype{7^{1,1,5}\e{1}}}; \node at (4.79,-4.8) {\redtype{7^{2,2,3}\e{3}}}; \node at (7.58,-4.8) {\redtype{7^{2,2,3}\e{2}}}; \node at (10.37,-4.8) {\redtype{7^{1,2,4}\e{1}}}; \node at (13.16,-4.8) {\redtype{7^{2,6,6}\e{6}}}; \node at (15.95,-4.8) {\redtype{7^{1,3,3}\e{1}}}; \node at (18.74,-4.8) {\redtype{7^{3,5,6}\e{6}}}; \node at (2,-5.4) {\redtype{7^{3,5,6}\e{5}}}; \node at (4.79,-5.4) {\redtype{7^{4,4,6}\e{6}}}; \node at (7.58,-5.4) {\redtype{7^{4,5,5}\e{5}}}; \node at (10.37,-5.4) {\redtype{7^{4,4,6}\e{4}}}; \node at (13.16,-5.4) {\redtype{7^{4,5,5}\e{4}}}; \node at (15.95,-5.4) {\redtype{7^{3,5,6}\e{3}}}; \node at (18.74,-5.4) {\redtype{7^{2,6,6}\e{2}}}; \node at (2,-6) {\redtype{8^{1,5,2}\e{5}}}; \node at (4.79,-6) {\redtype{8^{1,1,6}\e{1}}}; \node at (7.58,-6) {\redtype{8^{3,3,2}\e{3}}}; \node at (10.37,-6) {\redtype{8^{3,7,6}\e{7}}}; \node at (13.16,-6) {\redtype{8^{1,5,2}\e{1}}}; \node at (15.95,-6) {\redtype{8^{5,5,6}\e{5}}}; \node at (18.74,-6) {\redtype{8^{7,7,2}\e{7}}}; \node at (2,-6.6) {\redtype{8^{3,7,6}\e{3}}}; \node at (4.79,-6.6) {\redtype{9^{1,5,3}\e{5}}}; \node at (7.58,-6.6) {\redtype{9^{1,2,6}\e{2}}}; \node at (10.37,-6.6) {\redtype{9^{2,4,3}\e{4}}}; \node at (13.16,-6.6) {\redtype{9^{1,2,6}\e{1}}}; \node at (15.95,-6.6) {\redtype{9^{2,4,3}\e{2}}}; \node at (18.74,-6.6) {\redtype{9^{4,8,6}\e{8}}}; \node at (2,-7.2) {\redtype{9^{1,5,3}\e{1}}}; \node at (4.79,-7.2) {\redtype{9^{5,7,6}\e{7}}}; \node at (7.58,-7.2) {\redtype{9^{7,8,3}\e{8}}}; \node at (10.37,-7.2) {\redtype{9^{5,7,6}\e{5}}}; \node at (13.16,-7.2) {\redtype{9^{7,8,3}\e{7}}}; \node at (15.95,-7.2) {\redtype{9^{4,8,6}\e{4}}}; \node at (18.74,-7.2) {\redtype{10^{3,2,5}_D\e{3}}}; \node at (2,-7.8) {\redtype{10^{1,4,5}_D\e{1}}}; \node at (4.79,-7.8) {\redtype{10^{9,6,5}_D\e{9}}}; \node at (7.58,-7.8) {\redtype{10^{7,8,5}_D\e{7}}}; \node at (10.37,-7.8) {\redtype{12^{1,5,6}\e{5}}}; \node at (13.16,-7.8) {\redtype{12^{5,3,4}\e{5}}}; \node at (15.95,-7.8) {\redtype{12^{1,3,8}\e{1}}}; \node at (18.74,-7.8) {\redtype{12^{1,5,6}\e{1}}}; \node at (2,-8.4) {\redtype{12^{7,11,6}\e{11}}}; \node at (4.79,-8.4) {\redtype{12^{11,9,4}\e{11}}}; \node at (7.58,-8.4) {\redtype{12^{7,9,8}\e{7}}}; \node at (10.37,-8.4) {\redtype{12^{7,11,6}\e{7}}}; \node at (13.16,-8.4) {\redtype{14^{5,2,7}\e{5}}}; \node at (15.95,-8.4) {\redtype{14^{3,4,7}\e{3}}}; \node at (18.74,-8.4) {\redtype{14^{1,6,7}\e{1}}}; \node at (2,-9) {\redtype{14^{13,8,7}\e{13}}}; \node at (4.79,-9) {\redtype{14^{11,10,7}\e{11}}}; \node at (7.58,-9) {\redtype{14^{9,12,7}\e{9}}}; \end{tikzpicture}\\
\end{tabular}
\\[-1.4pt]
\begin{tabular}{|p{6cm}|p{6cm}|p{6cm}|}
\hline
Shape\>\>\>\>Types & Shape\>\>\>\>Types & Shape\>\>\>\>Types \Tcr
\begin{tikzpicture}[xscale=0.95,yscale=1,] \node at (0.15,-0.25) {}; \node[pPShapeChi,scale=0.8] (1) at (0.5,-0.25) {-5}; \draw[shorten >= 4pt] (1) -- ++(130:0.41) node[pPShapeWeight,scale=0.72] {1}; \draw[shorten >= 4pt] (1) -- ++(110:0.41) node[pPShapeWeight,scale=0.72] {1}; \draw[shorten >= 4pt] (1) -- ++(90:0.41) node[pPShapeWeight,scale=0.72] {1}; \draw[shorten >= 4pt] (1) -- ++(70:0.41) node[pPShapeWeight,scale=0.72] {1}; \draw[shorten >= 4pt] (1) -- ++(50:0.41) node[pPShapeWeight,scale=0.72] {1}; \node at (2.5,0) {\redtype{I_{g1}\e{1}\e{1}\e{1}\e{1}\e{1}}}; \node at (5.1,0) {\redtype{I_1\e{1}\e{1}\e{1}\e{1}\e{1}}}; \end{tikzpicture}&
\begin{tikzpicture}[xscale=0.95,yscale=1,] \node at (0.15,-0.25) {}; \node[pPShapeChi,scale=0.8] (1) at (0.5,-0.25) {-5}; \draw[shorten >= 4pt] (1) -- ++(110:0.41) node[pPShapeWeight,scale=0.72] {1}; \draw[shorten >= 4pt] (1) -- ++(90:0.41) node[pPShapeWeight,scale=0.72] {1}; \draw[shorten >= 4pt] (1) -- ++(70:0.41) node[pPShapeWeight,scale=0.72] {3}; \node at (2,0) {\redtype{IV\e{1}\e{1}\e{3}}}; \node at (4.1,0) {\redtype{IV^*\e{2}\e{2}\e{3}}}; \end{tikzpicture}&
\begin{tikzpicture}[xscale=0.95,yscale=1,] \node at (0.15,-0.25) {}; \node[pPShapeChi,scale=0.8] (1) at (0.5,-0.25) {-5}; \draw[shorten >= 4pt] (1) -- ++(110:0.41) node[pPShapeWeight,scale=0.72] {1}; \draw[shorten >= 4pt] (1) -- ++(90:0.41) node[pPShapeWeight,scale=0.72] {1}; \draw[shorten >= 4pt] (1) -- ++(70:0.41) node[pPShapeWeight,scale=0.72] {5}; \node at (2,0) {\redtype{5^{1,4}\e{1}\e{4}\e{5}}}; \node at (4.1,0) {\redtype{5^{2,3}\e{2}\e{3}\e{5}}}; \end{tikzpicture}\\
\hline
\end{tabular}
\par\bigskip
\textsc{Table} \customlabel{tabPS5}{PS5} of
176
principal types with $\chi=-5$, grouped by shape.
\end{center}

\vfill
\newpage

\begin{center}
\renewcommand{\arraystretch}{0}
%
\par\bigskip
\textsc{Table} \customlabel{tabPS6}{PS6} of
282
principal types with $\chi=-6$ and empty weights.
\end{center}

\vfill 

$$\begin{array}{|l|cccccccccc|}
\hline
\chi_T & -1 & -2 & -3 & -4 & -5 & -6 & -7 & -8 & -9 & -10\cr
\#\text{types} & 13 & 83 & 75 & 277 & 176 & 591 & 352 & 1068 & 684 & 1686\cr
\text{Table} & \text{\ref{tabPS1}} & \text{\ref{tabPS2}} & \text{\ref{tabPS3}} & \text{\ref{tabPS4}} & \text{\ref{tabPS5}} &&&&& \cr
\hline
\chi_T & -11 & -12 & -13 & -14 & -15 & -16 & -17 & -18 & -19 & -20\cr
\#\text{types} & 1162 & 3324 & 1880 & 4371 & 3400 & 6992 & 4677 & 11267 & 8189 & 17034\cr\hline
\chi_T & -21 & -22 & -23 & -24 & -25 & -26 & -27 & -28 & -29 & -30\cr
\#\text{types} & 12094 & 22320 & 18921 & 38691 & 26710 & 49919 & 44239 & 75846 & 57236 & 108656\cr\hline
\end{array}$$

\bigskip

\begin{center}
\textsc{Table} \customlabel{tabPN}{PN}.
The number of principal types with given $\chi=-1,-2,...,-30$.\\All reduction types of genus $g\ge 2$ are made out of
principal types with $2-2g\le \chi_T\le -1$.
\end{center}

\newpage

\begin{minipage}[c]{\textwidth}\begin{center}        
\cbox{$\redtype{I_{g2}}$}\quad \cbox{$\redtype{I_{1,g1}}$}\quad \cbox{$\redtype{I_{1,1}}$}\quad \cbox{$%
 $}
\vskip 8pt
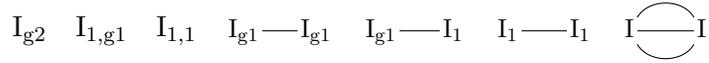
\textsc{Table} \customlabel{sem2table}{S2} of 7 semistable type families in genus 2
\end{center}\end{minipage}

\bigskip
\vfillrule

\begin{minipage}[c]{\textwidth}\begin{center}        
\cbox{$\redtype{I_{g3}}$}\quad \cbox{$\redtype{I_{1,g2}}$}\quad \cbox{$\redtype{I_{1,g1,1}}$}\quad \cbox{$\redtype{I_{1,1,1}}$}\quad \cbox{$%
 $}
\vskip 12pt
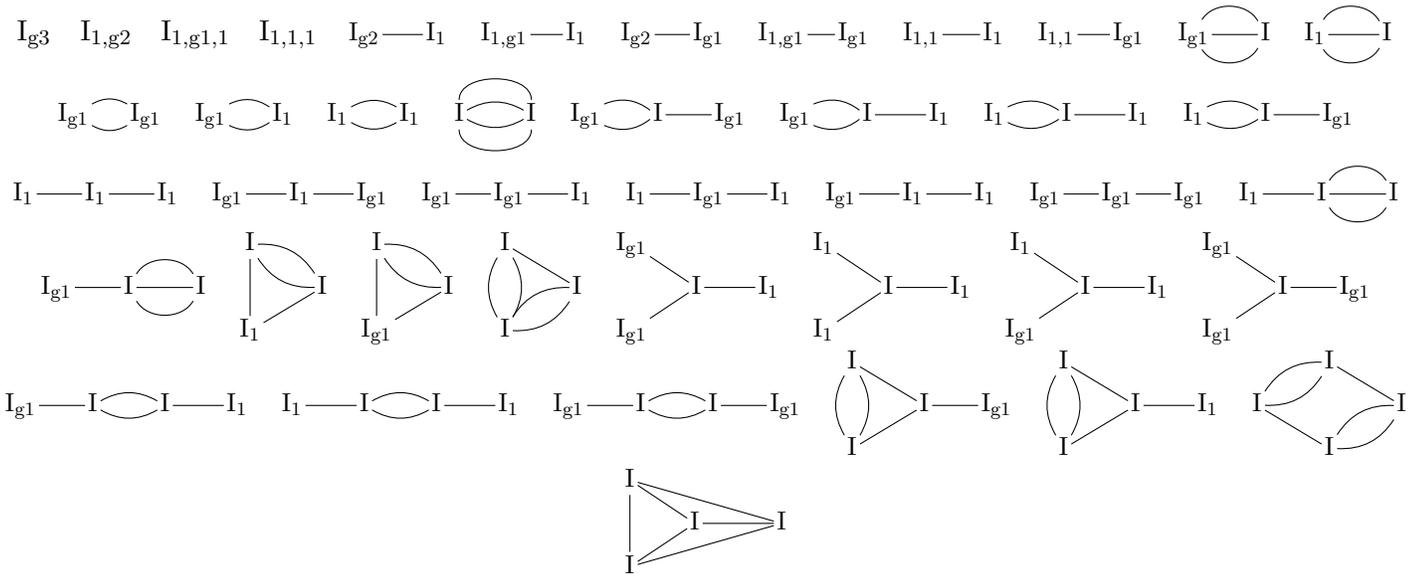
\textsc{Table} \customlabel{sem3table}{S3} of 42 semistable type families in genus 3
\end{center}\end{minipage}

\bigskip
\vfillrule

\begin{center}   
\cbox{%
}

\bigskip
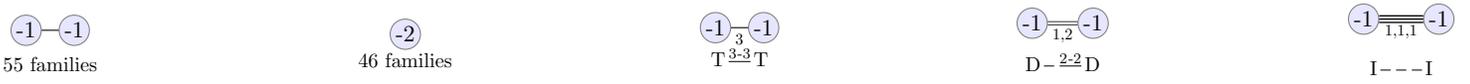
\textsc{Table} \customlabel{tabR2}{R2}. Shapes for all \gtwonumtypes{} reduction families in genus 2
\end{center}

\vfillrule

\begin{center}   
\cbox{%
}

\bigskip
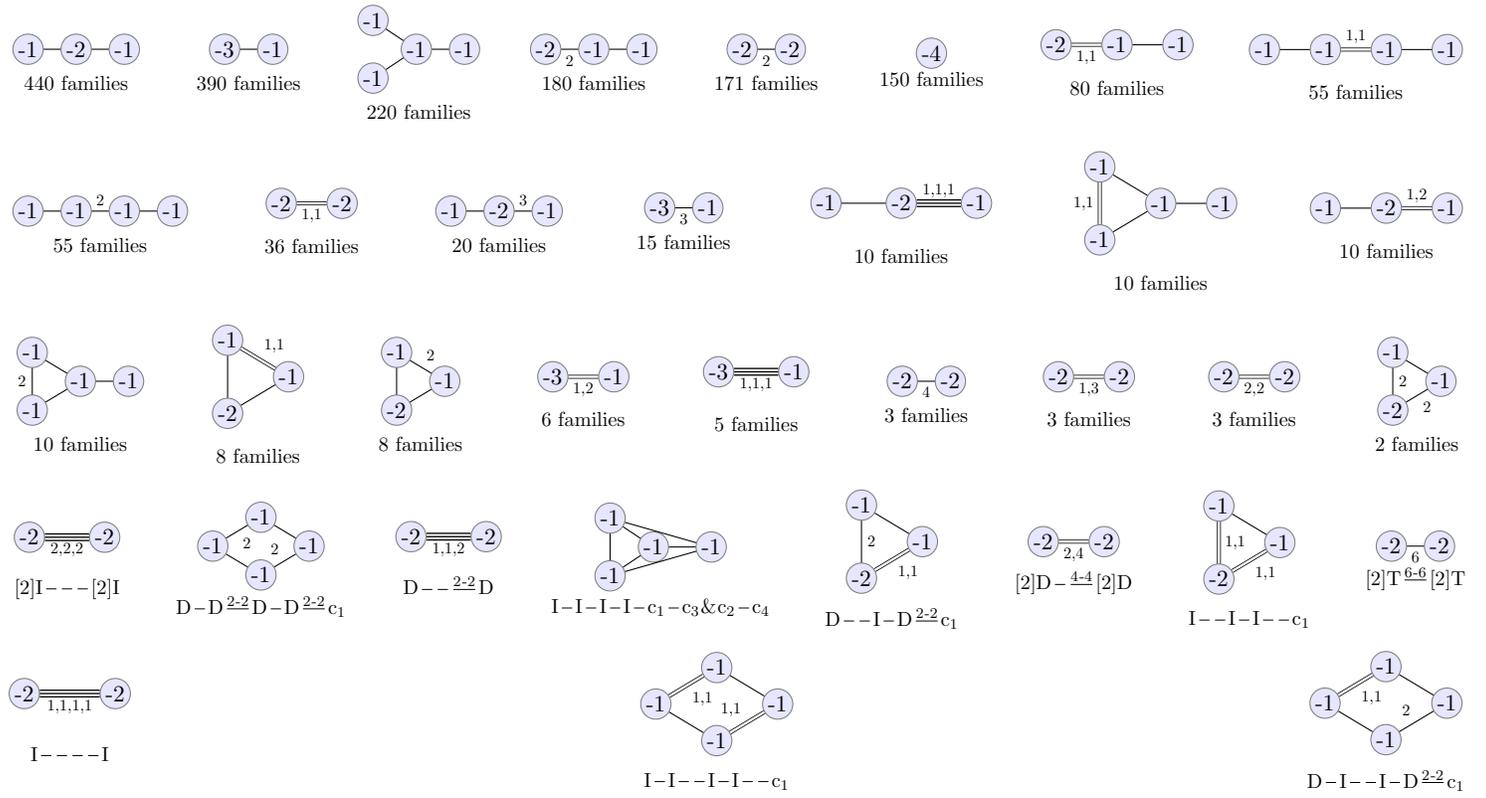
\textsc{Table} \customlabel{tabR3}{R3}. Shapes for all \gthreenumtypes{} reduction families in genus 3
\vfill
\end{center}

\vfill
\clearpage

\begin{center}   
\cbox{%
}

\bigskip
\textsc{Table} \customlabel{tabR4}{R4}. Shapes for all \gfournumtypes{} reduction families in genus 4
\vfill
\end{center}

\begin{center}   
\cbox{%
}

\vfill
\textsc{Table} \customlabel{tabR5}{R5}. Shapes with most reduction families in genus 5
\vfill\vfill
\end{center}

\vfill
\clearpage

\begin{center}    
\cbox{%
}

\bigskip
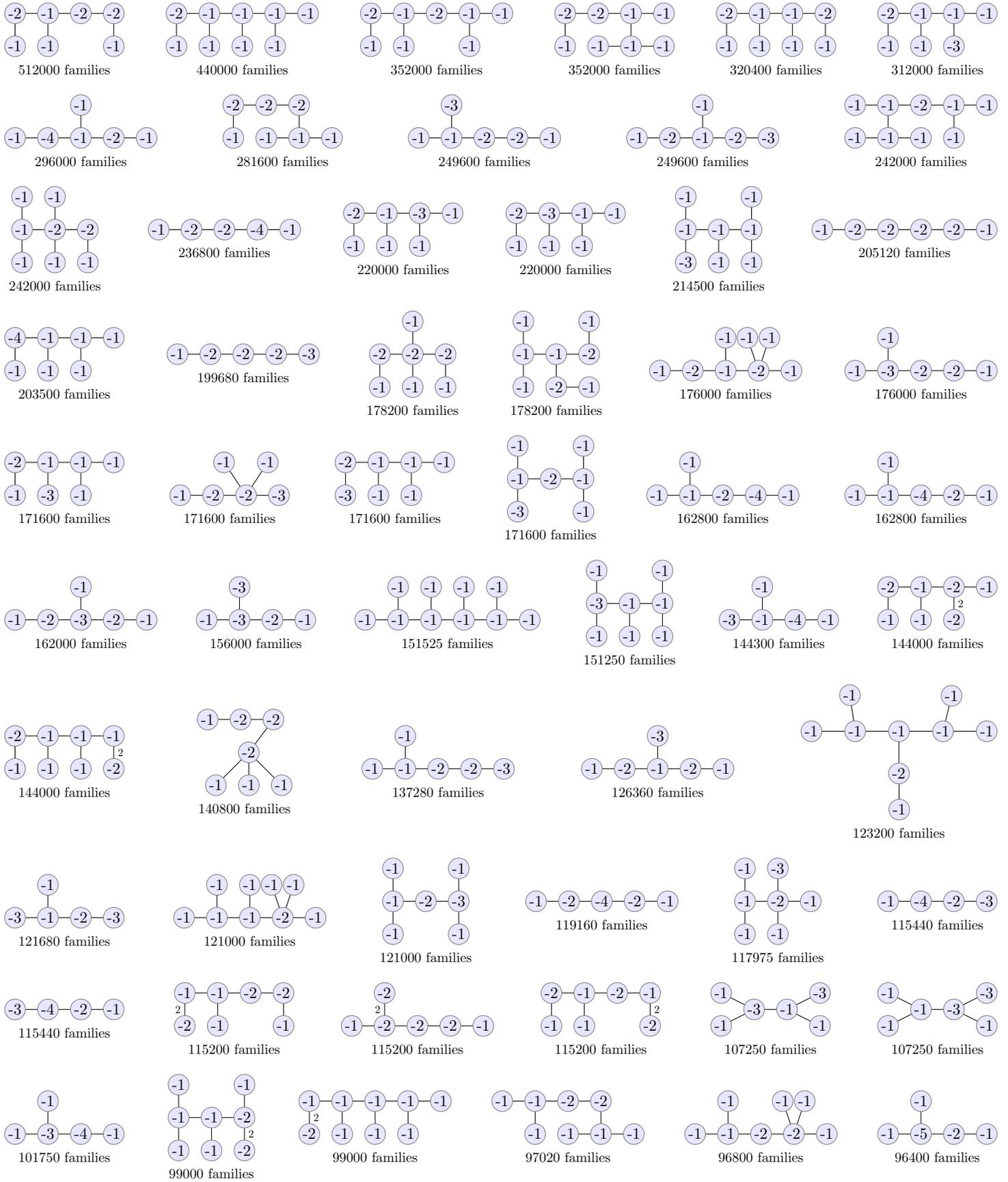
\textsc{Table} \customlabel{tabR6}{R6}. Shapes with most reduction families in genus 6
\vfill\vfill
\end{center}

\vfill
\clearpage
$$
  \begingroup
  \tikzstyle{smallmult}=[blue,scale=0.77,inner sep=1.5pt,left]
  \tikzstyle{prinmult}=[blue,scale=0.8,inner sep=2pt,above left]
  \tikzstyle{gamma}=[inner sep=2pt,below left=0pt and -1pt,scale=0.8]
  \tikzstyle{quotheader}=[scale=1.1,align=center]
  \tikzstyle{unitheader}=[scale=0.75]
  \tikzstyle{restfibre}=[below left=-3pt and -2.2em,scale=0.7]
  \begin{array}{l@{\>\>}l@{\>\>}l@{\qquad}c@{\qquad}l@{\qquad}l}
  \multicolumn{3}{c}{\begin{tabular}{c}Quotient $\frac{m_\G}{-\chi_\G}>2$,\\$\G$ principal component\end{tabular}} &
  \begin{tabular}{c}$\mgOLG=(nm_{\G_0},0,n\cO_{\G_0},n\cL_{\G_0})$\cr for some $n\ge 1$ and $\G_0$ as follows:
  \end{tabular} &
  \multicolumn{2}{c}{\text{parameters}}\\[10pt]\hline\\[-8pt]
   2+\frac{4}{2k-1}&=&6,3\frac{1}{3},2\frac{4}{5},... & \begin{tikzpicture}[xscale=1.75,yscale=1,baseline=0] \draw[thick] (0,0)--(2.7,0) node[prinmult] {$4k\text{+}2$} node[gamma] {$\G_0$}; \draw[shorten <=-3pt] (0.4,0)--node[smallmult] {$2k\text{+}1\,\text{-}\,u$} (0.4,-0.6); \draw[shorten <=-3pt] (1.2,0)--node[smallmult] {$2k\text{+}1$} (1.2,-0.6); \draw[shorten <=-3pt] (2,0)--node[smallmult] {$u$} (2,0.6); \end{tikzpicture} & u\in\bigl(\frac{\mathbb Z}{(4k\text{+}2)\mathbb Z}\bigr)^\times & k\ge 1\\[17pt] 2+\frac{3}{k-1}&=&5,3\frac{1}{2},3,... & \begin{tikzpicture}[xscale=1.75,yscale=1,baseline=0] \draw[thick] (0,0)--(2.7,0) node[prinmult] {$4k\text{+}2$} node[gamma] {$\G_0$}; \draw[shorten <=-3pt] (0.4,0)--node[smallmult] {$u$} (0.4,-0.6); \draw[shorten <=-3pt] (1.2,0)--node[smallmult] {$2k\text{+}1\,\text{-}\,u$} (1.2,-0.6); \draw[shorten <=-3pt] (2,0)--node[smallmult] {$2k\text{+}1$} (2,-0.6); \end{tikzpicture} & u\in\bigl(\frac{\mathbb Z}{(4k\text{+}2)\mathbb Z}\bigr)^\times & k\ge 2\\[17pt] 2+\frac{2}{k-1}&=&4,3,2\frac{2}{3},... & \begin{tikzpicture}[xscale=1.75,yscale=1,baseline=0] \draw[thick] (0,0)--(2.7,0) node[prinmult] {$4k$} node[gamma] {$\G_0$}; \draw[shorten <=-3pt] (0.4,0)--node[smallmult] {$u$} (0.4,-0.6); \draw[shorten <=-3pt] (1.2,0)--node[smallmult] {$2k\,\text{-}\,u$} (1.2,-0.6); \draw[shorten <=-3pt] (2,0)--node[smallmult] {$2k$} (2,-0.6); \end{tikzpicture} & u\in\bigl(\frac{\mathbb Z}{4k\mathbb Z}\bigr)^\times & k\ge 2\\
\end{array}
  \endgroup$$
  \par\bigskip
  \begin{center}
  \textsc{Table} \customlabel{tablargem}{M}.
  All possible principal components $\G$ of large multiplicity $\frac{m_\G}{-\chi_\G}>2$ in genus $>1$.\\[4pt]
  They are multiples $[n]\G_0$ ($n\ge 1$) with $\G_0$ as in the table.
  Outer multiplicities $o\in\cO_{\G_0}$ are shown pointing down, inner
  multiplicities $l\in\cL_{\G_0}$ pointing up; they are viewed as elements of $\Z/m_{\G_0}\Z$.
  \end{center}

\clearpage

\input{namuenofull.inc}  
\vskip 12pt
\begin{center}
\textsc{Table} \customlabel{tabNU}{G2}. Correspondence with Namikawa-Ueno-Liu classification in genus 2
\label{NUtable}
\end{center}

\restoregeometry

\end{document}